\documentclass{elsarticle}

\usepackage{amssymb}
\usepackage{amsmath,subfigure}
\usepackage{amsthm}
\usepackage{color}

\usepackage[left=1in,right=1in,top=0.82in,bottom=0.82in]{geometry}

\newtheorem{theorem}{Theorem}[section]
\newtheorem{lemma}{Lemma}[section]
\newtheorem{proposition}{Proposition}[section]

\newtheorem{example}{Example}

\def\hE{\mathcal{E}}

\newcommand{\bra}[1]{\left(#1\right)}
\newcommand{\brab}[1]{\big(#1\big)}

\newcommand{\brat}[1]{(#1)}
\newcommand{\kbra}[1]{\left[#1\right]}
\newcommand{\kbrab}[1]{\big[#1\big]}
\newcommand{\kbraB}[1]{\Big[#1\Big]}
\newcommand{\kbrat}[1]{[#1]}

\newcommand{\myinnerb}[1]{\big\langle#1\big\rangle}
\newcommand{\myinnerB}[1]{\Big\langle#1\Big\rangle}

\newcommand{\mynorm}[1]{\left\|#1\right\|}

\newcommand{\mynormt}[1]{\|#1\|}

\newcommand{\cte}{{\mathrm{T}}}

\newcommand{\cnt}{{\mathrm{N}}}

\def\qing#1{\textcolor{cyan}{#1}}

\def\lan#1{\textcolor{blue}{#1}}

\begin{document}

\begin{frontmatter}




\title{Original energy dissipation preserving corrections of integrating factor Runge-Kutta methods for gradient flow problems}


\author[1,2]{Hong-lin Liao\corref{cor1}}
\ead{liaohl@nuaa.edu.cn}\ead{liaohl@csrc.ac.cn}

\author[1]{Xuping Wang}
\ead{wangxp@nuaa.edu.cn}

\author[1]{Cao Wen}
\ead{0131122255@mail.imu.edu.cn}

\address[1]{School of Mathematics, Nanjing University of Aeronautics and Astronautics, Nanjing 211106, P. R. China}
\address[2]{Key Laboratory of Mathematical Modeling and High Performance Computing of Air Vehicles (NUAA), MIIT, Nanjing 211106, P. R. China}

\cortext[cor1]{Corresponding author}



\begin{abstract}
Explicit integrating factor Runge-Kutta methods are attractive and popular in developing high-order maximum bound principle preserving time-stepping schemes for Allen-Cahn type gradient flows. However, they always suffer from the non-preservation of steady-state solution and original energy dissipation law. To overcome these disadvantages, some new integrating factor methods are developed by using two classes of difference correction, including the telescopic correction and nonlinear-term translation correction, enforcing the preservation of steady-state solution. Then the original energy dissipation properties of the new methods are examined by using the associated differential forms and the differentiation matrices. As applications, some new integrating factor Runge-Kutta methods up to third-order maintaining the original energy dissipation law are constructed by applying the difference correction strategies to some popular explicit integrating factor methods in the literature. Extensive numerical experiments are presented to support our theory and to demonstrate the improved performance of new methods.
\end{abstract}




\begin{keyword}
gradient flow problem \sep integrating factor Runge-Kutta method \sep steady-state preserving correction \sep differentiation matrix \sep discrete energy dissipation law


\MSC 35K58 \sep 65L20 \sep 65M06 \sep 65M12
\end{keyword}

\end{frontmatter}



\section{Introduction}
We develop some correction schemes of explicit integrating factor Runge-Kutta (in short, IF) methods and examine their energy dissipation properties for the semi-discrete semilinear parabolic problem
\begin{align}\label{problem: autonomous}
  u_{h}^{\prime}(t)+L_{h} u_{h}(t)=g(u_{h}(t)), \quad u_{h}(t_{0})=u_{h}^{0},
\end{align}
where $L_{h}$ is a symmetric, positive definite matrix resulting from certain spatial discretization of stiff term, typically the Laplacian operator $-\Delta$ with periodic boundary conditions, and $g$ represents a nonlinear but non-stiff term. Assume that there exists a non-negative Lyapunov function $G$ such that $g(v)=-\frac{\delta}{\delta v} G(v)$. Then the semi-discrete semilinear problem \eqref{problem: autonomous} can be formulated into a gradient flow system
\begin{align}\label{eq: continuous energy}
     \frac{\mathrm{d} u_{h}}{\mathrm{d} t}=-\frac{\delta E}{\delta u_{h}} \quad \text { with } \quad E[u_{h}]:=\frac{1}{2} \myinnerb{u_{h}, L_{h} u_{h}} + \myinnerb{G(u_{h}), 1}.
\end{align}
Without losing the generality, finite difference method is assumed here to approximate spatial operators and we define the discrete $L^{2}$ inner product $\langle u_{h}, v_{h}\rangle:=v_{h}^{T} u_{h}$ and the $L^{2}$ norm $\|v_{h}\|:=\sqrt{\langle v_{h}, v_{h}\rangle}$. Also, $\|v_{h}\|_{\infty}$ represents the maximum norm on the discrete spatial mesh. 
The dynamics approaching the steady-state solution $u_{h}^{*}$, that is $L_{h} u_{h}^{*}=g(u_{h}^{*})$, of this gradient flow system \eqref{eq: continuous energy} satisfies the energy dissipation law
\begin{align}\label{eq: continuous energy law}
  \frac{\mathrm{d} E}{\mathrm{d} t}= \myinnerB{\frac{\delta E}{\delta u_{h}}, \frac{\mathrm{d} u_{h}}{\mathrm{d} t}} =- \myinnerB{\frac{\mathrm{d} u_{h}}{\mathrm{d} t}, \frac{\mathrm{d} u_{h}}{\mathrm{d} t}} \leq 0.
\end{align}

Let $u_{h}^{k}$ be the numerical approximation of $u_{h}(t_{k})$ at the point $t_{k}$ for $0 \leq k \leq N$. To integrate the semilinear parabolic problem \eqref{problem: autonomous} from the discrete time $t_{n-1} \, (n \geq 1)$ to the next grid point $t_{n}=t_{n-1}+\tau$, classical IF methods start from an integrating factor form of the problem \eqref{problem: autonomous},
\begin{align*}
  w_{h}^{\prime}(t) = e^{t L_{h}} g\kbrab{e^{-t L_{h}}w_{h}(t)}=:G\brab{w_{h}(t)} \quad\text{with }w_{h}(t):=e^{t L_{h}} u_{h}(t).
\end{align*}
Let $W^{n, i}$ be the approximation of $w_{h}(t_{n-1}+c_{i} \tau)$ at the abscissas $c_{1}:=0, c_{i} \in(0,1]$ for $2 \leq i \leq s$, and $c_{s+1}:=1$. Replacing $\tau$ by $c_{i} \tau$ to define the internal stages $t_{n,i}:=t_{n-1}+c_{i} \tau$ \lan{ (typically, $\tau$ also represents a variable-step
	size)}, one has the general class of explicit Runge-Kutta methods:
\begin{align*}
   & W^{n,i+1}= W^{n,1}+\tau \sum_{j=1}^{i} a_{i+1, j}(0) G(W^{n, j}) \quad\text{for $1\le i\le s-1$,} \\
  & W^{n,s+1}= W^{n,1}+\tau \sum_{i=1}^{s} b_{i}(0) G(W^{n, i}), 
\end{align*}
where $W^{n,1} = w_{h}^{n-1}$ and $w_{h}^{n}= W^{n,s+1}$. 
Always, we define 
\begin{align*}
 a_{s+1, j}(0):=b_{j}(0) \quad \text{for $1 \leq j \leq s$},
\end{align*}
and assume that the coefficient $a_{k+1,k}(0)\not=0$ for $1\leq k \leq s$.
Then the associate Butcher tableau reads
\begin{align}\label{Butcher tableau}
\begin{array}{c|ccccc}
  c_{1} & 0 & & & &  \\
  c_{2} & a_{2,1} & 0 & & &  \\
  c_{3} & a_{3,1} & a_{3,2} & 0 & &  \\
  \vdots & \vdots & \vdots & \ddots & \ddots &  \\
  c_{s} & a_{s,1} & a_{s,2} & \cdots & a_{s, s-1} & 0  \\
  \hline & a_{s+1,1} & a_{s+1,2} & \cdots & a_{s+1, s-1} & a_{s+1, s}
\end{array}\,,\end{align}
where the coefficients $a_{i+1,j}:=a_{i+1,j}(0)$ for $1\le j\le i\le s$. 
Returning to the original variables $U^{n, i}\approx u_{h}(t_{n,i})$, one obtains the following general IF  method or the so-called Lawson method \cite{HochbruckLeiboldOstermann:2020NM,Lawson:1967SINUM} for \eqref{problem: autonomous},
\begin{align}\label{scheme: IF}
  & U^{n,i+1}= e^{-c_{i+1} \tau L_{h}} U^{n,1}+\tau \sum_{j=1}^{i} a_{i+1, j}(-\tau L_{h}) g(U^{n, j}) \quad\text{for $1\le i\le s$,} 
\end{align}
where $U^{n,1} = u_{h}^{n-1}$ and $u_{h}^{n}= U^{n,s+1}$. The method coefficients $a_{i j}$ are defined by, cf. \cite{LiLiJuFeng:2021SISC,MasetZennaro:2009MCOM,ZhangYanQianChenSong:2022JSC},
\begin{align}\label{def: IF coefficients}
  a_{i+1,j}(-\tau L_{h}):=a_{i+1,j}(0)e^{-(c_{i+1}-c_{j})\tau L_{h}}\quad\text{for $1\le j\le i\le s$.}
\end{align}
The IF methods \eqref{scheme: IF} do not have stiff order one, but they perform to their non-stiff convergence order on some particular problems where the solutions satisfy certain smoothness properties, cf. \cite{MasetZennaro:2009MCOM}.
Here and hereafter, the exponentials $e^{-(c_{i+1}-c_{j})\tau L_{h}}$ are the matrix functions defined on the spectrum of $-{\tau}L_h$, see \cite[Theorem 1.13]{Higham:2008} for more properties on the function of matrix exponential. Typically, $f(-\tau L_{h})$ is a positive definite operator if the given function $f(z)$ is entire and positive. 
\lan{In general, choosing the appropriate method to calculate matrix exponentials depends on the properties and application scenarios of the matrix, cf. \cite{FasiGaudreaultLundSchweitzer:2024,MolerVanLoan:2003SIREV}. Always, we apply the MATLAB function `expm', which uses the scaling and squaring method internally, to calculate the matrix exponentials in our experiments.}

The IF scheme \eqref{scheme: IF} reduces 
to the classical explicit Runge-Kutta method with coefficients $a_{ij}:=a_{ij}(0)$ if we put $L_h=0$. 
The latter method is called the underlying explicit Runge-Kutta method such that
\begin{align}\label{cond: equilibria underlying RK} 
  \sum_{j=1}^{i}a_{i+1,j}(0)=c_{i+1}\quad\text{for $ 1\le i\le s$.}
\end{align}
Always, the formal order of an IF method is the same as that of the underlying
Runge-Kutta method \cite{MasetZennaro:2009MCOM}. 

To control the possible nonlinear instability introduced by the explicit discretization of the nonlinear term, we follow \cite{DuJuLiQiao:2019SINUM,DuJuLiQiao:2021SIREV} to introduce a stabilized parameter $\kappa\ge0$ and the linear stabilized term $\kappa u$, 
\begin{align}\label{def: stabilized parameter}
  L_{\kappa}:=L_h+\kappa I\quad\text{and}\quad g_{\kappa}(u):=\kappa u+g(u),
\end{align} 
such that the problem \eqref{problem: autonomous} becomes the stabilized version
\begin{align}\label{problem: stabilized version}
  u_h'(t)=-L_{\kappa}u_h(t)+g_{\kappa}(u_h),\quad u_h(t_0)=u_h^0.
\end{align}
It also preserves the original steady-state solution $u_{h}^{*}$ of the gradient flow system \eqref{eq: continuous energy}.
Thus, applying the IF method \eqref{scheme: IF} to the stabilized problem \eqref{problem: stabilized version}, we have the following stabilized IF method
\begin{align}\label{scheme: IF stabilized}
  U^{n,i+1}= e^{-c_{i+1} \tau L_{\kappa}} U^{n,1}+\tau \sum_{j=1}^{i} a_{i+1, j}(-\tau L_{\kappa}) g_{\kappa}(U^{n, j}) \quad\text{for $ 1\le i\le s$}.
\end{align}
Obviously, the stabilized IF method does not change the formal order of the IF method \eqref{scheme: IF}.
 It is worth noting that, cf. \cite{LiLiJuFeng:2021SISC,ZhangYanQianSong:2021ANM,ZhangYanQianSong:2022CMAME}, an appropriate $\kappa>0$ is always necessary for the stabilized IF methods \eqref{scheme: IF stabilized} to maintain the maximum bound principle of  Allen-Cahn type models.

Recently, Ju et al. \cite{JuLiQiaoYang:2021JCP} viewed each stage of IF methods \eqref{scheme: IF} as a convex combination of the exponential forward Euler substeps and established the maximum bound principle under the time-step requirements having the same magnitudes as the first-order IF scheme. They identified some concrete and practical IF schemes up to fourth-order. The idea of using the Shu-Osher form was also applied in \cite{ZhangYanQianSong:2021ANM} to develop maximum bound principle preserving IF methods up to fourth-order with certain time-step constraints related to the strong stability preserving property.
To remove the time-step requirements, Li et al. \cite{LiLiJuFeng:2021SISC} combined the linear stabilized technique to develop the stabilized IF methods \eqref{scheme: IF stabilized}. Under some additional conditions on the abscissas and nonnegative Butcher coefficients, they proved that the stabilized IF methods maintain the maximum bound principle unconditionally. 
The imposed conditions on the nonnegative Butcher coefficients in \cite{LiLiJuFeng:2021SISC} seem somewhat severe and can not meet with many classic strong stability-preserving stabilized IF schemes except the first-order one. Moreover, they recognized the so-called exponential effect \cite{DuJuLiQiao:2021SIREV,DuYang:2019JCP}, the maximum bounds of the stabilized IF solutions decay exponentially when the time-step size is moderately large, and attributed it to the magnitude of stabilized parameter $\kappa$. Recently, two improved techniques to relieve the above two disadvantages were suggested by Zhang et al. \cite{ZhangYanQianSong:2022CMAME} and the resulting two modified families of stabilized IF schemes \eqref{scheme: IF stabilized} were also shown to preserve the maximum bound principle unconditionally. Nonetheless, it is noticed that none of the above mentioned works established the energy dissipation properties of the proposed IF schemes.

We note that, the coefficient conditions in \eqref{cond: equilibria underlying RK} ensure that the underlying explicit Runge-Kutta method of \eqref{scheme: IF} or \eqref{scheme: IF stabilized} preserves the equilibria $u_{h}^{*}$; but they can not ensure the IF methods \eqref{scheme: IF} and \eqref{scheme: IF stabilized} to preserve the equilibria $u_{h}^{*}$, cf. \cite{CoxMatthews:2002JCP,FornbergDriscoll:1999JCP,HuangShu:2018JCP}. Actually, we obtain from \eqref{def: IF coefficients} that
\begin{align}\label{cond: non-equilibria}
  \sum_{j=1}^{i}a_{i+1,j}(z)=\sum_{j=1}^{i}a_{i+1,j}(0)e^{(c_{i+1}-c_{j})z}\neq\frac{e^{c_{i+1}z}-1}{z},
\end{align}
for $z\neq0$ and $1\le i\le s.$ It says that the IF methods \eqref{scheme: IF} and \eqref{scheme: IF stabilized} always suffer from the non-preservation of steady-state solution. For example, consider the Ginzburg-Landau double-well potential $G(u)=\frac1{4}(u^2-1)^2$ such that $g(u)=u-u^3$. Putting $u_h^k=u_h^*$ in the first-order IF (IF1) method \eqref{scheme: IF1} yields
\begin{align*}
  u_h^*= e^{- \tau L_{\kappa}} u_h^*+\tau  e^{- \tau L_{\kappa}} \kbrab{(1+\kappa)u_h^*-(u_h^*)^3}.
\end{align*}
Obviously, only the metastable state solution $u_h^*=0$ solves this equation exactly, 
while the desired equilibria $u_h^*$ only approximately satisfies the steady-state equation with an error of $O(\kappa^{2}\tau^2)$.
Another typical example is the Flory-Huggins potential $G(u):=\frac{\theta}{2}[(1+u)\ln(1+u)+(1-u)\ln(1-u)]-\frac{\theta_c}{2}u^2$ for $0<\theta<\theta_c$ such that $g(u)=\frac{\theta}{2}\ln\frac{1-u}{1+u}+\theta_c u$. Putting $u_h^k=u_h^*$ in the IF1 method \eqref{scheme: IF1} yields 
\begin{align*}
  u_h^*= e^{- \tau L_{\kappa}} u_h^*+\tau  e^{- \tau L_{\kappa}} \kbrab{(\theta_c+\kappa)u_h^*+\tfrac{\theta}{2}\ln\tfrac{1-u_h^*}{1+u_h^*}}.
\end{align*}
Also, the metastable state solution $u_h^*=0$ solves this equation exactly, 
while the desired equilibria $u_h^*$ only approximately satisfies the steady-state equation with an error of $O(\kappa^{2}\tau^2)$.
It is not expected that the IF1 solution always collapses to the trivial  solution, cf. {Figures \ref{fig: Compare IF and CIF in u with tau 0.05}-\ref{fig: Compare IF and CIF in u with tau 0.5}}, especially when some properly large time-steps are employed to accelerate the long-time dynamics.
Actually, this defect will occur in all stages of any IF methods \eqref{scheme: IF stabilized} due to the fact \eqref{cond: non-equilibria}. That is{ to say}, small time-step sizes are always required for IF methods 
to obtain acceptable equilibria $u_{h}^{*}$ satisfying the approximate steady-state equations at all stages.

Moreover, as pointed out in \cite[Section 3]{LiLiJuFeng:2021SISC} and the existing works mentioned above, the energy dissipation property of the IF method \eqref{scheme: IF stabilized} remains open up to now. Actually, it always can not 
preserve the original energy dissipation law
\eqref{eq: continuous energy law} although it may admit certain modified energy dissipation law, cf. Theorem \ref{thm: energy stable IF1} on the energy stability of the IF1 method. 
To overcome the above two disadvantages, we will develop some corrected IF methods by using two classes of difference correction, including the telescopic correction and nonlinear-term translation correction, enforcing the preservation of steady-state solution. Then the recent theoretical framework in \cite{Liao:2024arxiv} will be applied to evaluate the original energy dissipation properties of the corrected IF methods by using the associated differential forms and differentiation matrices. As applications, some new methods up to fourth-order are constructed by applying the difference correction strategies to some popular IF methods \cite{JuLiQiaoYang:2021JCP,LiLiJuFeng:2021SISC,
  ZhangYanQianSong:2021ANM,ZhangYanQianSong:2022CMAME} having nonnegative Butcher coefficients.  

In the next section, we introduce two different corrections for the IF1 method and examine their energy dissipation properties theoretically and numerically. The difference corrections of second-order IF (IF2) will be investigated in Section \ref{sec:CIF2} and further extensions to high-order IF methods will be addressed in Section \ref{sec:CIF3}.  
Numerical experiments are presented in Section \ref{sec:simulations} to demonstrate the improved performance of our methods. Some conclusions and open issues are included in the last section.

\section{Difference correction of IF1 scheme}
\label{sec:CIF1}

In this section, we will investigate the energy dissipation properties of the IF1 method and two different corrections. To do so, we recall some preliminary results. 

\begin{lemma}\cite[Lemma 2.1]{Liao:2024arxiv}\label{lemma: origional energy derivation}
  If the nonlinear function $g$ is Lipschitz-continuous with a constant $\ell_{g}>0$ and 
  the stabilized parameter $\kappa\ge2\ell_g$, then
  \begin{align*}
    \myinnerb{u-v,g_{\kappa}(v)-\tfrac12L_{\kappa} (u+v)} \leq E[v]-E[u],
  \end{align*}
  where the energy $E[\,\cdot\,]$ is defined in \eqref{eq: continuous energy}.
\end{lemma}
 \begin{lemma}\cite[Theorem 2.1]{Liao:2024arxiv}\label{lemma: EERK energy stability} 
  Assume that the nonlinear function $g$ is Lipschitz-continuous with a constant $\ell_{g}>0$ and  
  the stabilized parameter $\kappa$ in \eqref{def: stabilized parameter} is chosen properly large such that $\kappa\ge 2\ell_g$.
  Consider the following steady-state preserving explicit exponential Runge-Kutta method
  \begin{align}\label{Scheme: general EERK stabilized2}     
    U^{n,i+1}=U^{n,1}+\sum_{j=1}^{i}\hat{a}_{i+1,j}(-{\tau}L_{\kappa})\kbra{{\tau}g_{\kappa}(U^{n,j})-{\tau}L_{\kappa}U^{n,1}}
    \quad\text{for $1\le i\le s$.}
  \end{align}
  Let the lower triangular matrix \lan{$\widehat{A}(z):=\kbra{\hat{a}_{i+1,j}(z)}_{i,j=1}^s$ be the coefficient matrix}, $E_s:=(1_{i\ge j})_{s\times s}$ be the lower triangular matrix full of element 1 and $I$ be the identity matrix. Then the method \eqref{Scheme: general EERK stabilized2} has the following differential form
  \begin{align*}    
    \sum_{\ell=1}^{k}d_{k\ell}(-{\tau}L_{\kappa})\delta_{\tau}U^{n,\ell+1}
    = {\tau}g_{\kappa}(U^{n,k})-\frac{\tau}2L_{\kappa}(U^{n,k+1}+U^{n,k})
    \quad\text{for $1\le k\le s$,}
  \end{align*} 
  where the difference $\delta_{\tau}U^{n,\ell+1}:=U^{n,\ell+1}-U^{n,\ell}$ and the associated differentiation  matrix  $D=(d_{k\ell})_{s\times s}$ is defined as follows, also see \cite[Theorem 2.1]{FuShenYang:2024arxiv-ETDRK},
  \begin{align}\label{def: the differentiation matrix}      
    D(z):=\widehat{A}(z)^{-1}E_s+zE_s-\frac{z}{2}I.
  \end{align}
    If the symmetric part $\mathcal{S}(D;z):=\frac{1}{2}\kbrat{D(z)+D(z)^T}$ is positive (semi-)definite  for $z\le0$, 
  then the explicit exponential  method \eqref{Scheme: general EERK stabilized2} preserves 
  the original energy dissipation law \eqref{eq: continuous energy law}
  at all stages in the sense that
  \begin{align}\label{thmResult: stage energy laws}     
    E[U^{n,j+1}]-E[U^{n,1}]\le-\frac1{\tau}\myinnerb{\delta_{\tau}\vec{U}_{n,j+1},
      D_{j}(-{\tau}L_{\kappa})\delta_{\tau}\vec{U}_{n,j+1}}
    \quad\text{for $1\le j\le s$,}
  \end{align}
  where  $\delta_{\tau}\vec{U}_{n,j+1}:=(\delta_{\tau}U^{n,2},\delta_{\tau}U^{n,3},\cdots,\delta_{\tau}U^{n,j+1})^T$ and $D_{j}(z):=D[1:j,1:j]$ is the $j$-th sequential sub-matrix of the differentiation matrix $D(z)$. 
\end{lemma}


\lan{It is to remark that, to make the subsequent presentation more concise, Lemmas \ref{lemma: origional energy derivation}-\ref{lemma: EERK energy stability} assume that the nonlinear function $g$ is Lipschitz continuous with a constant $\ell_{g}>0$, cf. the recent discussions in \cite[Section 4]{FuShenYang:2024arxiv-ETDRK} or \cite[Subsection 2.2]{FuYang:2022JCP-ETDRK2}. This assumption is really limited for the problem \eqref{problem: autonomous} with the Ginzburg-Landau double-well potential or the  Flory-Huggins potential mentioned in Section 1. One of remedies for weakening this assumption is to establish the maximum bound principle of the involved numerical methods, see the open issue (a) in the last section. Currently, we verify the maximum bound principle numerically for our methods, see the numerical results in Subsection \ref{subsec:test IF1}, Subsection \ref{subsec:test IF2} and Section \ref{sec:simulations}; while we will explore it theoretically in a forthcoming report. As the closely related issue, determining the minimum stabilized parameter $\kappa$ is also practically useful although we impose a properly large $\kappa$ in Lemmas \ref{lemma: origional energy derivation}-\ref{lemma: EERK energy stability}. In this sense,  if a corrected IF method is proven to maintain the original energy dissipation law \eqref{eq: continuous energy law} unconditionally, we mean that this method can be stabilized by using the linear regularization approach \eqref{def: stabilized parameter} with a properly large parameter $\kappa$ (a large $\kappa$ may not be necessary in practical calculations).}

\subsection{IF1 method}
The well-known IF1 method reads
\begin{align}\label{scheme: IF1}
  u_{h}^{n}= e^{- \tau L_{\kappa}} u_{h}^{n-1}+\tau  e^{- \tau L_{\kappa}} g_{\kappa}(u_{h}^{n-1}) \quad\text{for }n\geq1.
\end{align}
As mentioned, it is not steady-state preserving and thus we can not apply Lemma \ref{lemma: EERK energy stability} to establish the energy dissipation law.
On the other hand, one can reformulate \eqref{scheme: IF1} into an equivalent form
\begin{align}\label{eq: IF1 DifferentialForm}
  \brab{I+\frac{\tau}{2}L_{\kappa}} \delta_{\tau} u_{h}^{n} + \brab{e^{\tau L_{\kappa}}-\tau L_{\kappa}-I} u_{h}^{n} = \tau g_{\kappa}(u_{h}^{n-1})- \frac{\tau}{2} L_{\kappa} (u_{h}^{n}+u_{h}^{n-1})
\end{align}
for $n\geq 1.$ With the help of Lemma \ref{lemma: origional energy derivation}, 
we have the following result.
\begin{theorem}\label{thm: energy stable IF1}
  The solution of IF1 scheme \eqref{scheme: IF1} \lan{with the stabilized parameter $\kappa\ge 2\ell_g$} satisfies
  \begin{align*}
  \hE[u_{h}^{n}]-\hE[u_{h}^{n-1}]\le-\frac{1}{2\tau}\myinnerb{\brat{I+e^{\tau L_{\kappa}}} \delta_{\tau}u_{h}^{n},\delta_{\tau} u_{h}^{n}} \quad\text{for $n\geq 1,$}
  \end{align*}
  where the modified energy functional $\hE[\,\cdot\,]$ is defined by
  \begin{align*}
    \hE[u_{h}^{n}]:=E[u_{h}^{n}]
   +\frac{1}{2\tau}\myinnerb{\brat{ e^{\tau L_{\kappa}}-\tau L_{\kappa}-I} u_{h}^{n},u_{h}^{n}}\quad\text{for $n\geq 0$}.
  \end{align*}
\end{theorem}
\begin{proof}
  Making the inner product of the equality \eqref{eq: IF1 DifferentialForm} with $\frac{1}{\tau}\delta_{\tau}u_{h}^{n}$ yields
  \begin{align*}
    & \frac{1}{\tau}\myinnerb{\brat{I+\frac{\tau}{2}L_{\kappa}} \delta_{\tau} u_{h}^{n}, \delta_{\tau}u_{h}^{n}} +\frac{1}{\tau} \myinnerb{\brat{e^{\tau L_{\kappa}}-\tau L_{\kappa}-I} u_{h}^{n},\delta_{\tau} u_{h}^{n}} \notag \\
    &\qquad\qquad = \myinnerb{g_{\kappa}(u_{h}^{n-1})- \frac{1}{2} L_{\kappa} (u_{h}^{n}+u_{h}^{n-1}), \delta_{\tau} u_{h}^{n}}\le E[u_{h}^{n-1}]-E[u_{h}^{n}] \quad\text{for $n\ge1$,}
  \end{align*}
  where Lemma \ref{lemma: origional energy derivation} was used.
  Since the operator $\beta(\tau L_{\kappa}):=e^{\tau L_{\kappa}}-\tau L_{\kappa}-I$ is positive definite, 
  one can apply the formula $2a(a-b)=a^{2}-b^{2}+(a-b)^{2}$ to find that
  \begin{align*}
    2\myinnerb{\beta(\tau L_{\kappa})u_{h}^{n}, \delta_{\tau} u_{h}^{n}} = \myinnerb{\beta(\tau L_{\kappa})u_{h}^{n},u_{h}^{n}}
    - \myinnerb{\beta(\tau L_{\kappa}) u_{h}^{n-1},u_{h}^{n-1}}+\myinnerb{\beta(\tau L_{\kappa}) \delta_{\tau}u_{h}^{n},\delta_{\tau}u_{h}^{n}}.
  \end{align*}
  It is easy to check that 
  \begin{align*}
    \hE[u_{h}^{n}]-\hE[u_{h}^{n-1}]   
    \le-\frac{1}{\tau}\myinnerb{\brat{I+\frac{\tau}{2}L_{\kappa}}\delta_{\tau}u_{h}^{n},\delta_{\tau} u_{h}^{n}}-\frac{1}{2\tau}\myinnerb{\beta(\tau L_{\kappa}) \delta_{\tau}u_{h}^{n},\delta_{\tau}u_{h}^{n}},
 \end{align*}
  which leads to the claimed inequality and completes the proof.
\end{proof}

Theorem \ref{thm: energy stable IF1} says that the IF1 scheme \eqref{scheme: IF1} is energy stable; but this does not necessarily guarantee the decrease of the original energy $E[u_{h}^{n}]$ since the modified energy $\hE[u_{h}^{n}]$ introduces an additional term of order $O(\kappa^{2}\tau)$ to the original energy, see Figures \ref{fig: accuracy and IF1 energy tau 0.2 and 0.5}(b)-(c). Inspired by Lemma \ref{lemma: EERK energy stability}, it is expected that one can remedy it by making simple modification to the IF1 scheme \eqref{scheme: IF1} such that the resulting corrected IF scheme is steady-state preserving. 

This idea is very simple, while the possible modifications are generally diverse. For instance, Krogstad \cite{Krogstad:2005JCP} proposed a generalization of \lan{the} IF method,
and in particular constructed multistep-type methods with several orders of magnitude improved accuracy. Here we are to correct the IF method \eqref{scheme: IF stabilized} within the framework of one-step method. For simplicity, we will restrict ourselves to consider the following two classes of modifications:
\begin{description}
  \item[($\mathrm{T}$-type)] We modify the coefficient of $u_h^n$ by  introducing an undetermined coefficient $\chi(-\tau L_{\kappa})$ such that the new scheme $\chi(-\tau L_{\kappa})u_{h}^{n}= e^{- \tau L_{\kappa}} u_{h}^{n-1}+\tau  e^{- \tau L_{\kappa}} g_{\kappa}(u_{h}^{n-1})$ is steady-state preserving. As seen, it is equivalent to correct the IF1 solution $u_h^n$ by a telescopic factor $\chi(-\tau L_{\kappa})$, so we call this modification as telescopic-type correction.
    \item[($\mathrm{N}$-type)] Also, we can modify the coefficient of nonlinear term  $g_{\kappa}(u_{h}^{n-1})$, that is, introduce an undetermined coefficient $\chi(-\tau L_{\kappa})$ such that the new scheme $u_{h}^{n}=e^{- \tau L_{\kappa}} u_{h}^{n-1}+\tau  \chi(-\tau L_{\kappa})g_{\kappa}(u_{h}^{n-1})$ is steady-state preserving. We call this modification as nonlinear-term translation correction.    
\end{description}

In the next subsection, we will investigate these two correction strategies in detail. 
\lan{In general, both modifications can preserve the original energy dissipation law \eqref{eq: continuous energy law}. However, we find in numerical simulations that the discrete energy curves generated by the two classes of corrections have obvious differences from the continuous (reference) energy curve. It is really interesting to pick out a ``better'' method that the associated discrete energy would be closer to the continuous energy with the same space-time discretization parameters; at the same time, this task is beyond our current scope of this article and will be explored in future works, see the open issue (b) in the last section.}
In fact, the reason we examine multiple correction strategies is that we don't know which one should be better.
As a preliminary attempt, we also consider a linear-term translation correction, like $u_{h}^{n}=\chi(-\tau L_{\kappa}) u_{h}^{n-1}+\tau e^{- \tau L_{\kappa}}   g_{\kappa}(u_{h}^{n-1})$; however, the resulting corrected IF methods seem practically and theoretically inferior to the presented $\mathrm{T}$-type and $\mathrm{N}$-type corrections. 

\subsection{Two corrections of IF1 scheme}\label{subsec:CIF1}

\subsubsection{Telescopic-type correction}
We introduce an undetermined coefficient $\chi_{\cte}^{(1)}(-\tau L_{\kappa})$ of $u_{h}^{n}$ at $t_{n}$ and the corrected method reads
\begin{align*}
  \chi_{\cte}^{(1)}(-\tau L_{\kappa}) u_{h}^{n} = e^{-\tau L_{\kappa}} u_{h}^{n-1} + \tau e^{-\tau L_{\kappa}} g_{\kappa}(u_{h}^{n-1})  \quad\text{for $n \ge 1$.}
\end{align*}
Requiring $u_{h}^{k}=u_{h}^{*}$ for all $k\geq0$ immediately yields 
 $\chi_{\cte}^{(1)}(z) = (1-z)e^{z}$.
The resulting  telescopic corrected IF1 (in short, TIF1) scheme reads
\begin{align}\label{scheme: TIF1}
 u_{h}^{n} = (I + \tau L_{\kappa})^{-1}u_{h}^{n-1} + \tau (I + \tau L_{\kappa})^{-1} g_{\kappa}(u_{h}^{n-1}) \quad\text{for $n \ge 1$}.
\end{align}
It is nothing new but  the stabilized semi-implicit Euler scheme. One can reformulate it into
\begin{align*}
  u_{h}^{n}=u_{h}^{n-1}+(I + \tau L_{\kappa})^{-1}\kbra{\tau g_{\kappa}(u_{h}^{n-1})-\tau L_{\kappa}u_{h}^{n-1}}  \quad\text{for $n \ge 1$.}
\end{align*}
This formulation takes the steady-state preserving form \eqref{Scheme: general EERK stabilized2} with the coefficient $\hat{a}_{21}(z)=(1-z)^{-1}$. Thus the definition \eqref{def: the differentiation matrix} gives the differentiation matrix (scalar) 
\begin{align*}
  D_{\cte}^{(1)}(z)= 1-\frac{z}{2}\ge1 \quad\text{for $z\le0$}.
\end{align*}
 Lemma \ref{lemma: EERK energy stability} confirms that the stabilized semi-explicit Euler scheme \eqref{scheme: TIF1} preserves the original energy dissipation law \eqref{eq: continuous energy law}. 
Here and hereafter, the superscript $(p)$ is always used to indicate the formal order of method, and the subscripts $\mathrm{T}$ and $\mathrm{N}$ represent different correction approaches. That is to say, $D_{\cte}^{(p)}$ denotes the differentiation matrix  of a formally $p$-th order TIF method.


\subsubsection{Nonlinear-term translation correction}
Introducing an undetermined coefficient $\chi_{\cnt}^{(1)}(-\tau L_{\kappa})$ of the nonlinear term $g_{\kappa}(u_{h}^{n-1})$ arrives at
\begin{align*}
  u_{h}^{n} = e^{-\tau L_{\kappa}} u_{h}^{n-1} + \tau \chi_{\cnt}^{(1)}(-\tau L_{\kappa}) g_{\kappa}(u_{h}^{n-1})  \quad\text{for $n \ge 1$}.
\end{align*}
Requiring $u_{h}^{k}=u_{h}^{*}$ for all $k\geq0$ gives the correction coefficient $\chi_{\cnt}^{(1)}(z) =(e^{z}-1)/z$.
The resulting nonlinear-term corrected IF1 (in short, NIF1) scheme reads
\begin{align}\label{scheme: NIF1}
  u_{h}^{n} = u_{h}^{n-1}+ (\tau L_{\kappa})^{-1} (I - e^{-\tau L_{\kappa}})\kbra{\tau g_{\kappa}(u_{h}^{n-1})-\tau L_{\kappa}u_{h}^{n-1}}  \quad\text{for $n \ge 1$,}
\end{align}
which is just the widespread ETD1 scheme \cite{CoxMatthews:2002JCP,DuJuLiQiao:2019SINUM,DuJuLiQiao:2021SIREV,LiLiJuFeng:2021SISC,Liao:2024arxiv}.
It takes the steady-state preserving form \eqref{Scheme: general EERK stabilized2} with the coefficient $\hat{a}_{21}(z)=(e^{z}-1)/z$. The definition \eqref{def: the differentiation matrix} gives the differentiation matrix (scalar) 
\begin{align*}
     D_{\cnt}^{(1)}(z) := \frac{z}{e^{z}-1}+\frac{z}{2}\ge1 \quad\text{for $z\le0$},
\end{align*}
due to the fact $\frac{z}{e^{b z}-1}\ge\frac{1}{b}$ for $b>0$ and $z\le0$.
 Lemma \ref{lemma: EERK energy stability} shows that the NIF1 scheme \eqref{scheme: NIF1} preserves the energy dissipation law \eqref{eq: continuous energy law}, cf. \cite[Corollary 2.1]{Liao:2024arxiv}.

As the end of this subsection, we summarize the above results as follows.
\begin{theorem}\label{thm: energy stability CIF1}
  Assume that $g$ is Lipschitz-continuous with a constant $\ell_{g}>0$ and the stabilized parameter $\kappa$ in \eqref{def: stabilized parameter} is chosen properly large such that $\kappa\ge 2\ell_g$. The TIF1  \eqref{scheme: TIF1} and NIF1  \eqref{scheme: NIF1} schemes preserve the original energy dissipation law \eqref{eq: continuous energy law} in the sense that
  \begin{align*}
    E[u_h^{n}]-E[u_h^{n-1}]\le-\frac1{\tau} \myinnerb{\delta_{\tau} u_h^{n}, D^{(1)}(-\tau L_{\kappa}) \delta_{\tau}u_h^{n}}\quad \text{for $n\ge1$},
  \end{align*}
  where $D^{(1)}$ represents $D_{\cte}^{(1)}$ and $D_{\cnt}^{(1)}$ for the two different corrections, respectively. 
\end{theorem} 

\subsection{Tests of corrected IF1 schemes}\label{subsec:test IF1}

Before we turn to construct new corrections for high-order IF methods, some preliminary tests are presented in this subsection to emphasize the practical effectiveness of two corrected IF1 methods.

\begin{example}
\cite{LiaoTangZhou:2020SINUM}
\label{example 1} Consider the Allen-Cahn model $\partial_tu=\epsilon^2\partial_{xx} u-u^3+u$ on $\Omega=(-1,1)$ with $\epsilon=0.1$ subject to the 2-periodic initial data $u_0=  -\tanh(((x-0.3)^2-0.2^2)/\epsilon)\tanh(((x+0.3)^2-0.2^2)/\epsilon)$. The center difference approximation is used to discrete the spatial operator with the length $h=0.01$.  \lan{The solution preserves the maximum bound principle, that is, $\mynorm{u_h(t)}_{\infty}\le1$ since $\mynorm{u_h^0}_{\infty}\le1$, cf. \cite[Section 2]{DuJuLiQiao:2021SIREV}}.
\end{example}

Following the suggestions in \cite{LiLiJuFeng:2021SISC,ZhangYanQianSong:2022CMAME,GaoZhangQianSong:ANM2023-IF,JuLiQiaoYang:2021JCP}, we always choose a stabilized parameter $\kappa>0$ in our tests. We run the IF1 \eqref{scheme: IF1}, TIF1 \eqref{scheme: TIF1} and NIF1 \eqref{scheme: NIF1} schemes with $\tau=0.001$ for $\kappa=1$, 2 and 4. The three schemes work well up to $T=40$, and the corresponding solution and energy curves (omitted here) are hard to distinguish from each other. 
\lan{The numerical solution $u_h^\star$ of the NIF1 scheme \eqref{scheme: NIF1} computed with the small time-step size $\tau=0.001$ is taken as the reference solution in the following tests.}

\begin{figure}[htb!]
  \centering
  \subfigure[\lan{Numerical errors}]{
    \includegraphics[width=2.05in,height=1.5in]{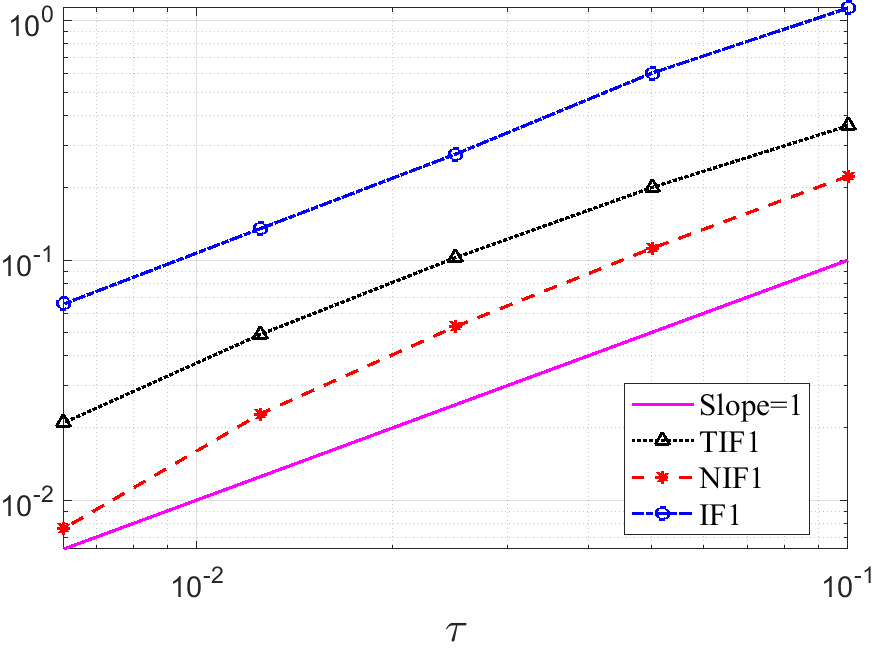}}
  \subfigure[Energy for $\tau=0.2$]{
    \includegraphics[width=2.05in,height=1.52in]{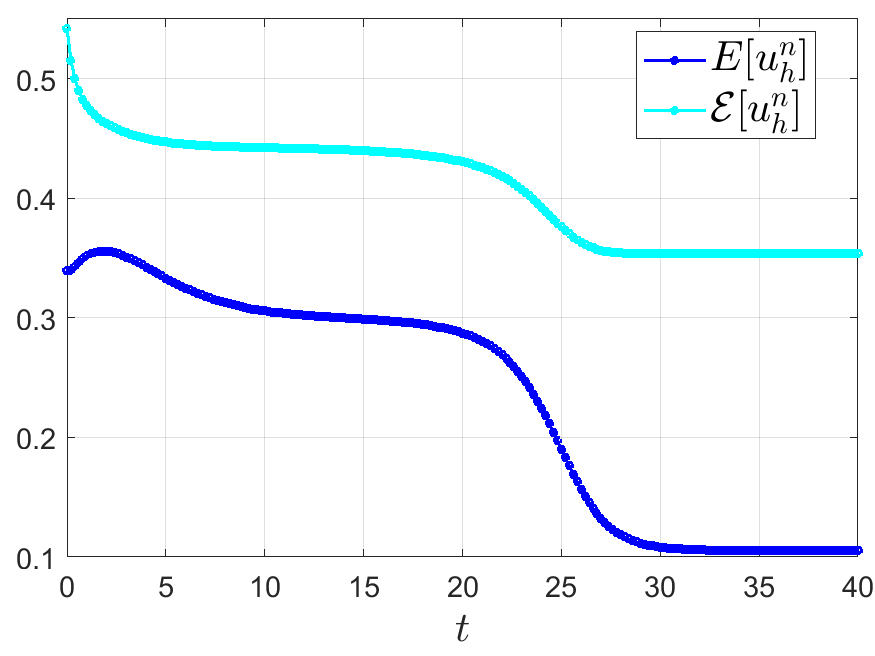}}
  \subfigure[Energy for $\tau=0.5$]{
    \includegraphics[width=2.05in,height=1.52in]{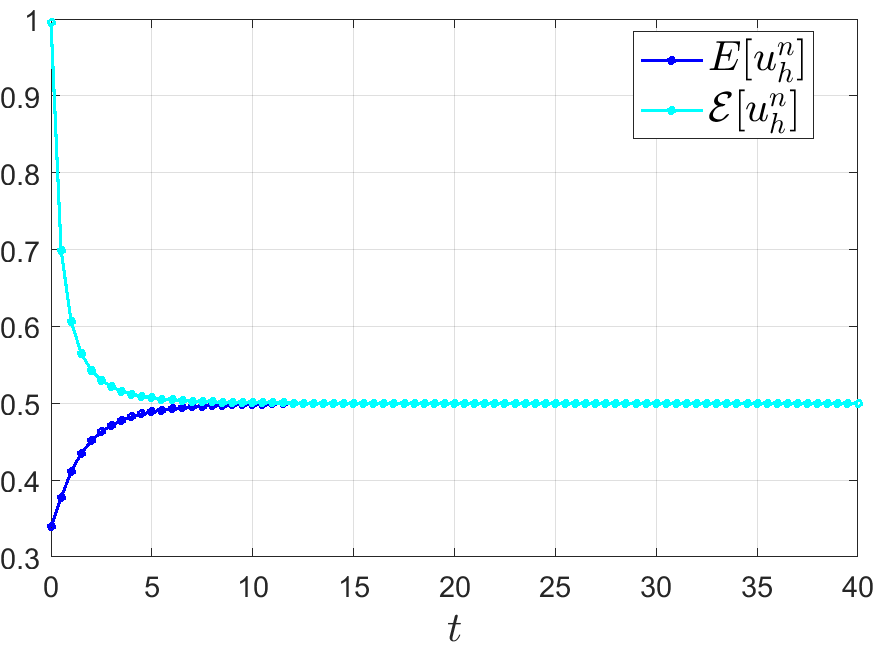}}
  \caption{\qing{Errors of corrected IF1 schemes and energy dissipation of IF1 scheme.}}
  \label{fig: accuracy and IF1 energy tau 0.2 and 0.5}
\end{figure}

At first, we run the IF1 and two corrected IF1 schemes for different time-steps $\tau=2^{-k}/10$ $(0\le k\le 4)$ with the final time $T=20$ and $\kappa=4$ to test the temporal accuracy (here we set a small spatial length $h = 1/500$ to make the spatial error negligible). The numerical errors \lan{$e(\tau) := \max_{1\leq n \leq N}\|u_h^n-u_h^\star\|_\infty$} of the three schemes, depicted in Figure \ref{fig: accuracy and IF1 energy tau 0.2 and 0.5}(a), confirm the first-order time accuracy.  Setting a stabilized parameter $\kappa=2$ and the spatial length $h=0.1$, we record the discrete modified energy $\mathcal{E}[u_h^n]$ (and original energy $E[u_h^n]$) of the IF1 scheme \eqref{scheme: IF1} in Figure \ref{fig: accuracy and IF1 energy tau 0.2 and 0.5}(b)-(c) for two different time-steps $\tau=0.2$ and $0.5$, respectively. The decrease of $\mathcal{E}[u_h^n]$ supports Theorem \ref{thm: energy stable IF1}; however, the increase of original energy for properly large $\tau=0.5$ is surprising and mysterious to us. It is to note that, with a fixed time-step size $\tau=0.5$, one can observe similar behaviors (omitted here) in Figure \ref{fig: accuracy and IF1 energy tau 0.2 and 0.5}(b)-(c)  for different stabilized parameters $\kappa=1$ and 2. 

\begin{figure}[htb!]
  \centering
  \subfigure[Final solution $u_h^N$]{
    \includegraphics[width=2.05in,height=1.5in]{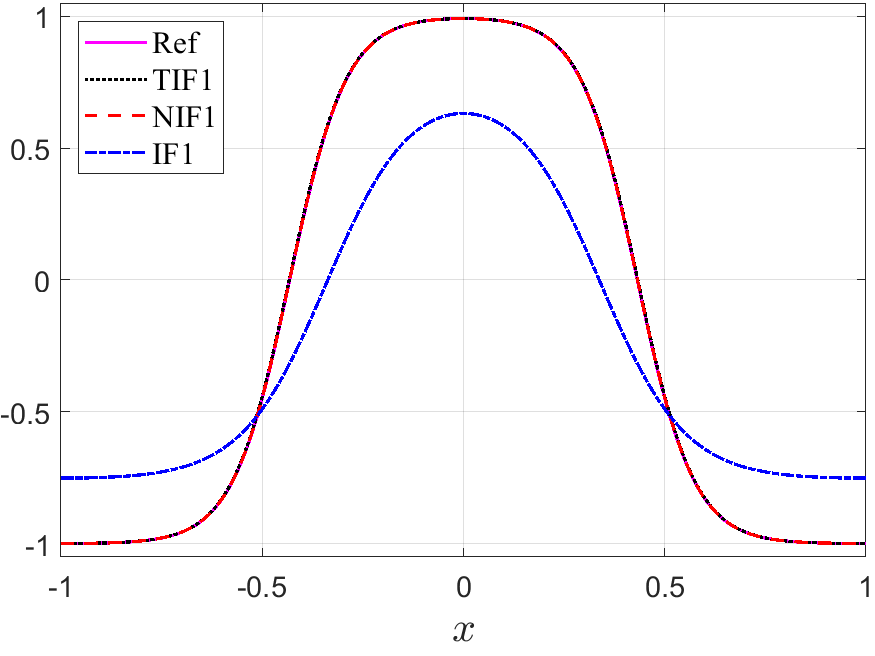}}
  \subfigure[Discrete energy $E\kbrat{u_h^n}$]{
    \includegraphics[width=2.05in,height=1.5in]{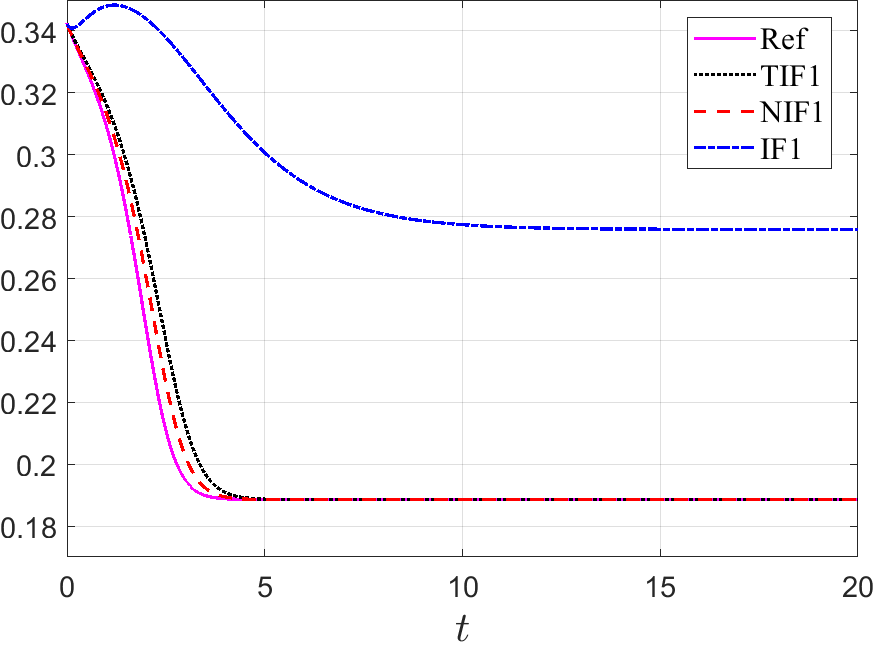}}
  \subfigure[Maximum norm $\mynormt{u_h^n}_{\infty}$]{
    \includegraphics[width=2.05in,height=1.52in]{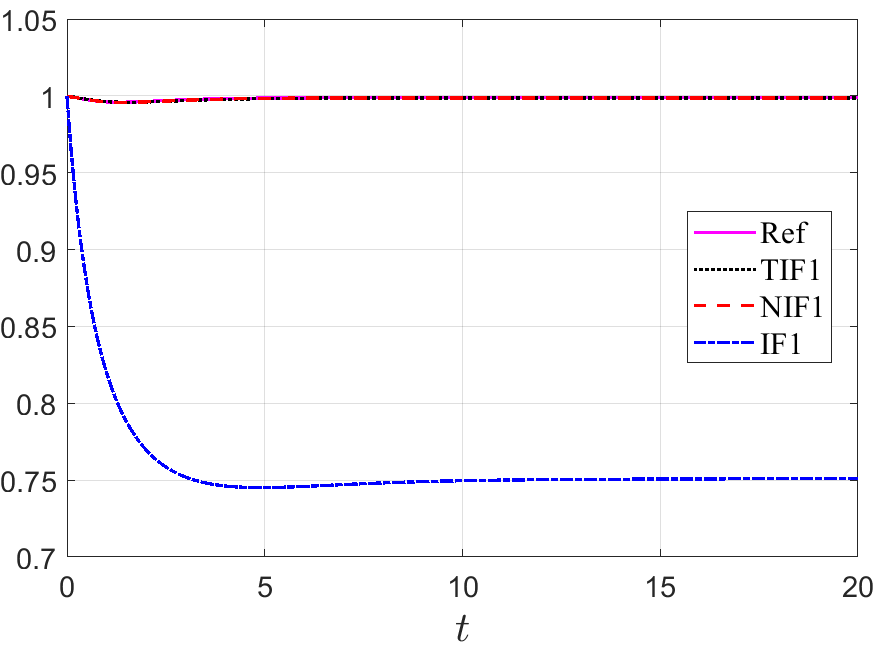}}
  \caption{Comparisons of IF1 and corrected IF1 schemes for $\tau=0.05$.}
  \label{fig: Compare IF and CIF in u with tau 0.05}
\end{figure}
  
\begin{figure}[htb!]
  \centering
  \subfigure[Final solution $u_h^N$]{
    \includegraphics[width=2.05in,height=1.5in]{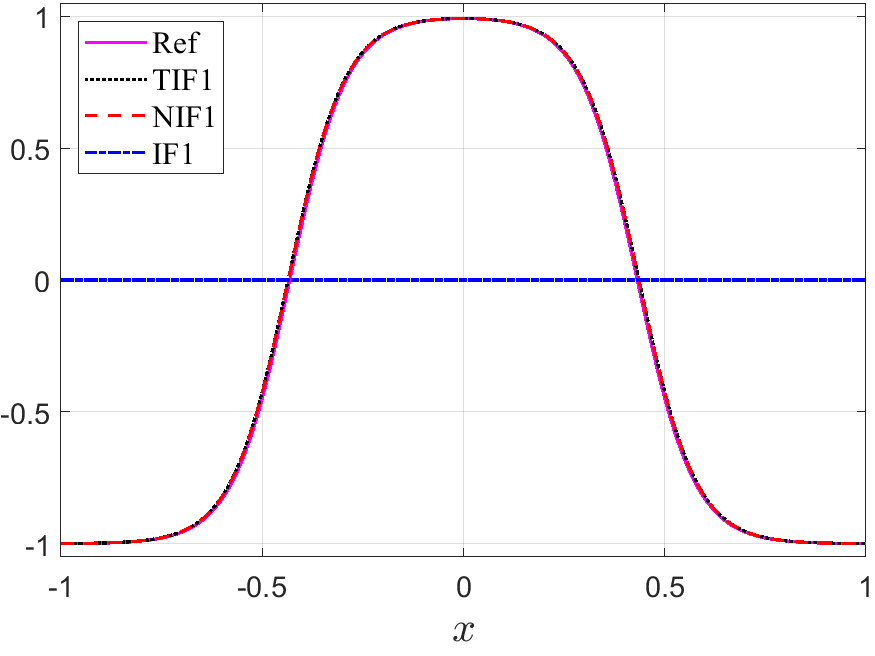}}
  \subfigure[Discrete energy $E\kbrat{u_h^n}$]{
    \includegraphics[width=2.05in,height=1.5in]{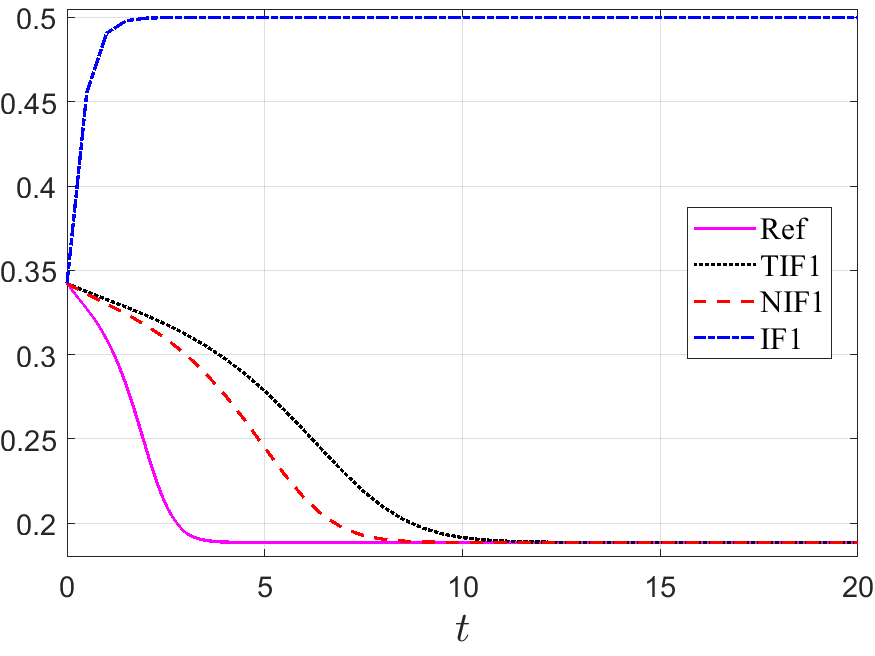}}
  \subfigure[Maximum norm $\mynormt{u_h^n}_{\infty}$]{
    \includegraphics[width=2.05in,height=1.5in]{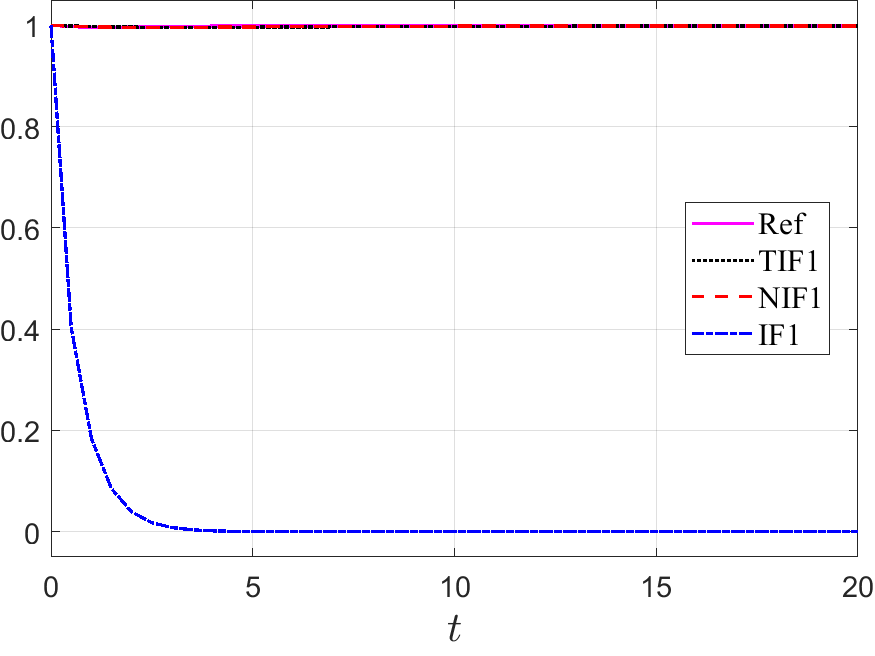}}
  \caption{Comparisons of IF1 and corrected IF1 schemes for $\tau=0.5$.}
  \label{fig: Compare IF and CIF in u with tau 0.5}
\end{figure}

Secondly, we run the \lan{three} first-order schemes to the final time $T=20$ with $\kappa=4$ for two different time-steps $\tau=0.05$ and $0.5$, see Figure \ref{fig: Compare IF and CIF in u with tau 0.05} and Figure \ref{fig: Compare IF and CIF in u with tau 0.5}, respectively. In these two groups of figures, we depict the profiles of final solution $u_h^N$, the discrete energy $E[u_h^n]$ and the maximum norm $\|u_h^n\|_{\infty}$. As expected, two corrected IF1 methods maintain the steady-state solution well, cf. \lan{Figures \ref{fig: Compare IF and CIF in u with tau 0.05}(a)-\ref{fig: Compare IF and CIF in u with tau 0.5}(a)}, while the final solution $u_h^N$ of IF1 scheme \eqref{scheme: IF1} has a bit collapse for the time-step $\tau=0.05$, and it gradually collapses into the metastable state as $\tau$ increases to $0.5$, see Figures \ref{fig: Compare IF and CIF in u with tau 0.05}(a)-\ref{fig: Compare IF and CIF in u with tau 0.5}(a). The similar solution behaviors can be obtained with different stabilized parameters $\kappa=0.5$ and 4 for the fixed step size $\tau=0.5$. This phenomenon would be attributed to the fact that the discrete equilibria $u_h^*$ of IF1 scheme \eqref{scheme: IF1} only approximately satisfies the steady-state equation with an error of $O(\kappa^{2}\tau^2)$. Experimentally, we observe that the TIF1 \eqref{scheme: TIF1} and NIF1 \eqref{scheme: NIF1} schemes always perform well, at least for this example.

\begin{figure}[htb!]
  \centering
  \includegraphics[width=2.4in]{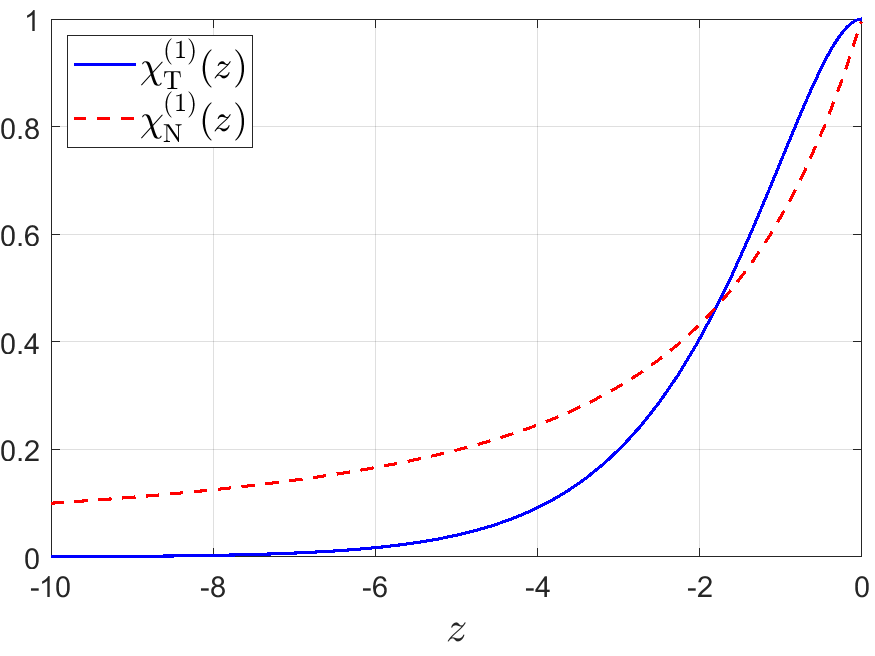}
  \caption{Correction coefficients of two corrected IF1 schemes.}
  \label{fig: corrected IF1 coefficients rate}
\end{figure}

It seems that our corrected IF1 methods also maintain the maximum bound principle of the Allen-Cahn model, see Figures \ref{fig: Compare IF and CIF in u with tau 0.05}(c)-\ref{fig: Compare IF and CIF in u with tau 0.5}(c), at least in the current experimental setting. We know that the contractivity of IF methods is essential to preserve the maximum bound principle of semilinear parabolic problems, cf. \cite{DuJuLiQiao:2019SINUM,DuJuLiQiao:2021SIREV,LiLiJuFeng:2021SISC,ZhangYanQianChenSong:2022JSC,ZhangYanQianSong:2022CMAME} and references therein. By the form \eqref{scheme: IF1}, it is easy to find that the IF1 method is unconditionally contractive in the sense of \cite{MasetZennaro:2009MCOM}. 
It is \lan{worth  mentioning} that the above two correction coefficients $\chi_{\cte}^{(1)}(z)$ and $\chi_{\cnt}^{(1)}(z)$ are always positive and smaller than 1, cf. Figure \ref{fig: corrected IF1 coefficients rate}. They suggest that these corrected IF1 schemes may be also unconditionally contractive. Actually, it is not difficult to show that the stabilized semi-implicit scheme \eqref{scheme: TIF1} and the ETD1 scheme \eqref{scheme: NIF1} are contractive, while detailed discussions are out of our current scope.
  
\section{A unified framework and corrections of IF2 methods}\label{sec:CIF2}
\setcounter{equation}{0}

IF2 methods require two internal stages at least, and have the following general form for solving the stabilized gradient flow problem \eqref{problem: stabilized version}
\begin{subequations}\label{scheme: IF2}
\begin{align}
  U^{n,2}=&\; e^{-c_{2} \tau L_{\kappa}} U^{n,1}+\tau a_{2,1}(-\tau L_{\kappa}) g_{\kappa}(U^{n, 1}), \label{scheme: IF2-1} \\
  U^{n,3}=&\; e^{- \tau L_{\kappa}} U^{n,1}+\tau \sum_{j=1}^2a_{3,j}(-\tau L_{\kappa}) g_{\kappa}(U^{n, j}), \label{scheme: IF2-2}
\end{align}\end{subequations}
where the coefficients $a_{i+1,j}(z):=a_{i+1,j}(0)e^{-(c_{i+1}-c_{j})\tau L_{\kappa}} $ for $1\le j\le i\le 2$.

Following the correction process described in the Subsection \ref{subsec:CIF1}, one can perform  the $\mathrm{T}$-type and $\mathrm{N}$-type corrections to modify the second-stage scheme \eqref{scheme: IF2-1}, and three possible corrections to modify the third-stage scheme \eqref{scheme: IF2-2} since it involves two nonlinear terms.
That is to say, one has a total of $2 \times 3=6$ corrected IF2 schemes, at least.

For the simplicity of presentation, we consider two classes of corrected IF2 schemes by utilizing the same correction strategy at each stage.
In other words, if we correct the second-stage \eqref{scheme: IF2-1} by $\mathrm{T}$-type ($\mathrm{N}$-type) correction, we also do the same type modification to \eqref{scheme: IF2-2}.
Specifically, for the nonlinear-term translation ($\mathrm{N}$-type) correction, we always modify the last nonlinear-term $g_{\kappa}(U^{n, i})$ at the stage $t_{n,i+1}$. In this section, we will present a unified theoretical framework for $s$-stage IF methods with a fixed $s\ge2$ and perform the $\mathrm{T}$-type and $\mathrm{N}$-type corrections to two widespread IF2 methods having nonnegative method coefficients.

\subsection{A unified theoretical framework for difference corrections}\label{sec:unified framework}


\textbf{(Telescopic correction)}
For the stabilized IF methods \eqref{scheme: IF stabilized},
introduce an undetermined coefficient 
$\chi_{\cte,i+1}^{(s)}(-\tau L_{\kappa})$ of $U^{n,i+1}$ at the stage $t_{n,i+1}$ for $1\le i\le s$ such that the telescopic correction of \eqref{scheme: IF stabilized} reads
\begin{align*}
  \chi_{\cte,i+1}^{(s)}(-\tau L_{\kappa})U^{n,i+1}= e^{-c_{i+1} \tau L_{\kappa}} U^{n,1}+\tau \sum_{j=1}^{i} a_{i+1, j}(-\tau L_{\kappa}) g_{\kappa}(U^{n, j})
\end{align*}
for $ 1\le i\le s$. Requiring $U^{n,i}=u_{h}^{*}$ for all $n\geq1$ and $1\le i\le s$ immediately yields the telescopic coefficients
\begin{align}\label{def: general T-type coefficients}
\chi_{\cte,i+1}^{(s)}(z):= e^{c_{i+1}z} -z \sum_{j=1}^ia_{i+1,j}(z) \quad\text{for $ 1\le i\le s$}.
\end{align}
The resulting telescopic corrected IF (in short, TIF) scheme reads
\begin{align}\label{eq: general T-type correction}
  U^{n,i+1}=&\; U^{n,1} + \sum_{j=1}^i\hat{a}^{\cte}_{i+1,j}(-\tau L_{\kappa}) \kbrab{\tau g_{\kappa}(U^{n, j}) - \tau L_{\kappa} U^{n,1}}\quad\text{for $1\le i\le s$,}
\end{align}
which takes the steady-state preserving form \eqref{Scheme: general EERK stabilized2} with the coefficients
\begin{align}\label{eq: general T-type steady-state coefficients}
  \hat{a}^{\cte}_{i+1,j}(z):= \frac{a_{i+1,j}(z)}{\chi_{\cte,i+1}^{(s)}(z)}
  =\frac{a_{i+1,j}(0)e^{-c_{j}z}}{1 -z \sum_{\ell=1}^ia_{i+1,\ell}(0)e^{-c_{\ell}z}}\quad\text{for $1\le j\le i\le s$}.
\end{align}
Thus the definition \eqref{def: the differentiation matrix} gives the associated differentiation matrix, 
\begin{align}\label{eq: general T-type differentiation matrix}    
  D_{\cte}^{(s)}(z):=\widehat{A}_{\cte}^{(s)}(z)^{-1}E_s+zE_s-\frac{z}{2}I,
\end{align} 
where \lan{$\widehat{A}_{\cte}^{(s)}(z):=\kbra{\hat{a}^{\cte}_{i+1,j}(z)}_{i,j=1}^s$ is the coefficient matrix of the TIF scheme \eqref{eq: general T-type correction}}.
Then Lemma \ref{lemma: EERK energy stability} says that if the differentiation matrix $D_{\cte}^{(s)}$ in \eqref{eq: general T-type differentiation matrix} is positive semi-definite, the TIF scheme \eqref{eq: general T-type correction} preserves the original energy dissipation law \eqref{eq: continuous energy law} at all stages.

\qing{Note that, the TIF scheme \eqref{eq: general T-type correction} has the same formal order of original IF method.} Actually, the TIF scheme \eqref{eq: general T-type correction} has the same underlying explicit Runge-Kutta methods to the stabilized IF methods \eqref{scheme: IF stabilized} since the modified coefficients $\hat{a}^{\cte}_{i+1,j}(z)$ in \eqref{eq: general T-type steady-state coefficients} tend to $a_{i+1,j}(0)$ as $z\to0$.

\textbf{(Nonlinear-term translation correction)}
Introduce an undetermined coefficient 
$\chi_{\cnt,i+1}^{(s)}(-\tau L_{\kappa})$ of the nonlinear term $g_{\kappa}(U^{n, i})$ 
at each stage such that the $\mathrm{N}$-type correction of \eqref{scheme: IF stabilized} reads
\begin{align*}
  U^{n,i+1}= e^{-c_{i+1} \tau L_{\kappa}} U^{n,1} +\tau \sum_{j=1}^{i-1} a_{i+1, j}(-\tau L_{\kappa}) g_{\kappa}(U^{n, j}) +\tau \chi_{\cnt,i+1}^{(s)}(-\tau L_{\kappa}) g_{\kappa}(U^{n, i}) \quad\text{for $ 1\le i\le s$}.
\end{align*}
Requiring $U^{n,i}=u_{h}^{*}$ for all $n\geq1$ and $1\le i\le s$ yields the nonlinear-term corrected coefficients
\begin{align}\label{def: general N-type coefficients}
  \chi_{\cnt,i+1}^{(s)}(z):= \frac{e^{c_{i+1}z}-1}{z}-\sum_{j=1}^{i-1}a_{i+1,j}(z)
  \quad\text{for $ 1\le i\le s$}.
\end{align}
The resulting nonlinear-term corrected IF (NIF) scheme reads
\begin{align}\label{eq: general N-type correction}
  U^{n,i+1}=&\; U^{n,1} + \sum_{j=1}^i\hat{a}^{\cnt}_{i+1,j}(-\tau L_{\kappa}) \kbrab{\tau g_{\kappa}(U^{n, j}) - \tau L_{\kappa} U^{n,1}}\quad\text{for $1\le i\le s$,}
\end{align}
which takes the steady-state preserving form of \eqref{Scheme: general EERK stabilized2} with the coefficients
\begin{subequations}\label{eq: general N-type steady-state coefficients}
\begin{align}
  \hat{a}^{\cnt}_{i+1,j}(z):=&\, a_{i+1,j}(z)
  =a_{i+1,j}(0)e^{(c_{i+1}-c_{j})z}
  \quad\text{for $2\le j+1\le i\le s,$}\\
  \hat{a}^{\cnt}_{i+1,i}(z):=&\,\frac{e^{c_{i+1}z}-1}{z}
  -\sum_{j=1}^{i-1}a_{i+1,j}(0)e^{(c_{i+1}-c_{j})z}
  \quad\text{for $1\le i\le s$.}
\end{align}
\end{subequations}
Thus the definition \eqref{def: the differentiation matrix} gives the associated differentiation matrix
\begin{align}\label{eq: general N-type differentiation matrix}    
  D_{\cnt}^{(s)}(z):=\widehat{A}_{\cnt}^{(s)}(z)^{-1}E_s+zE_s-\frac{z}{2}I,
\end{align} 
where \lan{$\widehat{A}_{\cnt}^{(s)}(z):=\kbra{\hat{a}^{\cnt}_{i+1,j}(z)}_{i,j=1}^s$ is the coefficient matrix of the NIF scheme \eqref{eq: general N-type correction}}.
Lemma \ref{lemma: EERK energy stability} says that if the differentiation matrix $D_{\cnt}^{(s)}$ in \eqref{eq: general N-type differentiation matrix} is positive semi-definite, the NIF scheme \eqref{eq: general N-type correction} preserves the original energy dissipation law \eqref{eq: continuous energy law} at all stages.

It is easy to check that, as $z\to0$, the NIF coefficients in \eqref{eq: general N-type steady-state coefficients} tend to $a_{i+1,i}(0)$ for $1\le j\le i\le s$, due to the consistency condition \eqref{cond: equilibria underlying RK} and the simple fact, $\lim_{z\to0}\frac{e^{c_{i+1}z}-1}{z} = c_{i+1}$. It means that the NIF scheme \eqref{eq: general N-type correction} has the same underlying explicit Runge-Kutta methods to the stabilized IF methods \eqref{scheme: IF stabilized} \qing{so that the NIF scheme \eqref{eq: general N-type correction} has the same formal order of original IF method.}

\subsection{Corrections of Heun's IF2 method}

This subsection considers the corrections of second-order Heun's IF2 scheme \cite{JuLiQiaoYang:2021JCP,ZhangYanQianSong:2021ANM} with the following Butcher tableau of the underlying explicit Heun's method \cite{Heun:1900ZMP,Runge:1895MA}
\begin{align*}
  \text{Second-order Heun:}\quad \begin{array}{c|cc}
    0 & \\
    1 & 1 \\
    \hline\\[-10pt] & \frac{1}{2} & \frac{1}{2}
  \end{array}\;.
\end{align*}

(\textbf{$\mathrm{T}$-type correction}) By the definition \eqref{def: general T-type coefficients}, the telescopic correction coefficients read
\begin{align}\label{def: Heun IF2 T-type coefficients}
  \chi_{\cte,2}^{(2,H)}(z):= (1-z)e^{z},\quad
  \chi_{\cte,3}^{(2,H)}(z):= e^{z}-\tfrac{z}{2}(1+e^{z}).
\end{align}
By  \eqref{eq: general T-type steady-state coefficients}, the resulting TIF2-Heun scheme has
the following method coefficients
\begin{align}\label{def: Heun TIF2 coefficients}
  \hat{a}^{\cte}_{21}(z)=\tfrac{1}{1-z},\quad
  \hat{a}^{\cte}_{31}(z)=\tfrac{1}{2 -z(1+e^{-z})},\quad
  \hat{a}^{\cte}_{32}(z)=\tfrac{1}{2e^{z} -z(1+e^{z})}\,.
\end{align}
Then using \eqref{eq: general T-type differentiation matrix}, we have the associated differentiation matrix
\begin{align*}    
  D_{\cte}^{(2,H)}(z)=\begin{pmatrix}
 1-\frac{z}{2} & 0 \\[2pt]
 e^z & e^z (2-z)-\frac{z}{2}
  \end{pmatrix}\,.
\end{align*}
Obviously, the first leading principal minor of the symmetric matrix $\mathcal{S}(D_{\cte}^{(2,H)};z)$ is greater than 1 for $z\le0$. Also, the second leading principal minor is positive, that is,
\begin{align*}
  \mathrm{Det}\kbrab{\mathcal{S}(D_{\cte}^{(2,H)};z)}
  =&\, \tfrac{1}{4} \kbrab{(2 e^z+1) z^2-2z(4 e^z+1)+ {e^z(8-e^{z})}}
  \ge\tfrac{7}{4}\quad\text{for $z\le0$.}
\end{align*}
Thus the differentiation matrix $D_{\cte}^{(2,H)}(z)$ is positive definite for $z\le0$.

(\textbf{$\mathrm{N}$-type correction}) By the definition \eqref{def: general N-type coefficients},
one can obtain the nonlinear-term correction coefficients
\begin{align}\label{def: Heun IF2 N-type coefficients}
  \chi_{\cnt,2}^{(2,H)}(z)= \tfrac{e^{z}-1}{z},\quad
  \chi_{\cnt,3}^{(2,H)}(z)= \tfrac{e^{z}-1}{z}-\tfrac{1}{2}e^{z}\,.
\end{align}
By the definition \eqref{eq: general N-type steady-state coefficients}, the resulting NIF2-Heun scheme has
the method coefficients
\begin{align}\label{def:  Heun NIF2 coefficients}
  \hat{a}^{\cnt}_{21}(z)=\tfrac{e^{z}-1}{z},\quad
  \hat{a}^{\cnt}_{31}(z)=\tfrac{1}{2}e^{z},\quad
  \hat{a}^{\cnt}_{32}(z)=\tfrac{e^{z}-1}{z}-\tfrac{1}{2}e^{z}\,.
\end{align}
One can follow \eqref{eq: general N-type differentiation matrix} to find the associated differentiation matrix
\begin{align*}    
  D_{\cnt}^{(2,H)}(z)=\begin{pmatrix}
 \frac{z}{e^z-1}+\frac{z}{2} & 0 \\[2pt]
 \frac{e^z z}{e^z-1} & \frac{z^2e^z - 2z(e^z+1)}{2 e^z (z-2)+4}
  \end{pmatrix}\,.
\end{align*}
It is not difficult to check that $\mathrm{Det}\kbrab{\mathcal{S}(D_{\cnt,1}^{(2,H)};z)}=\frac{z}{e^z-1}+\frac{z}{2}\ge1$ for $z\le0$, and
\begin{align*}
  {\mathrm{Det}\kbrab{\mathcal{S}(D_{\cnt}^{(2,H)};z)}=\tfrac{z^2}{4 (e^z-1)^2}
    \tfrac{2-z-4 e^{z}+2e^{-z}}{ z-2+2e^{-z}}>0\quad\text{for $z\le0$.}}
\end{align*}
Thus the differentiation matrix $D_{\cnt}^{(2,H)}(z)$ is positive definite for $z\le0$.

In summary, one can apply Lemma \ref{lemma: EERK energy stability} to get the following result.
\begin{theorem}\label{thm: energy stability CIF2-Heun}
 Assume that $g$ is Lipschitz-continuous with a constant $\ell_{g}>0$ and the stabilized parameter $\kappa$ in \eqref{def: stabilized parameter} is chosen properly large such that $\kappa\ge 2\ell_g$.  In solving the gradient flow \eqref{problem: stabilized version}, the TIF2-Heun \eqref{def: Heun TIF2 coefficients} and NIF2-Heun  \eqref{def:  Heun NIF2 coefficients} schemes preserve the original energy dissipation law \eqref{eq: continuous energy law} at all stages.
\end{theorem}

\begin{figure}[htb!]
  \centering
  \subfigure[TIF2-Heun]{\includegraphics[width=2.4in]{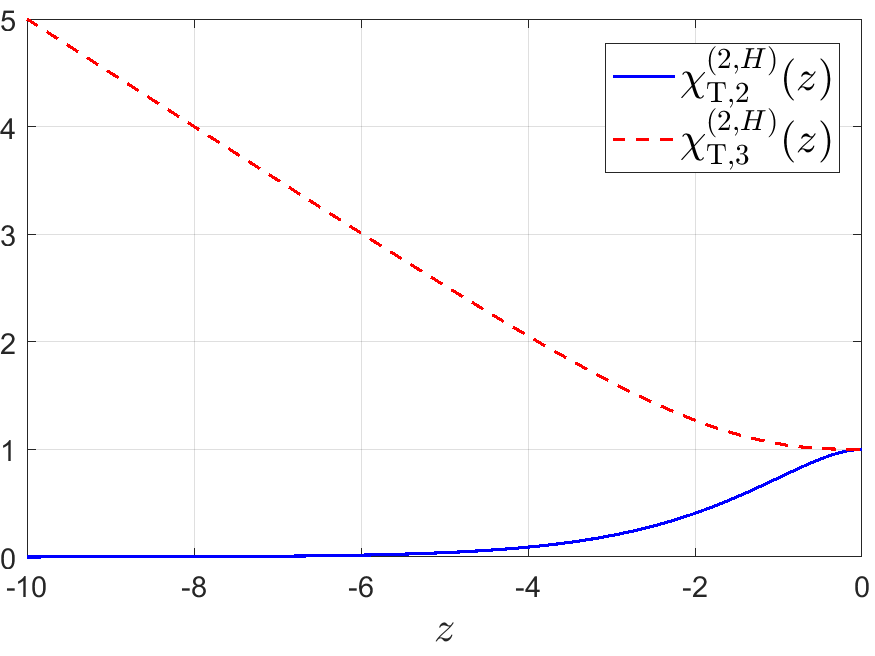}
  \label{subfig:chi-TIF2-Heun}}
\hspace{1cm}
  \subfigure[NIF2-Heun]{\includegraphics[width=2.4in]{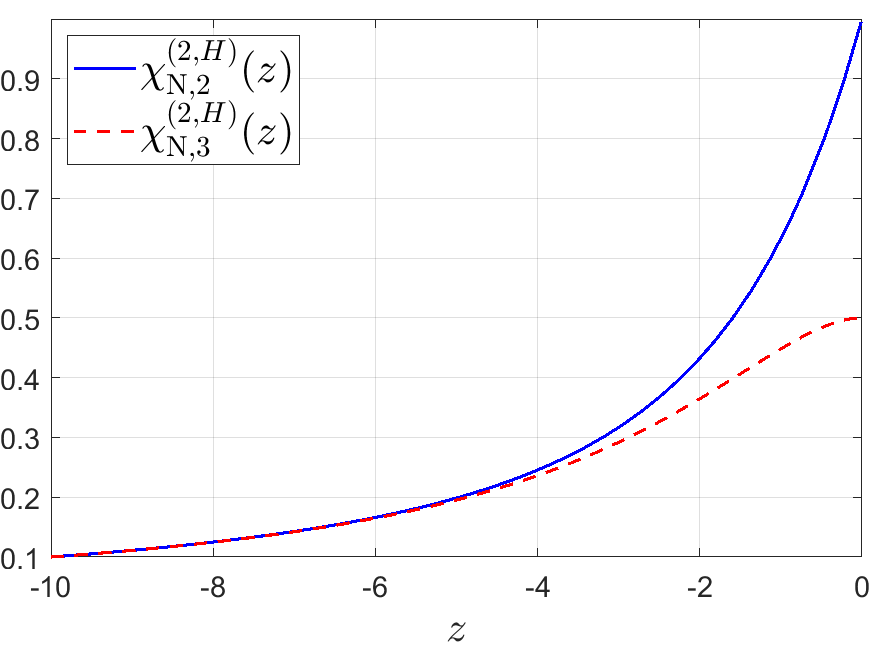}
  \label{subfig:chi-NIF2-Heun}}
  \caption{Correction coefficients of two corrected IF2-Heun schemes.}
  \label{fig: corrected IF2-Heun coefficients rate}
\end{figure}

To end this subsection, Figure \ref{fig: corrected IF2-Heun coefficients rate} depicts the correction coefficients in \eqref{def: Heun IF2 T-type coefficients} and \eqref{def: Heun IF2 N-type coefficients}. As seen from Figure \ref{fig: corrected IF2-Heun coefficients rate}, the correction coefficients for the TIF2-Heun \eqref{def: Heun TIF2 coefficients} and NIF2-Heun  \eqref{def:  Heun NIF2 coefficients} schemes are positive.

\subsection{Corrections of Ralston's IF2 method}

Now consider the $\mathrm{T}$-type and $\mathrm{N}$-type corrections for the second-order Ralston's IF2 scheme \cite{ZhangYanQianSong:2022CMAME} with the following Butcher tableau of the underlying explicit method \cite{Ralston:1962MCOM}
\begin{align*}
  \text{Second-order Ralston:}\quad \begin{array}{c|cc}
    0 & \\
    \frac{2}{3} & \frac{2}{3} \\[2pt]
    \hline\\[-2ex] & \frac{1}{4} & \frac{3}{4}
  \end{array}\;.
\end{align*}

(\textbf{$\mathrm{T}$-type correction}) By the definition \eqref{def: general T-type coefficients}, it is easy to get the telescopic correction coefficients
\begin{align}\label{def: Ralston IF2 T-type coefficients}
  \chi_{\cte,2}^{(2,R)}(z):=  e^{\frac{2 z}{3}} \brab{1-\tfrac{2}{3} z},\quad
  \chi_{\cte,3}^{(2,R)}(z):= \tfrac{1}{4} (4-z) e^z-\tfrac{3}{4} z e^{\frac{z}{3}}.
\end{align}
By the definition \eqref{eq: general T-type steady-state coefficients}, the resulting TIF2-Ralston scheme has
the following method coefficients
\begin{align}\label{def: Ralston TIF2 coefficients}
  \hat{a}^{\cte}_{21}(z)=\tfrac{2}{3-2z},\quad
  \hat{a}^{\cte}_{31}(z)=\tfrac{e^{\frac{2z}{3}}}{(4-z)e^{\frac{2z}{3}}-3z},\quad
  \hat{a}^{\cte}_{32}(z)=\tfrac{3}{(4-z)e^{\frac{2z}{3}}-3z}\,.
\end{align}
Then using \eqref{eq: general T-type differentiation matrix}, we have the associated differentiation matrix
\begin{align*}    
  D_{\cte}^{(2,R)}(z)=\begin{pmatrix}
 \frac{3}{2}-\frac{z}{2} & 0 \\
 \frac{5}{6} e^{\frac{2 z}{3}} & \frac{4-z}{3}e^{\frac{2 z}{3}}-\frac{z}{2}
  \end{pmatrix}\,.
\end{align*}
It is easy to get that $$\mathrm{Det}\kbrab{\mathcal{S}(D_{\cte,1}^{(2,R)};z)} = \tfrac{3-z}{2}\ge\tfrac{3}{2}\quad\text{ for $z\le0$.}$$ The second leading principal minor can be bounded by
\begin{align*}
  \mathrm{Det}\kbrab{\mathcal{S}(D_{\cte}^{(2,R)};z)} = e^{\frac{2 z}{3}}\brab{2-\tfrac{25}{144} e^{\frac{2 z}{3}}} + \brab{\tfrac{1}{6} e^{\frac{2 z}{3}}+\tfrac{1}{4}} z^2-\brab{\tfrac{7}{6} e^{\frac{2 z}{3}}+\tfrac{3}{4}} z\ge\tfrac{263}{144}
\end{align*}
for $z\le0$. It follows that the matrix $D_{\cte}^{(2,R)}(z)$ is positive definite for $z\le0$.

(\textbf{$\mathrm{N}$-type correction}) By the definition \eqref{def: general N-type coefficients},
we obtain the nonlinear-term correction coefficients
\begin{align}\label{def:  Ralston IF2 N-type coefficients}
  \chi_{\cnt,2}^{(2,R)}(z)= \tfrac{e^{\frac{2 z}{3}}-1}{z},\quad
  \chi_{\cnt,3}^{(2,R)}(z)= \tfrac{e^z-1}{z}-\tfrac{1}{4}e^z\,.
\end{align}
By the definition \eqref{eq: general N-type steady-state coefficients}, the resulting NIF2-Ralston scheme has
the method coefficients
\begin{align}\label{def: Ralston NIF2 coefficients}
  \hat{a}^{\cnt}_{21}(z)=\tfrac{e^{\frac{2 z}{3}}-1}{z},\quad
  \hat{a}^{\cnt}_{31}(z)=\tfrac{1}{4}e^z,\quad
  \hat{a}^{\cnt}_{32}(z)=\tfrac{e^z-1}{z}-\tfrac{1}{4}e^z\,.
\end{align}
One has the associated differentiation matrix
\begin{align*}    
  D_{\cnt}^{(2,R)}(z)=\begin{pmatrix}
 \frac{z}{e^{\frac{2 z}{3}}-1}+\frac{z}{2} & 0 \\[5pt]
 \frac{e^z [e^{\frac{2 z}{3}} (z-4)+4] z}{(e^{\frac{2 z}{3}}-1) [e^z (z-4)+4]} & \frac{[e^z (z-4)-4] z}{2 e^z (z-4)+8}
  \end{pmatrix}\,.
\end{align*}
To examine the positive definiteness of $D_{\cnt}^{(2,R)}(z)$, we need the following result.
\begin{proposition}\label{prop: auxiliary Ralston NIF2}
  For any $z\le0$, it holds that
  \begin{align*}
    g_{\cnt}^{(2,R)}(z):=-e^{2z}z^2+16(1-e^{\frac{2z}{3}}) (1-e^{2z}) +8(1-e^{\frac{2z}{3}}) e^{2z}\big[ z+2(e^{-\frac{4z}{3}}-1) \big]\ge0.
  \end{align*}
\end{proposition}
\begin{proof}
  Notice that 
 \begin{align*}
   1-e^{\frac{2z}{3}} \ge 1-\frac{1}{1-\frac{2z}{3}}=\frac{\frac{2z}{3}}{1-\frac{2z}{3}}, \quad 1-e^{2z} \ge 1-\frac{1}{1-2z}=\frac{2z}{1-2z}.
 \end{align*}
Also, the function $z+2(e^{-\frac{4z}{3}}-1) $ is decreasing for $z<0$. It is easy to get that
  \begin{align*}
    g_{\cnt}^{(2,R)}(z)\ge -e^{2z}z^2+16(1-e^{\frac{2z}{3}}) (1-e^{2z})\ge z^2\kbraB{-e^{2z}+\tfrac{64}{3(1-\frac{2z}{3}) (1-2z)}}\ge 0.
  \end{align*}
  It completes the proof.
\end{proof}
With the auxiliary function $g_{\cnt}^{(2,R)}(z)$ in Proposition \ref{prop: auxiliary Ralston NIF2}, one can check that the differentiation matrix $D_{\cnt}^{(2,R)}(z)$ is positive semi-definite for $z\le0$ due to the following facts: the first leading principal minor $\mathrm{Det}\kbrab{\mathcal{S}(D_{\cnt,1}^{(2,R)};z)}=\frac{z}{e^{\frac{2 z}{3}}-1}+\frac{z}{2}\ge\frac{3}{2}$ and the second leading principal minor
\begin{align*}
  \mathrm{Det}\kbrab{\mathcal{S}(D_{\cnt}^{(2,R)};z)}= \frac{z^2g_{\cnt}^{(2,R)}(z)}{4 (e^{\frac{2 z}{3}}-1)^2[e^z (z-4)+4]^2}\ge0\quad\text{for $z\le0$.}
\end{align*}

As a result, one has the following result by applying Lemma \ref{lemma: EERK energy stability}.
\begin{theorem}\label{thm: energy stability CIF2-Ralston}
  Assume that $g$ is Lipschitz-continuous with a constant $\ell_{g}>0$ and the stabilized parameter $\kappa$ in \eqref{def: stabilized parameter} is chosen properly large such that $\kappa\ge 2\ell_g$. In solving the gradient flow \eqref{problem: stabilized version}, the TIF2-Ralston \eqref{def: Ralston TIF2 coefficients} and NIF2-Ralston \eqref{def: Ralston NIF2 coefficients} schemes preserve the original energy dissipation law \eqref{eq: continuous energy law} at all stages.
\end{theorem}

\begin{figure}[htb!]
  \centering
  \subfigure[TIF2-Ralston]{\includegraphics[width=2.4in]{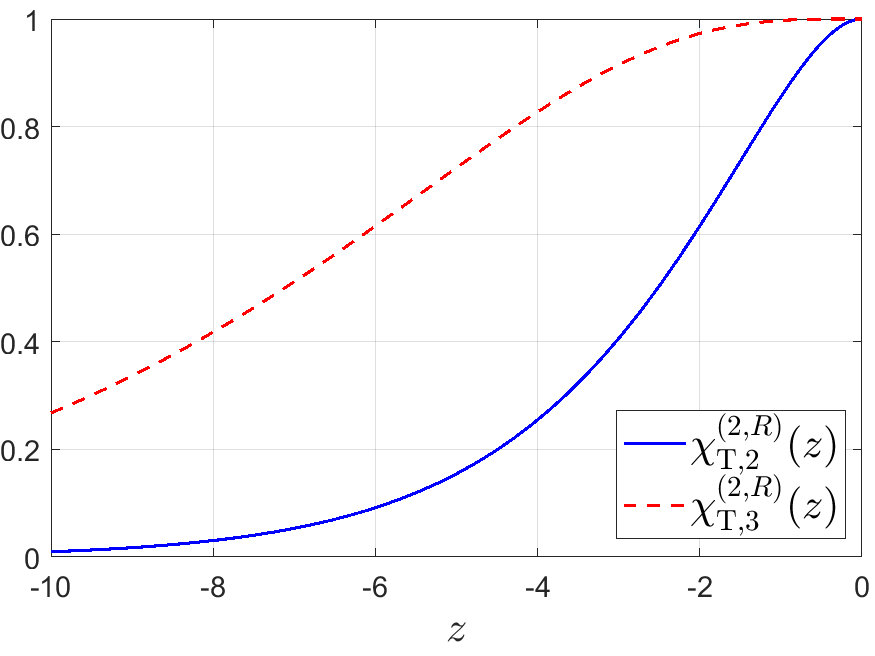}
  \label{subfig:chi-TIF2-Ralston}}
\hspace{1cm}
  \subfigure[NIF2-Ralston]{\includegraphics[width=2.4in]{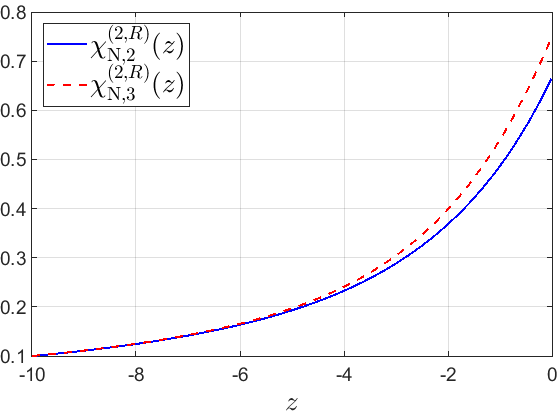}
  \label{subfig:chi-NIF2-Ralston}}
  \caption{Correction coefficients of two corrected IF2-Ralston schemes.}
  \label{fig: corrected IF2-Ralston coefficients rate}
\end{figure}

To end this subsection, we mention that the correction coefficients defined in \eqref{def: Ralston IF2 T-type coefficients} and \eqref{def: Ralston IF2 N-type coefficients} of two corrected IF2-Ralston schemes are positive and smaller than 1, cf. Figure \ref{fig: corrected IF2-Ralston coefficients rate}. The behaviors are quite similar to those of the correction coefficients $\chi_{\cte}^{(1)}(z)$ and $\chi_{\cnt}^{(1)}(z)$ for the corrected IF1 schemes, cf. Figure \ref{fig: corrected IF1 coefficients rate}. 

\subsection{Tests of corrected IF2 schemes}\label{subsec:test IF2}

This subsection uses Example \ref{example 1} in \lan{Subsection} \ref{subsec:test IF1} to examine the accuracy and numerical behaviors of the IF2-Heun and IF2-Ralston methods and the corresponding corrections, including the TIF2-Heun \eqref{def: Heun TIF2 coefficients}, NIF2-Heun \eqref{def:  Heun NIF2 coefficients}, TIF2-Ralston \eqref{def: Ralston TIF2 coefficients} and NIF2-Ralston \eqref{def: Ralston NIF2 coefficients} schemes.

\begin{figure}[htb!]
  \centering
  \subfigure[\lan{Numerical errors of corrected IF2-Heun}]{
    \includegraphics[width=2.4in]{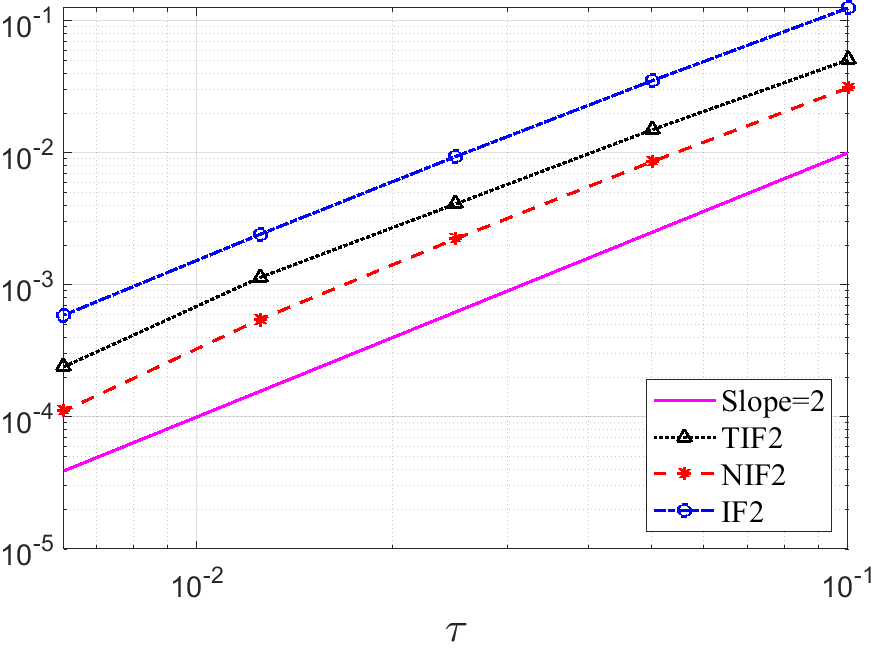}}
  \hspace{1cm}
  \subfigure[\lan{Numerical errors of corrected IF2-Ralston}]{
    \includegraphics[width=2.4in]{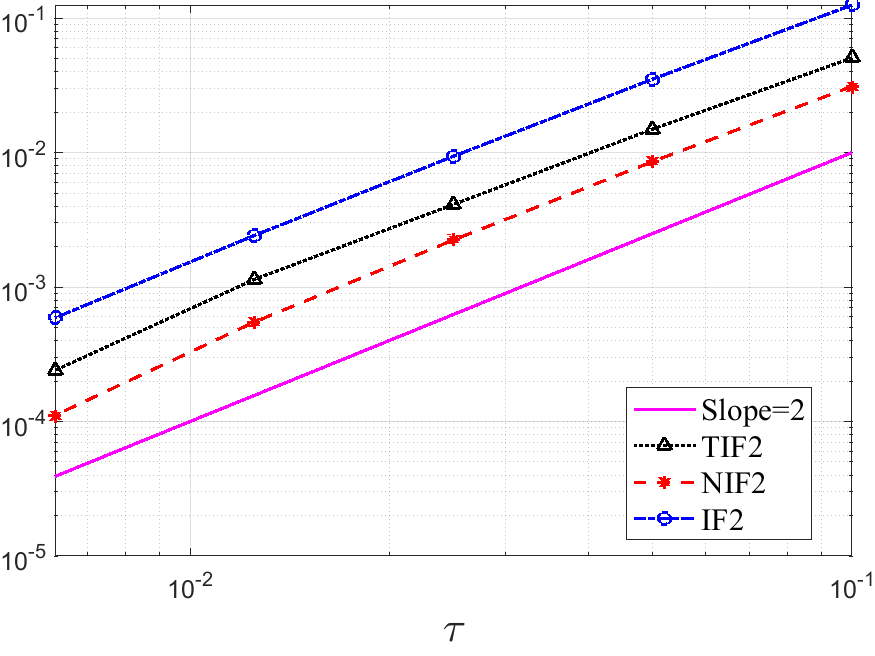}}
  \caption{\qing{Errors of two IF2 methods and their corrections.}}
  \label{fig: accuracy and IF2 and CIF2}
\end{figure}

First of all, we run the mentioned schemes for different time-step sizes $\tau=2^{-k}/10$ $(0\le k\le 4)$ up to $T=20$, with $h=1/500$ and the stabilized parameter $\kappa=4$ to test the temporal convergence. The solution errors are listed in Figure \ref{fig: accuracy and IF2 and CIF2} (a)-(b), where the reference solutions are computed by the NIF2-Heun and NIF2-Ralston schemes, respectively, with a small time-step $\tau=0.001$. As expected, these schemes are second-order accurate for small time-step sizes. Also, we observe that the T-type corrected schemes generate a bit less accurate solution than N-type schemes. The classic IF2 methods generate bigger errors than our corrected schemes.

\begin{figure}[htb!]
  \centering
  \subfigure[Final solution $u_h^N$]{
    \includegraphics[width=2.05in,height=1.5in]{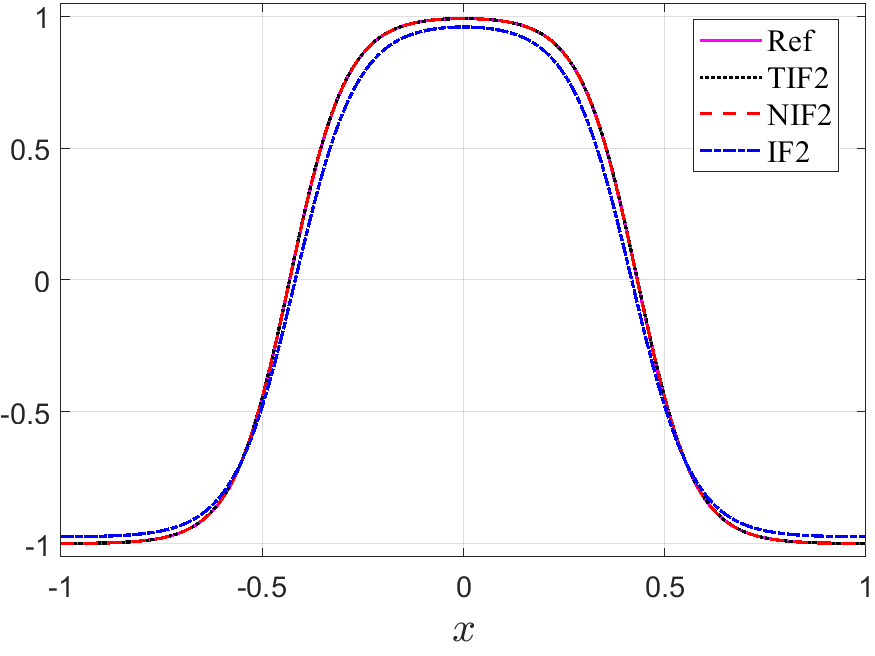}} 
  \subfigure[Discrete energy $E\kbrat{u_h^n}$]{
    \includegraphics[width=2.05in,height=1.53in]{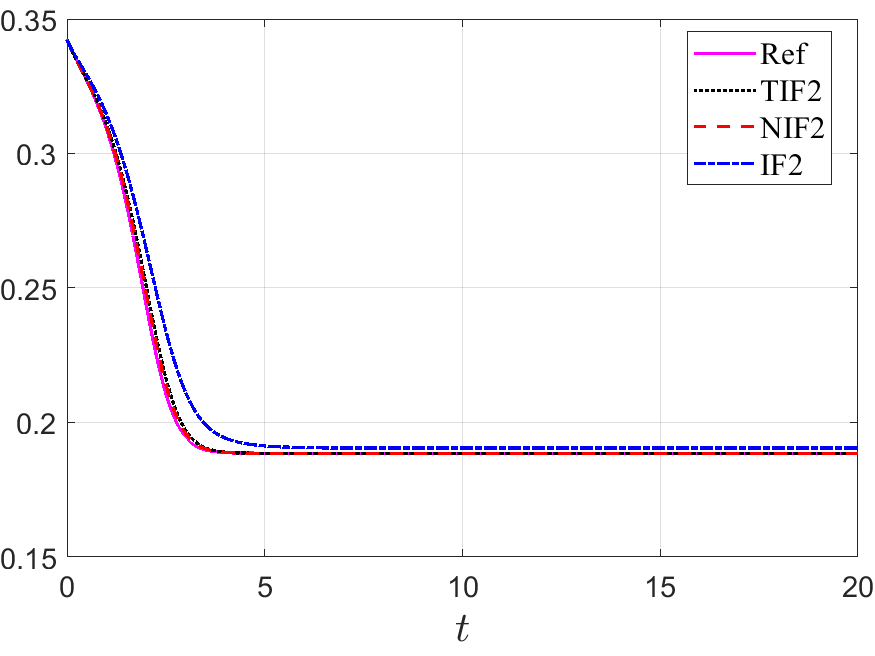}}
  \subfigure[Maximum norm $\|u_h^n\|_{\infty}$]{
    \includegraphics[width=2.05in,height=1.5in]{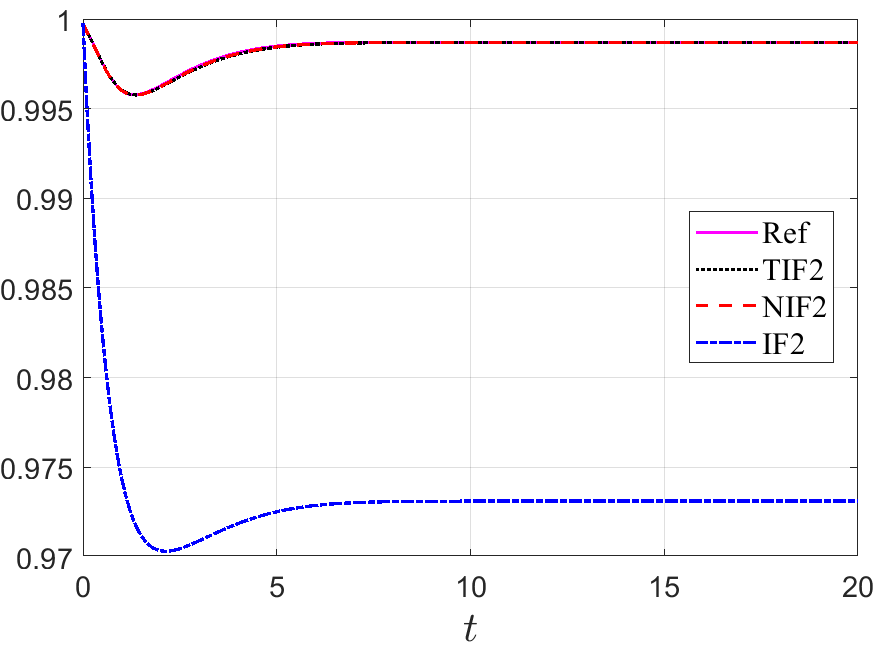}}
  \caption{Comparisons of IF2-Ralston and corrected IF2-Ralston schemes for $\tau=0.1$.}
  \label{fig: Compare IF2 and CIF2 with tau = 0.1}
\end{figure}
\begin{figure}[htb!]
  \centering
  \subfigure[Final solution $u_h^N$]{
    \includegraphics[width=2.05in,height=1.5in]{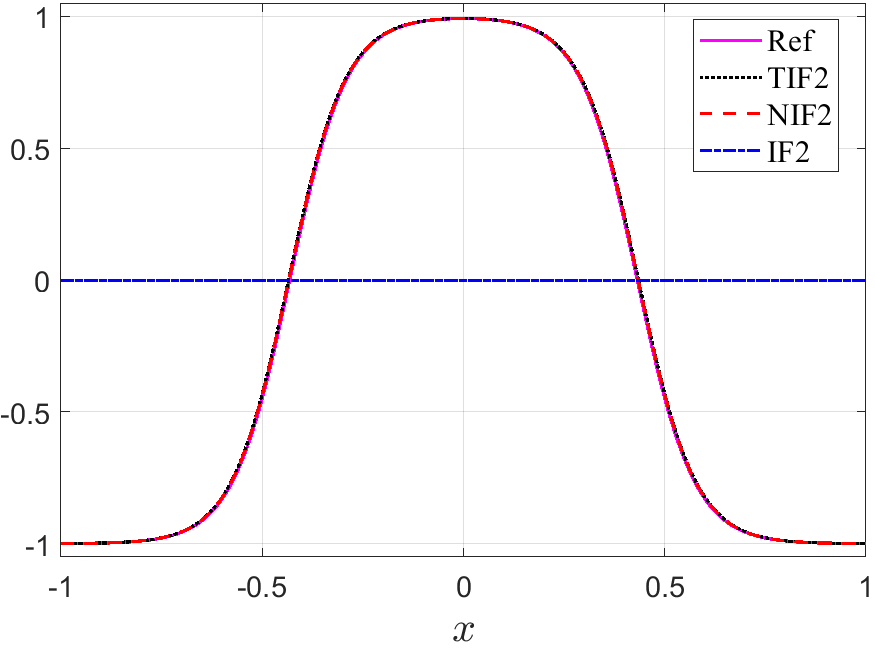}} 
  \subfigure[Discrete energy $E\kbrat{u_h^n}$]{
    \includegraphics[width=2.05in,height=1.53in]{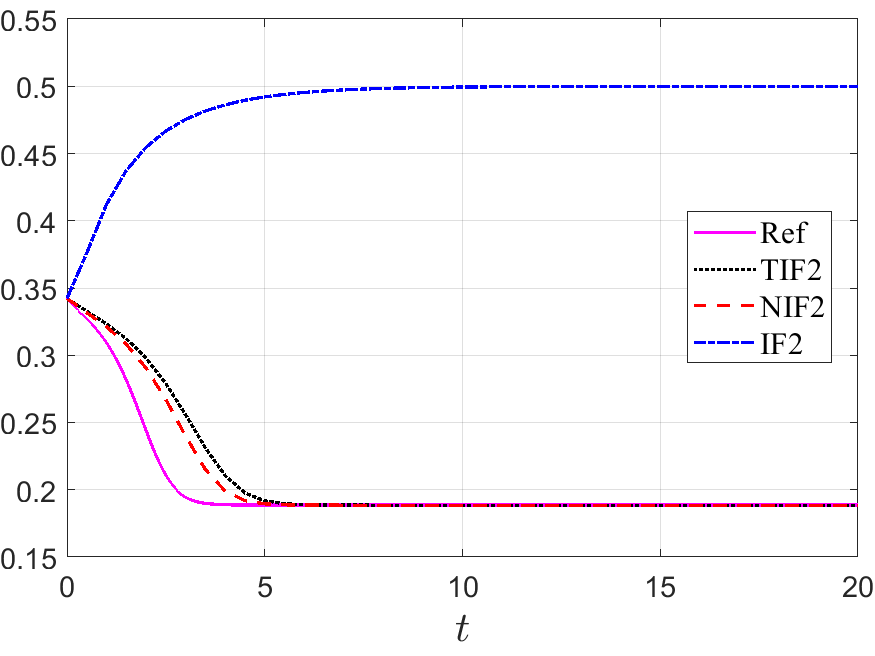}}
  \subfigure[Maximum norm $\|u_h^n\|_{\infty}$]{
    \includegraphics[width=2.05in,height=1.5in]{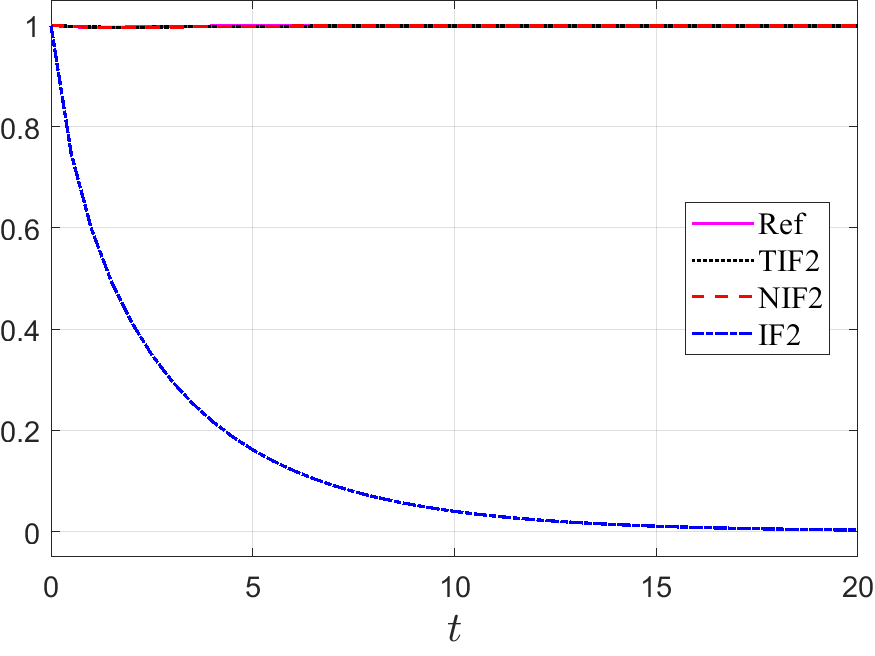}}
  \caption{Comparisons of IF2-Ralston and corrected IF2-Ralston schemes for $\tau=0.5$.}
  \label{fig: Compare IF2 and CIF2 with tau = 0.5}
\end{figure}

Next, taking the IF2-Ralston method for example, we run the three second-order schemes to the final time $T=20$ for two different time-steps $\tau=0.1$ and $0.5$, as shown in Figures \ref{fig: Compare IF2 and CIF2 with tau = 0.1}-\ref{fig: Compare IF2 and CIF2 with tau = 0.5}, respectively. As expected, two corrected IF2 methods maintain the steady-state well while the final solution $u_h^N$ of IF2-Ralston scheme gradually collapses into the trivial solution as the time-step $\tau$ increases to $0.5$, see Figures \ref{fig: Compare IF2 and CIF2 with tau = 0.1}(a)-\ref{fig: Compare IF2 and CIF2 with tau = 0.5}(a). For this example, we observe that the two corrected IF2-Ralston schemes 
\eqref{def: Ralston TIF2 coefficients} and \eqref{def: Ralston NIF2 coefficients} always perform well.

Interestingly, the two corrected IF2-Ralston methods also preserve the maximum bound principle similar to the IF2-Ralston method.
Maybe, this is closely related to the positivity (cf. Figure \ref{fig: corrected IF2-Ralston coefficients rate}) of the corresponding correction coefficients defined in \eqref{def: Ralston IF2 T-type coefficients} and \eqref{def: Ralston IF2 N-type coefficients}.
Note that, similar numerical behaviors are also observed for the IF2-Heun, \lan{TIF2-Heun \eqref{def: Heun TIF2 coefficients} and NIF2-Heun \eqref{def:  Heun NIF2 coefficients} schemes}, and we omit relevant presentations.

\section{Extensions to high-order IF methods}\label{sec:CIF3}
\setcounter{equation}{0}

As described at the beginning of Section \ref{sec:CIF2}, one can apply the $\mathrm{T}$-type and $\mathrm{N}$-type corrections to each stage of the $s$-stage stabilized IF method \eqref{scheme: IF stabilized} and obtain a total of $(s+1)!$ corrected schemes for a specific $s$-stage IF method. At least, one has 24 corrected schemes for a 3-stage IF method and 120 corrections for a 4-stage IF method. 
As done in Section \ref{sec:CIF2}, we only modify the $s$-stage IF method by utilizing the same correction strategy at each stage. Specifically, for the nonlinear-term translation ($\mathrm{N}$-type) correction, we always modify the last nonlinear-term $g_{\kappa}(U^{n,i})$ at the stage $t_{n,i+1}$. 

\subsection{Corrections of Heun's IF3 method}
At first, we examine the $\mathrm{T}$-type and $\mathrm{N}$-type corrections for the third-order Heun's IF3 scheme \cite{LiLiJuFeng:2021SISC,ZhangYanQianSong:2022CMAME} with the following Butcher tableau of the underlying explicit Runge-Kutta method \cite{Heun:1900ZMP}
\begin{align*}
  \text{Third-order Heun:}\quad \begin{array}{c|ccc}
    0 & \\
    \frac{1}{3} & \frac{1}{3} \\[2pt]
    \frac{2}{3} & 0 & \frac{2}{3} \\[2pt]
    \hline\\[-2ex] & \frac{1}{4} & 0 & \frac{3}{4}
  \end{array}\,.
\end{align*}

By following the theoretical framework in Section \ref{sec:unified framework}, one can \lan{find that}
the TIF3-Heun and NIF3-Heun schemes have the following coefficient matrices, respectively, 
\begin{align}
  \widehat{A}_{\cte}^{(3,H)}(z) =&\;\begin{pmatrix}
 \frac{1}{3-z} \\[2pt]
 0 & \frac{2}{3 e^{\frac{z}{3}}-2 z}  \\[2pt]
 \frac{e^{\frac{2 z}{3}}}{e^{\frac{2 z}{3}} (4-z)-3 z} & 0 & \frac{3}{e^{\frac{2 z}{3}} (4-z)-3 z}
  \end{pmatrix}\,,\label{def: Heun TIF3 coefficients} \\
  \widehat{A}_{\cnt}^{(3,H)}(z)=&\; \begin{pmatrix}
 \frac{e^{\frac{z}{3}}-1}{z} \\[2pt]
 0 & \frac{e^{\frac{2 z}{3}}-1}{z} \\[2pt]
 \frac{e^z}{4} & 0 & \frac{e^z-1}{z}-\frac{e^z}{4}
  \end{pmatrix}\,.\label{def: Heun NIF3 coefficients}
\end{align}

\ref{app:IF3-Heun} shows that  the associated differentiation matrices of the TIF3-Heun \eqref{def: Heun TIF3 coefficients} and NIF3-Heun \eqref{def: Heun NIF3 coefficients} are positive (semi-)definite for $z\le0$. Then we obtain the following theorem by using Lemma \ref{lemma: EERK energy stability}.
\begin{theorem}\label{thm: energy stability CIF3-Heun}
  Assume that $g$ is Lipschitz-continuous with a constant $\ell_{g}>0$ and the stabilized parameter $\kappa$ in \eqref{def: stabilized parameter} is chosen properly large such that $\kappa\ge 2\ell_g$. In solving the gradient flow \eqref{problem: stabilized version}, the TIF3-Heun \eqref{def: Heun TIF3 coefficients} and NIF3-Heun \eqref{def: Heun NIF3 coefficients} schemes preserve the original energy dissipation law \eqref{eq: continuous energy law} at all stages.
\end{theorem}

\begin{figure}[htb!]
  \centering
  \subfigure[TIF3-Heun]{\includegraphics[width=2.4in]{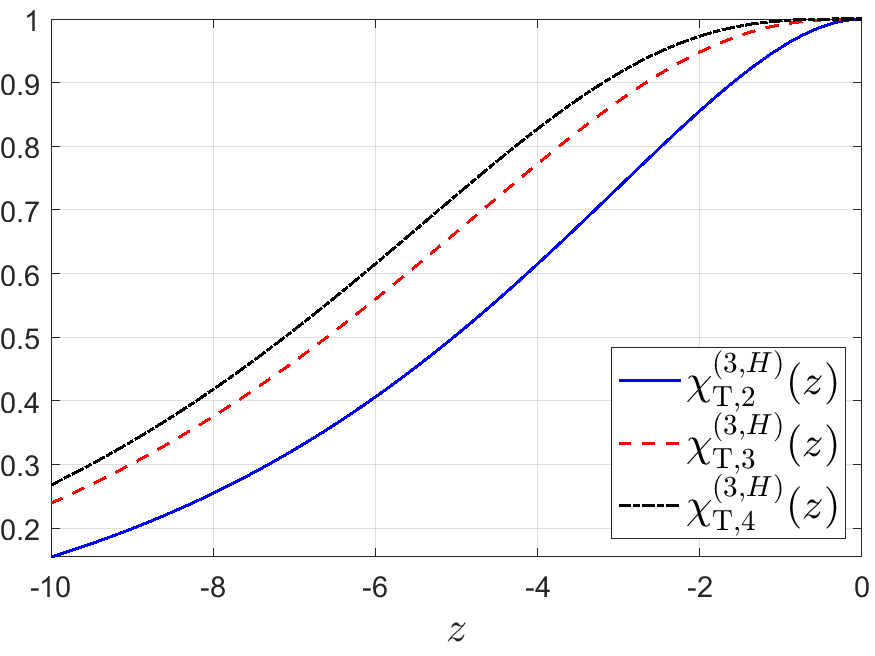}
  \label{subfig:chi-TIF3-Heun}}
\hspace{1cm}
  \subfigure[NIF3-Heun]{\includegraphics[width=2.4in]{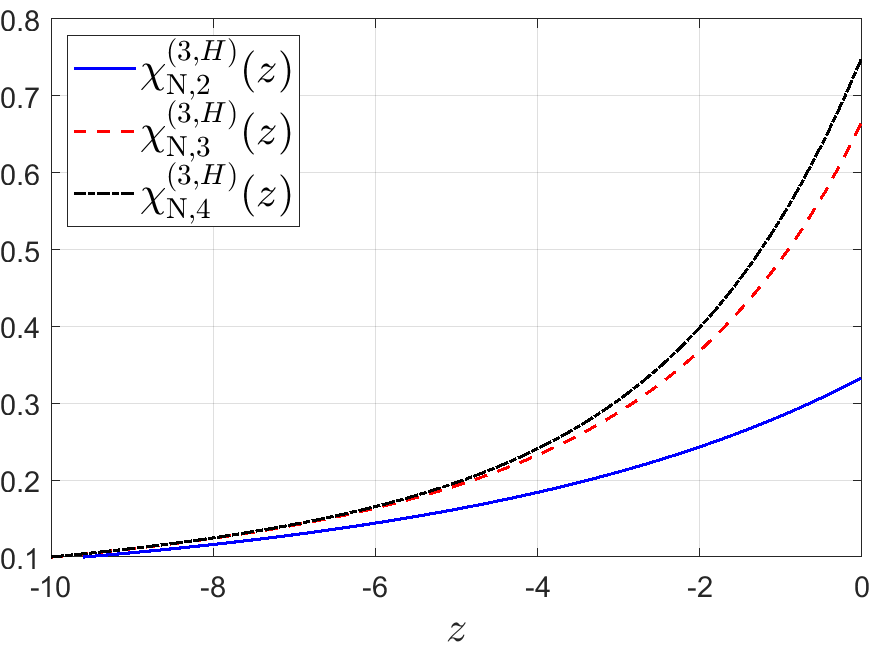}
  \label{subfig:chi-NIF3-Heun}}
  \caption{Correction coefficients of two corrected IF3-Heun schemes.}
  \label{fig: corrected IF3-Heun coefficients rate}
\end{figure}

Also, we mention that all the correction coefficients at three stages of two corrected IF3-Heun schemes are positive and smaller than 1, cf. Figure \ref{fig: corrected IF3-Heun coefficients rate}. 

\subsection{Corrections of Ralston's IF3 method}
As the second example, we consider the $\mathrm{T}$-type and $\mathrm{N}$-type corrections for the Ralston's IF3 scheme with the following Butcher tableau of the underlying explicit Ralston's method \cite{Ralston:1962MCOM}
\begin{align*}
  \text{Third-order Ralston:}\quad \begin{array}{c|ccc}
    0 & \\
    \frac{1}{2} & \frac{1}{2} \\[2pt]
    \frac{3}{4} & 0 & \frac{3}{4} \\[2pt]
    \hline\\[-2ex] & \frac{2}{9} & \frac{1}{3} & \frac{4}{9}
  \end{array}\,.
\end{align*}
The resulting TIF3-Ralston and NIF3-Ralston schemes have the following coefficient matrices, respectively, 
\begin{align}
  \widehat{A}_{\cte}^{(3,R)}(z) =&\;  \begin{pmatrix}
 \frac{1}{2-z}  \\[2pt]
 0 & \frac{3}{4 e^{\frac{z}{2}}-3 z}  \\[2pt]
 \frac{2 e^{\frac{3 z}{4}}}{e^{\frac{3 z}{4}} (9-2 z)-3 e^{\frac{z}{4}} z-4 z} & \frac{3 e^{\frac{z}{4}}}{e^{\frac{3 z}{4}} (9-2 z)-3 e^{\frac{z}{4}} z-4 z} & \frac{4}{e^{\frac{3 z}{4}} (9-2 z)-3 e^{\frac{z}{4}} z-4 z}
  \end{pmatrix}\,,\label{def: Ralston TIF3 coefficients} \\
  \widehat{A}_{\cnt}^{(3,R)}(z)=&\; \begin{pmatrix}
 \frac{e^{\frac{z}{2}}-1}{z}  \\[2pt]
 0 & \frac{e^{\frac{3 z}{4}}-1}{z}  \\[2pt]
 \frac{2 e^z}{9} & \frac{e^{\frac{z}{2}}}{3} & \frac{e^z-1}{z}-\frac{e^{\frac{z}{2}}}{3}-\frac{2 e^z}{9}
  \end{pmatrix}\,.\label{def: Ralston NIF3 coefficients}
\end{align}

\begin{figure}[htb!]
     \centering
     \subfigure[TIF3-Ralston]{\includegraphics[width=2.4in]{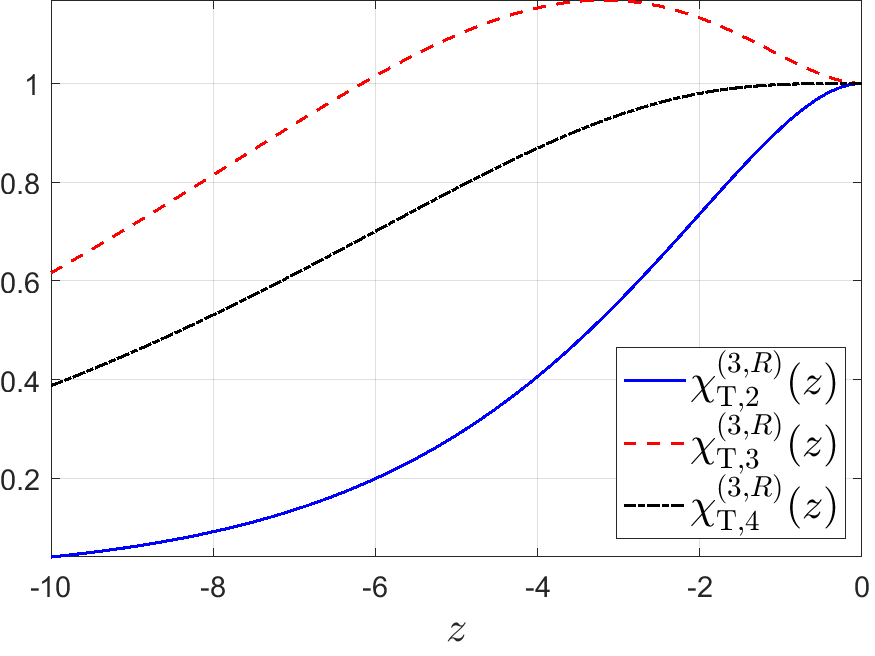}
          \label{subfig:chi-TIF3-Ralston}}
     \hspace{1cm}
     \subfigure[NIF3-Ralston]{\includegraphics[width=2.4in]{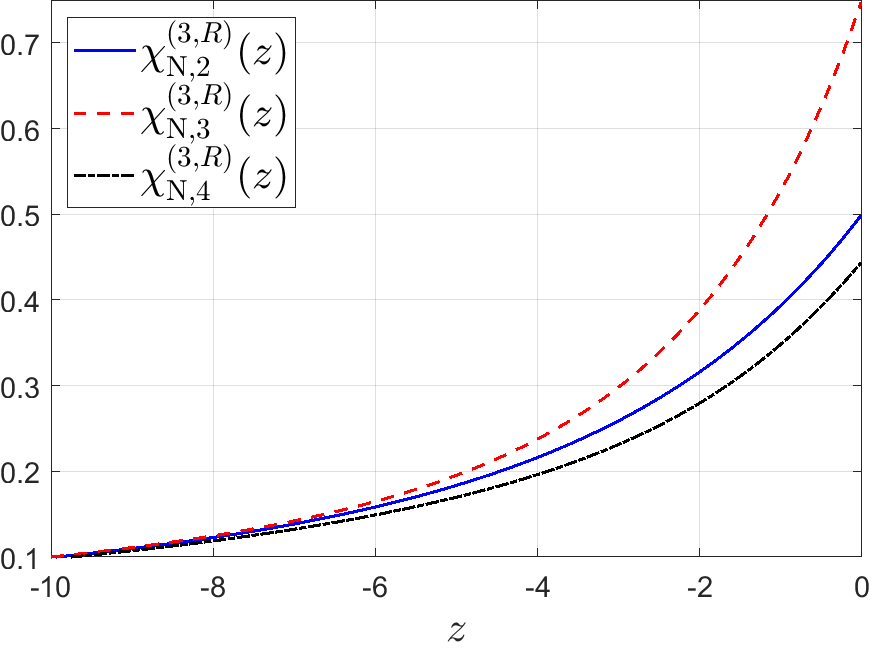}
          \label{subfig:chi-NIF3-Ralston}}
     \caption{Correction coefficients of two corrected IF3-Ralston schemes.}
     \label{fig: corrected IF3-Ralston coefficients rate}
\end{figure}

In \ref{app:IF3-Ralston}, we show that  the associated differentiation matrices of the TIF3-Ralston \eqref{def: Ralston TIF3 coefficients} and NIF3-Ralston \eqref{def: Ralston NIF3 coefficients} are positive (semi-)definite for $z\le0$. Thus one has the following result.
\begin{theorem}\label{thm: energy stability CIF3-Ralston}
 Assume that $g$ is Lipschitz-continuous with a constant $\ell_{g}>0$ and the stabilized parameter $\kappa$ in \eqref{def: stabilized parameter} is chosen properly large such that $\kappa\ge 2\ell_g$. In solving the gradient flow \eqref{problem: stabilized version}, the TIF3-Ralston \eqref{def: Ralston TIF3 coefficients} and NIF3-Ralston \eqref{def: Ralston NIF3 coefficients} schemes preserve the energy dissipation law \eqref{eq: continuous energy law} at all stages.
\end{theorem}

Figure \ref{fig: corrected IF3-Ralston coefficients rate} depicts the correction coefficients for the two corrected IF3-Ralston schemes. We see from Figure \ref{fig: corrected IF3-Ralston coefficients rate}(b) that the correction coefficients for the NIF3-Ralston scheme \eqref{def: Ralston NIF2 coefficients} are positive and smaller than 1. Most of the correction coefficients for the TIF3-Ralston \eqref{def: Ralston TIF2 coefficients} scheme are positive and smaller than 1, while the 3-stage correction coefficients $\chi_{\cte,3}^{(3,R)}>1$ near $z=-4$, cf. Figure \ref{fig: corrected IF3-Ralston coefficients rate}(a).

\subsection{Corrections of an IF4 method}\label{sec:CIF4}

To end this section, we consider the $\mathrm{T}$-type and $\mathrm{N}$-type corrections for Kutta's IF4 scheme \cite{ZhangYanQianSong:2022CMAME} with the following Butcher tableau of the underlying explicit Runge-Kutta method \cite{Kutta:1901ZMP}
\begin{align*}
  \text{Fourth-order Kutta:} \quad
  \begin{array}{c|cccc}
    0 & \\
    \frac{1}{2} & \frac{1}{2} \\[2pt]
    \frac{1}{2} & 0 & \frac{1}{2} \\[2pt]
    1 & 0 & 0 & 1 \\
    \hline\\[-2ex] & \frac{1}{6} & \frac{1}{3} & \frac{1}{3} & \frac{1}{6}
  \end{array}\,.
\end{align*}

\begin{figure}[htb!]
  \centering
  \subfigure[TIF4-Kutta]{\includegraphics[width=2.4in]{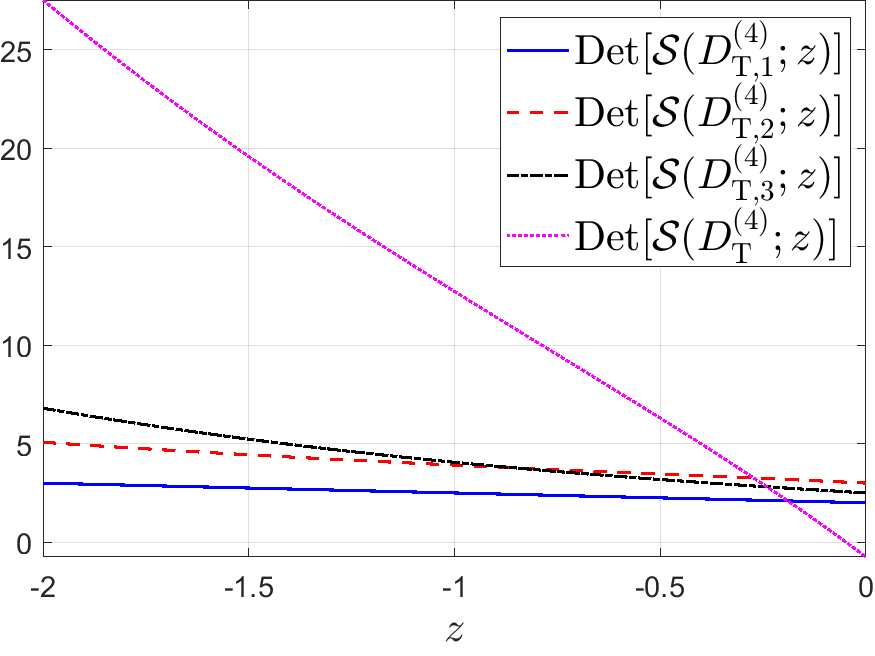}
    \label{subfig:chi-TIF4-Kutta}}
  \hspace{1cm}
  \subfigure[NIF4-Kutta]{\includegraphics[width=2.4in]{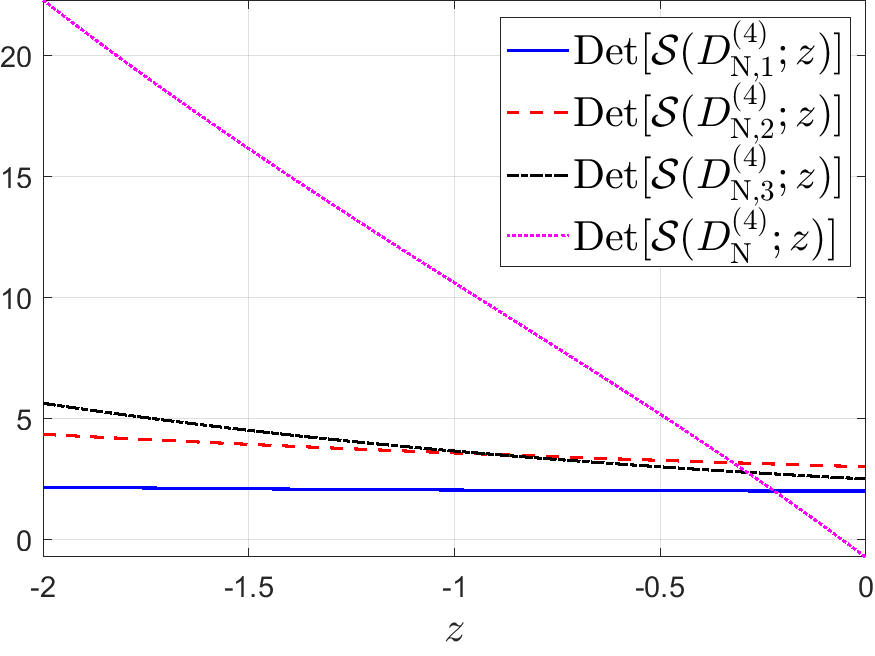}
    \label{subfig:chi-NIF4-Kutta}}
  \caption{Principal minors of differentiation matrices for corrected IF4-Kutta schemes.}
     \label{fig: corrected IF4-Kutta determinants}
\end{figure}

The $\mathrm{T}$-type and $\mathrm{N}$-type corrections will arrive at the TIF4-Kutta and NIF4-Kutta schemes having the following coefficient matrices, respectively,
\begin{align}
  \widehat{A}_{\cte}^{(4)}(z) =&\; \begin{pmatrix}
   \frac{1}{2-z}  \\
   0 & \frac{1}{2 e^{\frac{z}{2}}-z}  \\
   0 & 0 & \frac{1}{e^{\frac{z}{2}}-z}  \\
   \frac{e^z}{e^z (6-z)-4 e^{\frac{z}{2}} z-z} & \frac{2 e^{\frac{z}{2}}}{e^z (6-z)-4 e^{\frac{z}{2}} z-z} & \frac{2 e^{\frac{z}{2}}}{e^z (6-z)-4 e^{\frac{z}{2}} z-z} & \frac{1}{e^z (6-z)-4 e^{\frac{z}{2}} z-z}
  \end{pmatrix}\,,\label{def: Kutta TIF4 coefficients} \\
  \widehat{A}_{\cnt}^{(4)}(z) =&\; \begin{pmatrix}
   \frac{e^{\frac{z}{2}}-1}{z}  \\
   0 & \frac{e^{\frac{z}{2}}-1}{z}  \\
   0 & 0 & \frac{e^z-1}{z}  \\
   \frac{e^z}{6} & \frac{e^{\frac{z}{2}}}{3} & \frac{e^{\frac{z}{2}}}{3} & \frac{e^z (6-z)-4 e^{\frac{z}{2}} z-6}{6 z}
  \end{pmatrix}\,.\label{def: Kutta NIF4 coefficients}
\end{align}
Unfortunately, these differentiation matrices $D_{\cte}^{(4)}(z)$ and $D_{\cnt}^{(4)}(z)$  are not always positive definite for $z\le 0$ (the determinants $\mathrm{Det}[\mathcal{S}(D_{\cte}^{(4)};z)]$ and $\mathrm{Det}[\mathcal{S}(D_{\cnt}^{(4)};z)]$ are negative near $z=0$), cf. Figure \ref{fig: corrected IF4-Kutta determinants}, in which the leading principal minors of their symmetrical parts  are depicted. The above two corrected IF4-Kutta schemes \eqref{def: Kutta TIF4 coefficients} and \eqref{def: Kutta NIF4 coefficients} are undesirable for our aim. We believe that there must be some fourth-order correction schemes preserving the original energy dissipation law \eqref{eq: continuous energy law} from 120 possible corrections; however, we have not yet found them up to now.

\section{Numerical experiments}\label{sec:simulations}
\setcounter{equation}{0}

\subsection{Tests of corrected IF3 schemes}\label{subsec:test IF3}

\begin{figure}[htb!]
  \centering
  \subfigure[\lan{Numerical errors of corrected IF3-Heun}]{
    \includegraphics[width=2.4in]{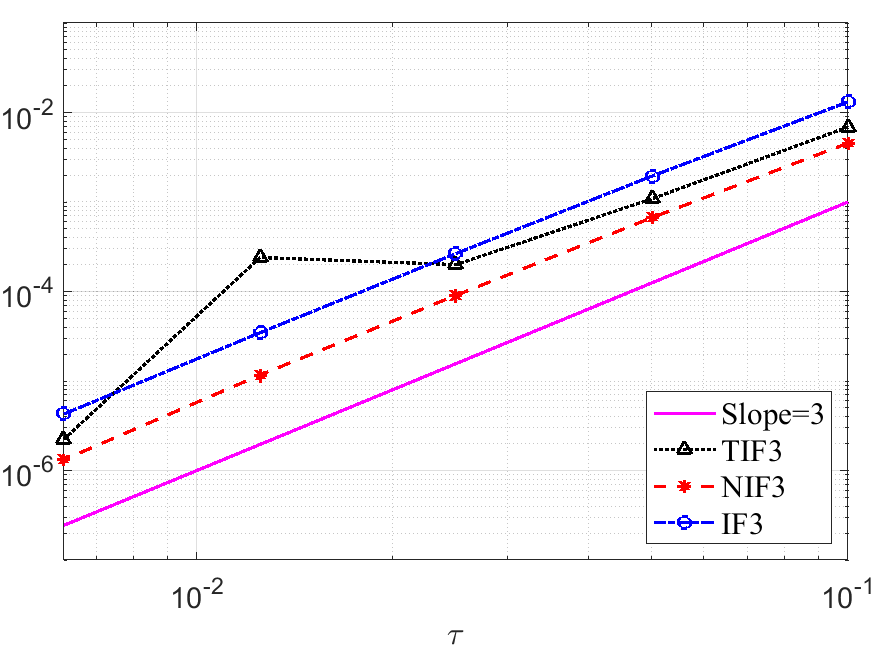}
    \label{subfig:IF3-Heun convergence}}
\hspace{1cm}
  \subfigure[\lan{Numerical errors of corrected IF3-Ralston}]{
    \includegraphics[width=2.4in]{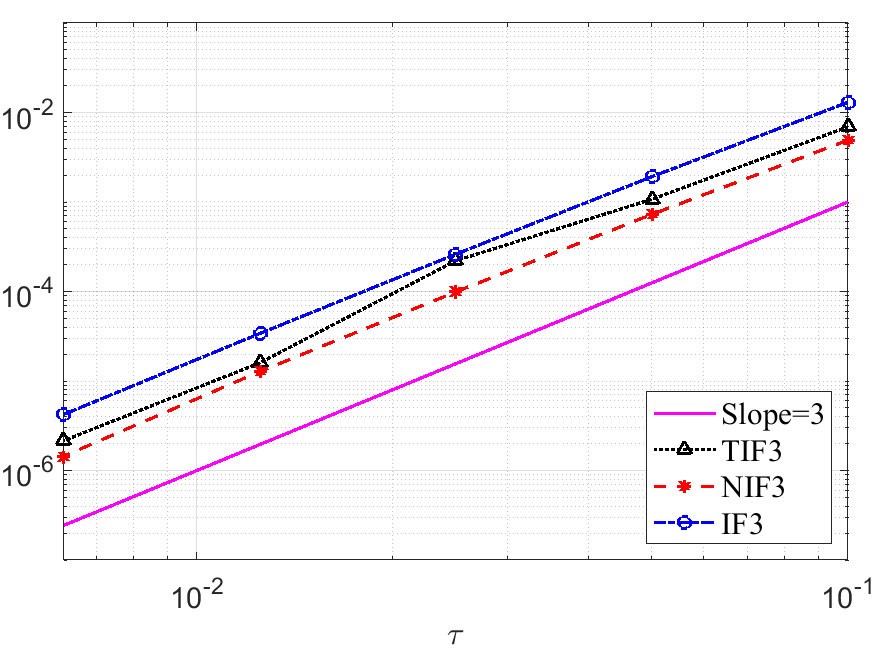}
    \label{subfig:IF3-Ralston convergence}}
  \caption{\qing{Errors of two IF3 methods and their corrections.}}
  \label{fig: accuracy and IF3 and CIF3}
\end{figure}

\lan{First of all, we use Example \ref{example 1} in Subsection \ref{subsec:test IF1} to examine the accuracy of the IF3-Heun and IF3-Ralston methods and the corresponding corrections, including the TIF3-Heun \eqref{def: Heun TIF3 coefficients}, NIF3-Heun \eqref{def:  Heun NIF3 coefficients}, TIF3-Ralston \eqref{def: Ralston TIF3 coefficients} and NIF3-Ralston \eqref{def: Ralston NIF3 coefficients} schemes. We run the mentioned six schemes for different time-step sizes $\tau=2^{-k}/10$ $(0\le k\le 4)$ up to $T=20$, with $h=1/500$ and the stabilized parameter $\kappa=4$. The solution errors are listed in Figure \ref{fig: accuracy and IF3 and CIF3} with the reference solution computed by a small time-step $\tau=0.01$ from the NIF3-Heun and NIF3-Ralston schemes, respectively. As expected, these schemes are third-order accurate. Also, we observe that the T-type corrected schemes generate a bit less accurate solution than N-type corrections. The classic IF3 methods generate bigger errors than our corrected schemes.}

\begin{example}\label{example 3} Consider the Allen-Cahn model $\partial_tu=\epsilon^2\partial_{xx} u-u^3+u$ with the interface parameter $\epsilon=0.1$ on $\Omega=(0,2\pi)$  subject to the following initial data from Chebfun (Matlab package), see \cite[guide19]{DriscollHaleTrefethen:2014Chebfun}, 
  $$u_0=\frac13\tanh(2\sin x) - e^{-23.5(x-\frac{\pi}{2})^2}+e^{-27(x-4.2)^2} + e^{-38(x-5.4)^2}.$$
  Always, we use the second-order center difference approximation with the spacing $h=\pi/320$. \lan{The solution preserves the maximum bound principle, that is, $\mynorm{u_h(t)}_{\infty}\le1$ since $\mynorm{u_h^0}_{\infty}\le1$, cf. \cite[Section 2]{DuJuLiQiao:2021SIREV}}.
\end{example}
\begin{figure}[htb!]
  \centering
  \subfigure[Final solution $u_h^N$]{
    \includegraphics[width=2.05in,height=1.5in]{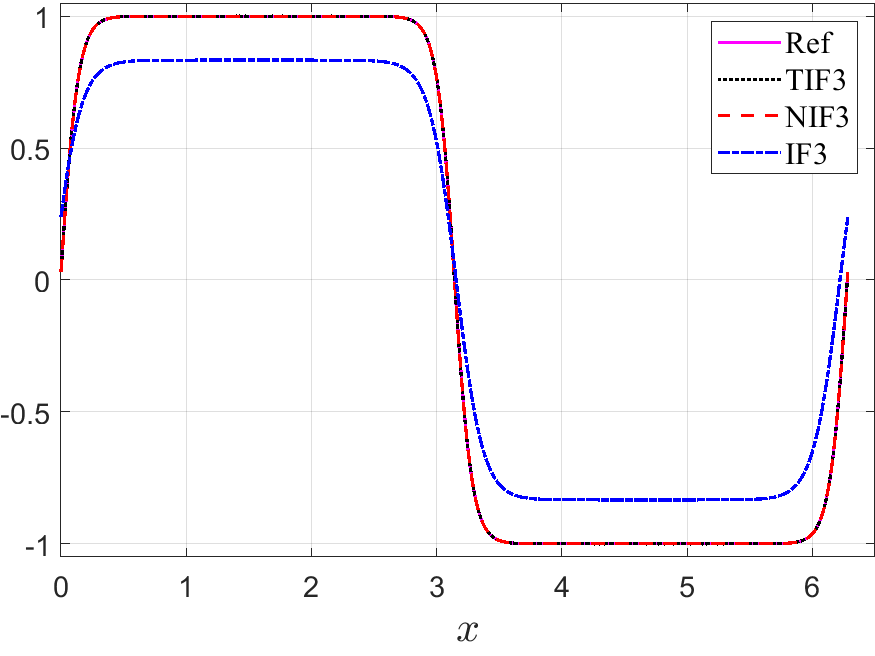}} 
  \subfigure[Discrete energy $E\kbrat{u_h^n}$]{
    \includegraphics[width=2.05in,height=1.53in]{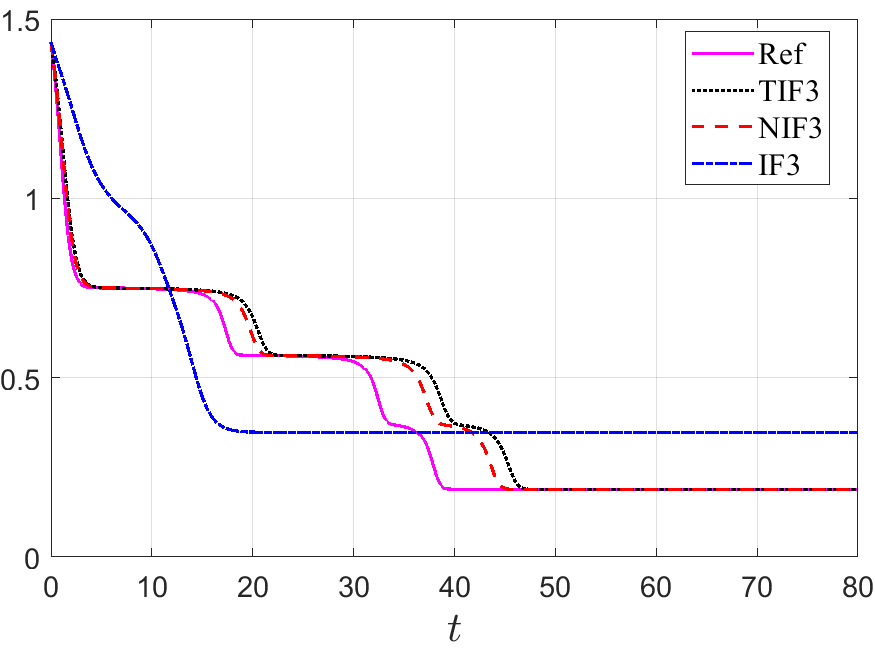}}
  \subfigure[Maximum norm $\|u_h^n\|_{\infty}$]{
    \includegraphics[width=2.05in,height=1.5in]{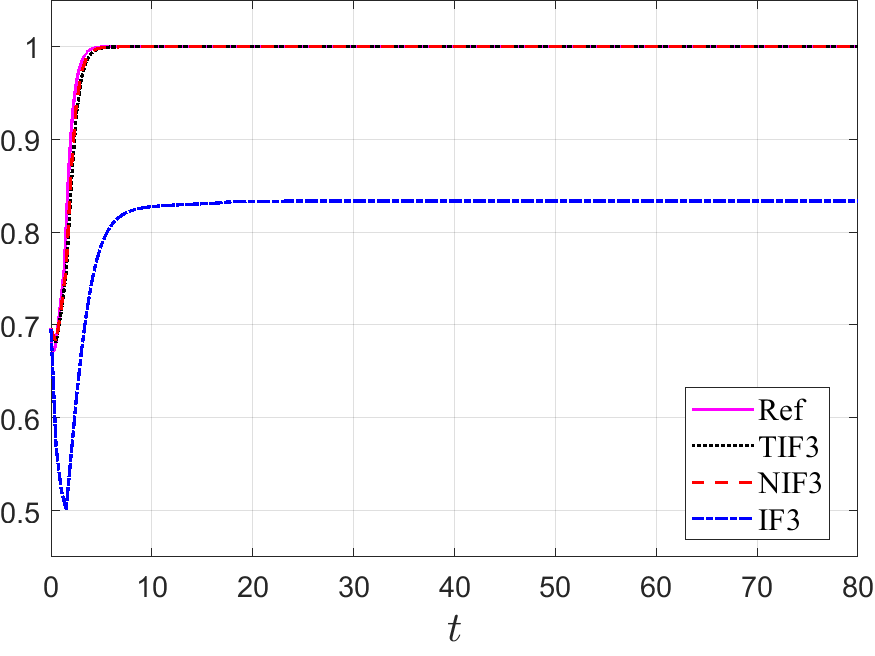}}
  \caption{Comparisons of IF3-Heun and corrected IF3-Heun schemes for $\tau=0.5$.}
  \label{fig: Compare Heun IF3 and CIF3 with tau = 0.5}
\end{figure}
\begin{figure}[htb!]
  \centering
  \subfigure[Final solution $u_h^N$]{
    \includegraphics[width=2.05in,height=1.5in]{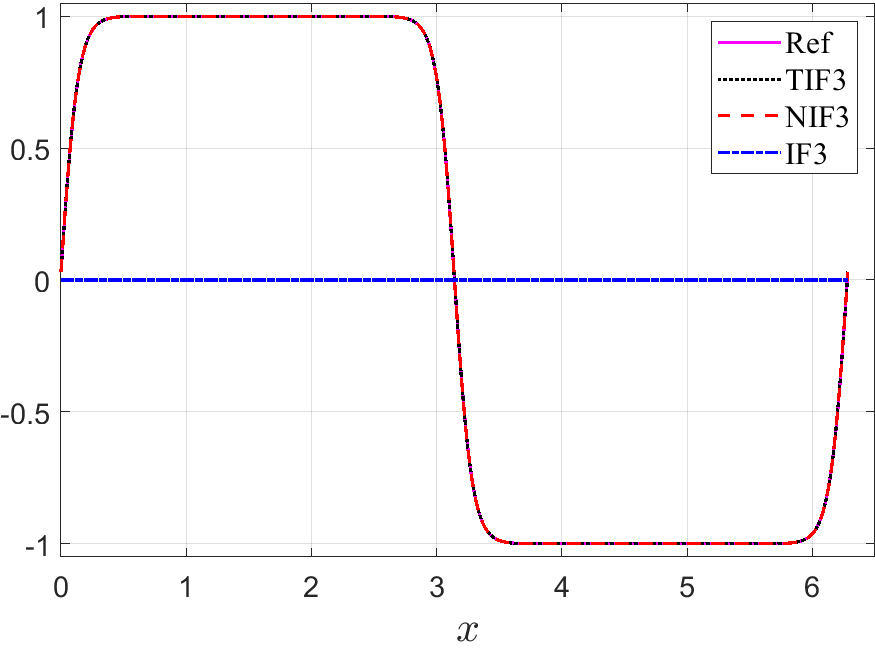}} 
  \subfigure[Discrete energy $E\kbrat{u_h^n}$]{
    \includegraphics[width=2.05in,height=1.53in]{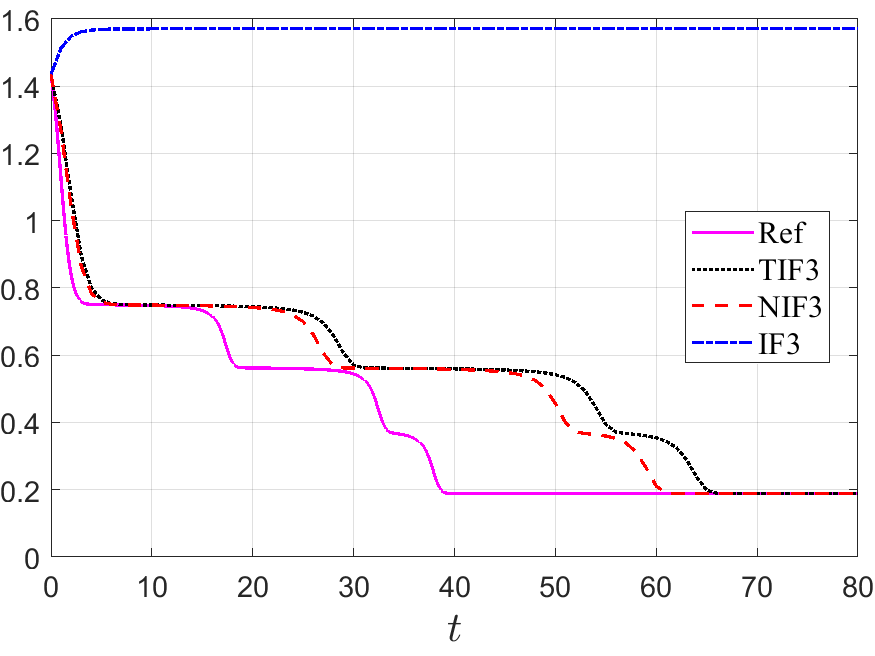}}
  \subfigure[Maximum norm $\|u_h^n\|_{\infty}$]{
    \includegraphics[width=2.05in,height=1.5in]{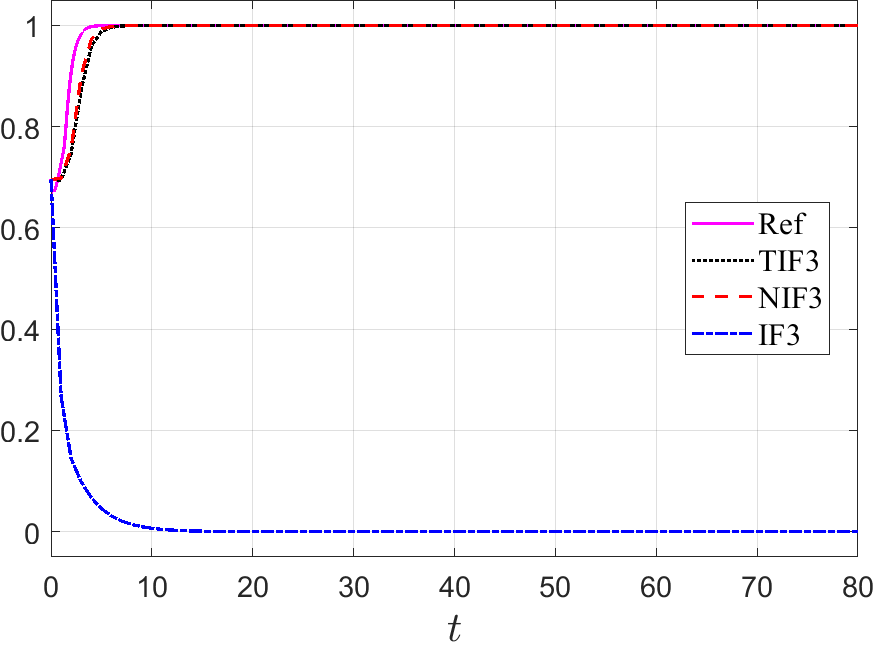}}
  \caption{Comparisons of IF3-Heun and corrected IF3-Heun schemes for $\tau=1$.}
  \label{fig: Compare Heun IF3 and CIF3 with tau = 1}
\end{figure}

We run the IF3-Heun, TIF3-Heun \eqref{def: Heun TIF3 coefficients} and NIF3-Heun \eqref{def: Heun NIF3 coefficients} schemes up to the final time $T=80$ with the stabilized parameter $\kappa=4$. We depict the final solution $u_h^N$, the discrete energy $E[u_h^n]$ and the maximum norm $\|u_h^n\|_{\infty}$ in Figures \ref{fig: Compare Heun IF3 and CIF3 with tau = 0.5} and \ref{fig: Compare Heun IF3 and CIF3 with tau = 1} for two different time-steps $\tau=0.5$ and $1$, respectively. The numerical solution (energy) of the NIF3-Heun scheme computed by $\tau=0.01$ is taken as the reference solution (energy). As seen, the two corrected IF3-Heun methods maintain the steady-state well while the final solution $u_h^N$ of IF3-Heun scheme gradually collapses into the metastable state solution as the parameter $\tau$ increases to $1$, see Figures \ref{fig: Compare Heun IF3 and CIF3 with tau = 0.5}(a)-\ref{fig: Compare Heun IF3 and CIF3 with tau = 1}(a). We notice that there are obvious differences in the energy curves of two corrected IF3-Heun schemes for $\tau=1$ although they are decreasing  
over the time and approach the same steady-state.

Also, we \lan{note} that the two corrected IF3-Heun methods preserve the maximum bound principle. 
This property may be closely related to the positivity of the corresponding correction coefficients, cf. Figure \ref{fig: corrected IF3-Heun coefficients rate}.

It is to mention that, under the same parameter settings, the IF3-Ralston, TIF3-Ralston \eqref{def: Ralston TIF3 coefficients} and NIF3-Ralston \eqref{def: Ralston NIF3 coefficients} schemes generate similar numerical behaviors and we omit relevant presentations.

\lan{\begin{example}\label{example 4} \cite{WangKouCai:2020,WuFengHeQian:2023} Consider the Allen-Cahn model with the Flory-Huggins potential,  $\partial_tu=\epsilon^2\partial_{xx} u-\tfrac{\theta}{2}\ln\tfrac{1-u}{1+u}-\theta_cu$ on $\Omega=(0,2\pi)$, subject to the same initial data in Example \ref{example 3}. Take the interface parameter $\epsilon=0.1$,  $\theta_c=1$  and $\theta=0.8$. Always, the second-order difference approximation with the spacing $h=\pi/320$ is used for the spatial discretization. \lan{The space-discrete solution preserves the maximum bound principle, that is, $\mynorm{u_h(t)}_{\infty}\le1$ since $\mynorm{u_h^0}_{\infty}\le1$, cf. \cite[Section 2]{DuJuLiQiao:2021SIREV}}.
\end{example}}

\begin{figure}[htb!]
  \centering
  \subfigure[Final solution $u_h^N$]{
    \includegraphics[width=2.05in,height=1.5in]{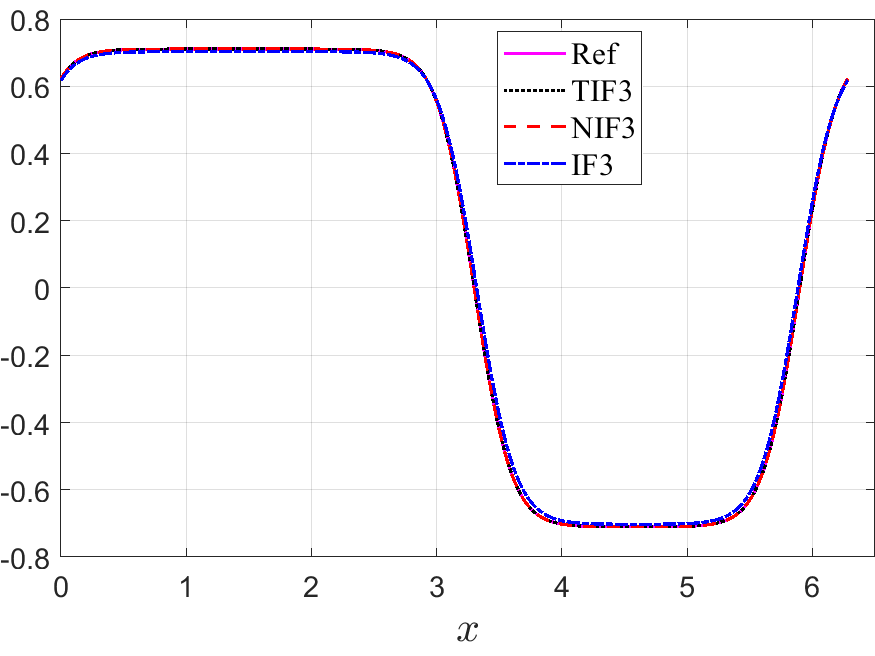}} 
  \subfigure[Discrete energy $E\kbrat{u_h^n}$]{
    \includegraphics[width=2.05in,height=1.5in]{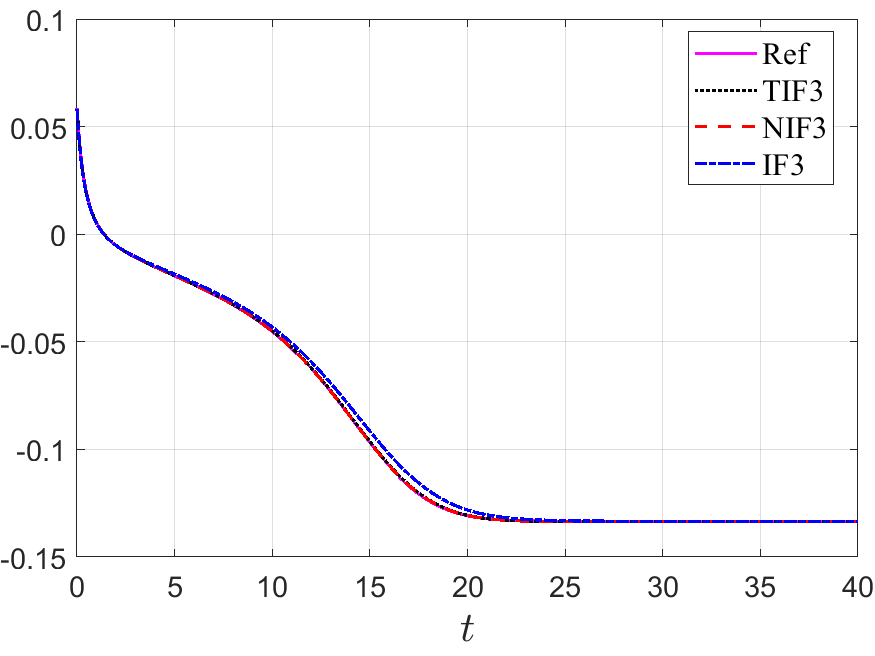}}
  \subfigure[Maximum norm $\|u_h^n\|_{\infty}$]{
    \includegraphics[width=2.05in,height=1.5in]{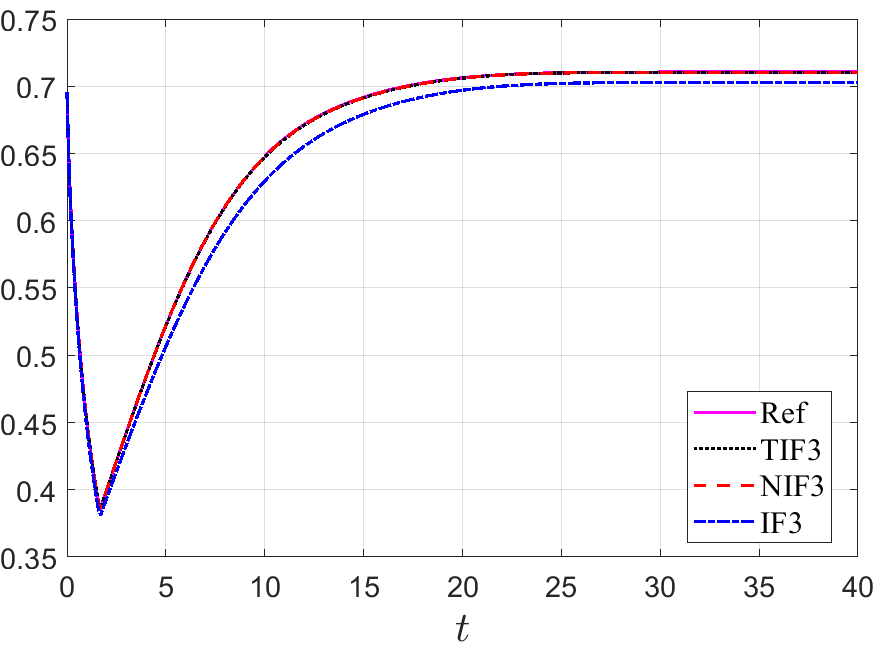}}
  \caption{Comparisons of IF3-Heun and corrected IF3-Heun schemes for $\tau=0.1$.}
  \label{fig: Compare Heun IF3 and CIF3 with tau = 0.1 FH}
\end{figure}
\begin{figure}[htb!]
  \centering
  \subfigure[Final solution $u_h^N$]{
    \includegraphics[width=2.05in,height=1.5in]{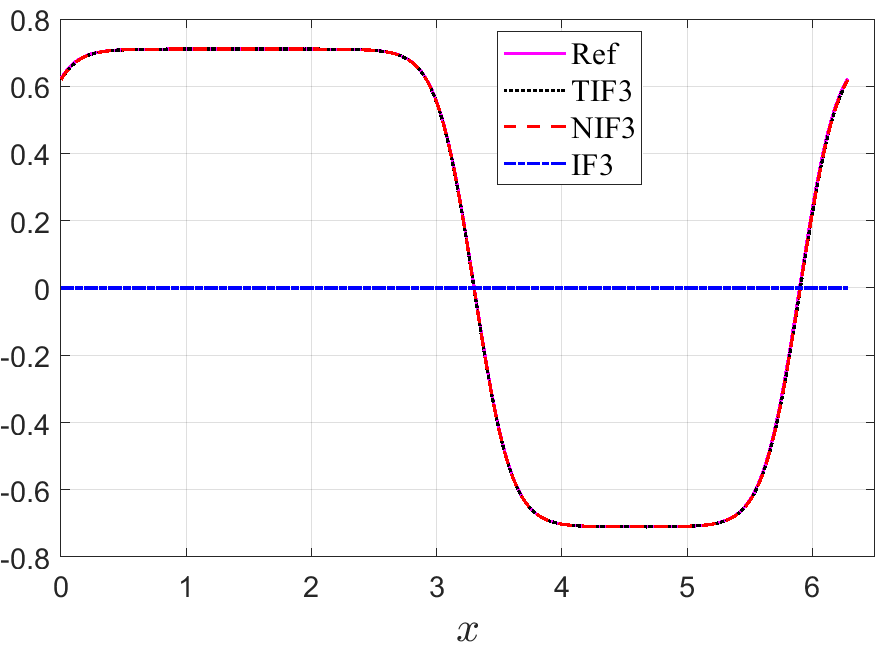}} 
  \subfigure[Discrete energy $E\kbrat{u_h^n}$]{
    \includegraphics[width=2.05in,height=1.5in]{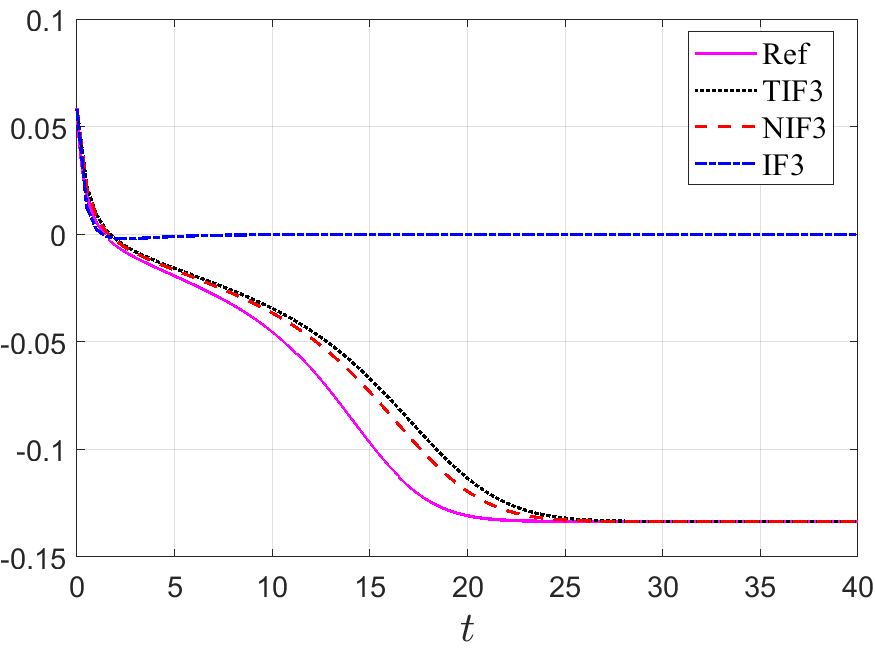}}
  \subfigure[Maximum norm $\|u_h^n\|_{\infty}$]{
    \includegraphics[width=2.05in,height=1.5in]{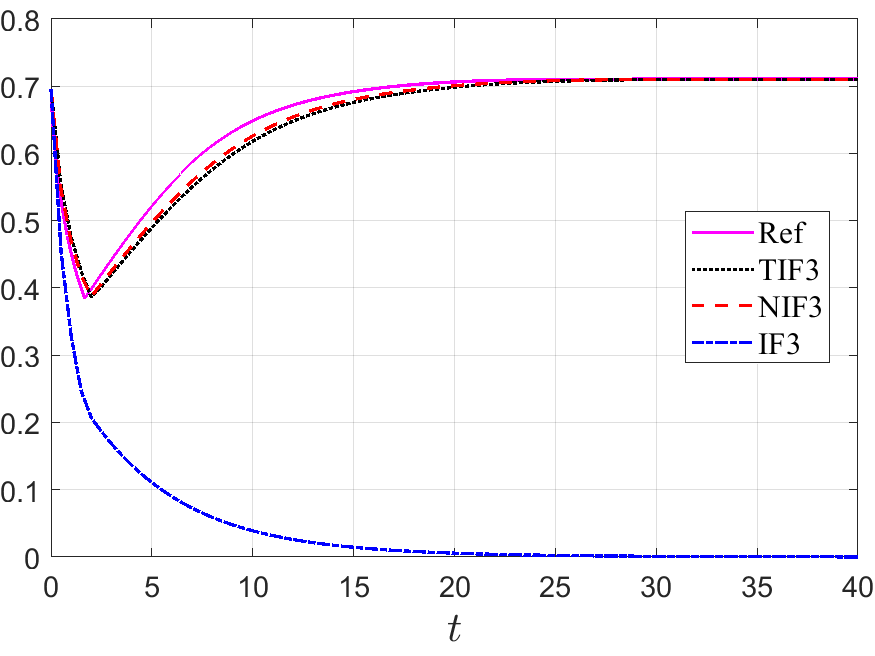}}
  \caption{Comparisons of IF3-Heun and corrected IF3-Heun schemes for $\tau=0.5$.}
  \label{fig: Compare Heun IF3 and CIF3 with tau = 0.5 FH}
\end{figure}

\lan{We run the IF3-Heun, TIF3-Heun \eqref{def: Heun TIF3 coefficients} and NIF3-Heun \eqref{def: Heun NIF3 coefficients} schemes up to the final time $T=40$ with the stabilized parameter $\kappa=4$. We depict the final solution $u_h^N$, the discrete energy $E[u_h^n]$ and the maximum norm $\|u_h^n\|_{\infty}$ in Figures \ref{fig: Compare Heun IF3 and CIF3 with tau = 0.1 FH} and \ref{fig: Compare Heun IF3 and CIF3 with tau = 0.5 FH} for two different time-steps $\tau=0.1$ and $0.5$, respectively. The numerical solution (energy) of the NIF3-Heun scheme computed by $\tau=0.01$ is taken as the reference solution (energy). As seen, the two corrected IF3-Heun methods maintain the steady-state well while the final solution $u_h^N$ of IF3-Heun scheme gradually collapses into the metastable state solution as the time-step size $\tau$ increases to $0.5$, see Figures \ref{fig: Compare Heun IF3 and CIF3 with tau = 0.1 FH}(a)-\ref{fig: Compare Heun IF3 and CIF3 with tau = 0.5 FH}(a). We notice that there are a bit differences in the energy curves of two corrected IF3-Heun schemes for $\tau=0.5$ although they are decreasing over the time and approach the same steady-state.}

\lan{Also, we note that the two corrected IF3-Heun methods preserve the maximum bound principle. This property may be closely related to the positivity of the corresponding correction coefficients, cf. Figure \ref{fig: corrected IF3-Heun coefficients rate}. It is to mention that, under the same parameter settings, the IF3-Ralston, TIF3-Ralston \eqref{def: Ralston TIF3 coefficients} and NIF3-Ralston \eqref{def: Ralston NIF3 coefficients} schemes generate similar numerical behaviors and we omit relevant presentations.}

\subsection{Simulations of bubbles merging}\label{subsec:test 2D}

\begin{example}\cite{LiaoTangZhou:2020SINUM}\label{example 2D}
  Consider the 2D Allen-Cahn model $\partial_tu=\epsilon^2\Delta u-u^3+u$ with interface parameter $\epsilon=0.05$ on $\Omega=(-1,1)^{2}$  subject to the 2-periodic initial data
  \begin{align*}
    u_0=&\, -\tanh\kbra{\bra{(x-0.3)^2+y^{2}-0.2^2}/\epsilon}\tanh\kbra{\bra{(x+0.3)^2+y^{2}-0.2^2}/\epsilon} \\
    &\,\times \tanh\kbra{\bra{x^{2}+(y-0.3)^2-0.2^2}/\epsilon}\tanh\kbra{\bra{x^{2}+(y+0.3)^2-0.2^2}/\epsilon}.
  \end{align*}
\end{example}

\begin{figure}[htb!]
     \centering
     \subfigure[$t=0$]{
          \includegraphics[width=2.0in]{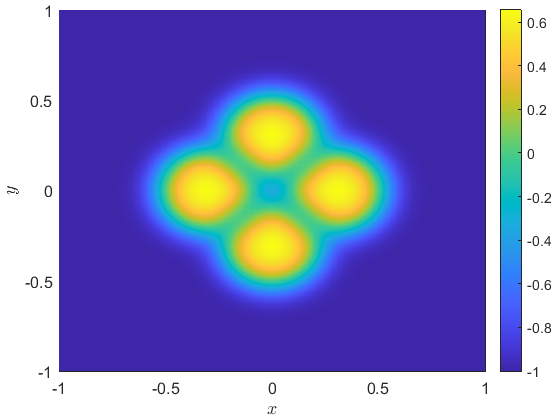}} 
     \subfigure[$t=5$]{
          \includegraphics[width=2.0in]{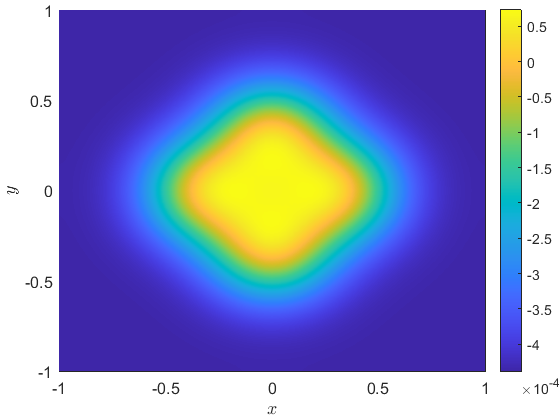}}
     \subfigure[$t=10$]{
          \includegraphics[width=2.0in]{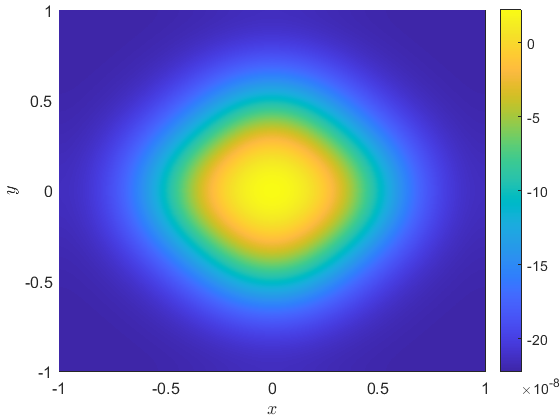}}
     \caption{Solution profiles generated by IF3-Ralston scheme for Example \ref{example 2D}.}
     \label{fig: 2D evolution process of IF3-Ralston}
\end{figure}
\begin{figure}[htb!]
     \centering
     \subfigure[$t=0$]{
          \includegraphics[width=2.0in]{Ralston3_2D_tau0.1InitialData.png}} 
     \subfigure[$t=5$]{
          \includegraphics[width=2.0in]{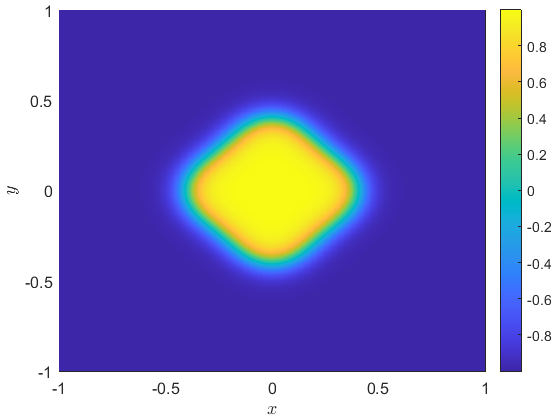}}
     \subfigure[$t=10$]{
          \includegraphics[width=2.0in]{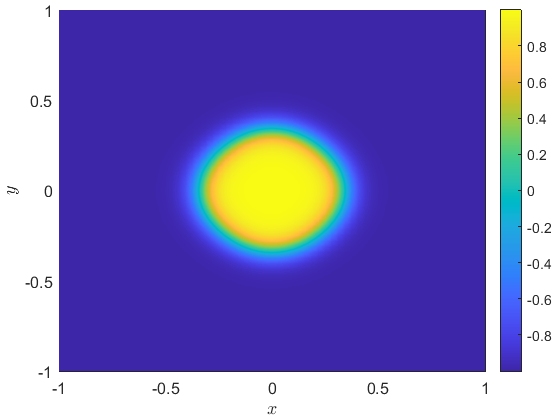}}
     \caption{Solution profiles generated by NIF3-Ralston scheme for Example \ref{example 2D}.}
     \label{fig: 2D evolution process of NIF3-Ralston}
\end{figure}

\begin{figure}[htb!]
   \centering
   \subfigure[Discrete energy $E\kbrat{u_h^n}$]{
       \includegraphics[width=2.4in]{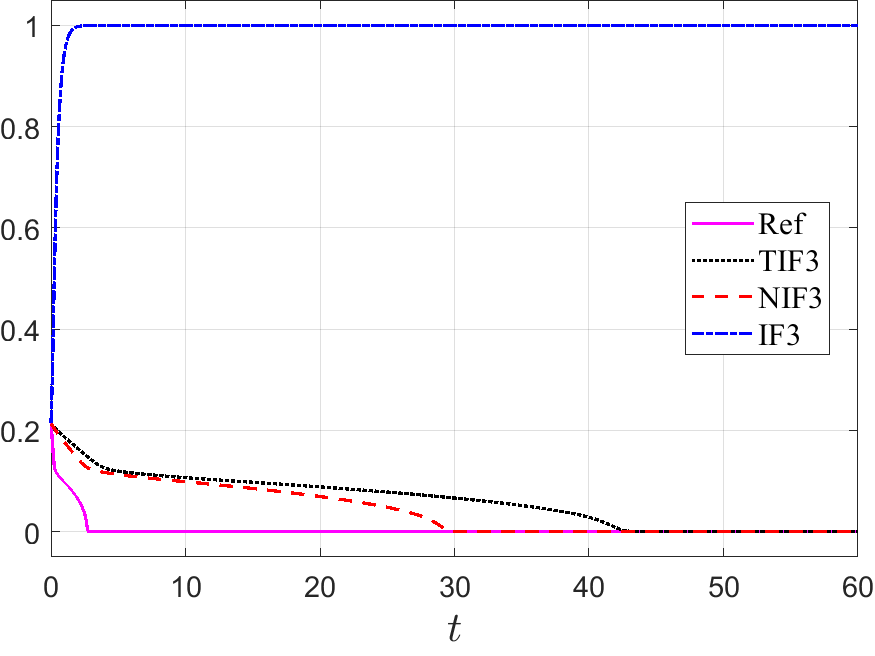}} 
 \hspace{1cm}
   \subfigure[Maximum norm $\|u_h^n\|_{\infty}$]{
       \includegraphics[width=2.4in]{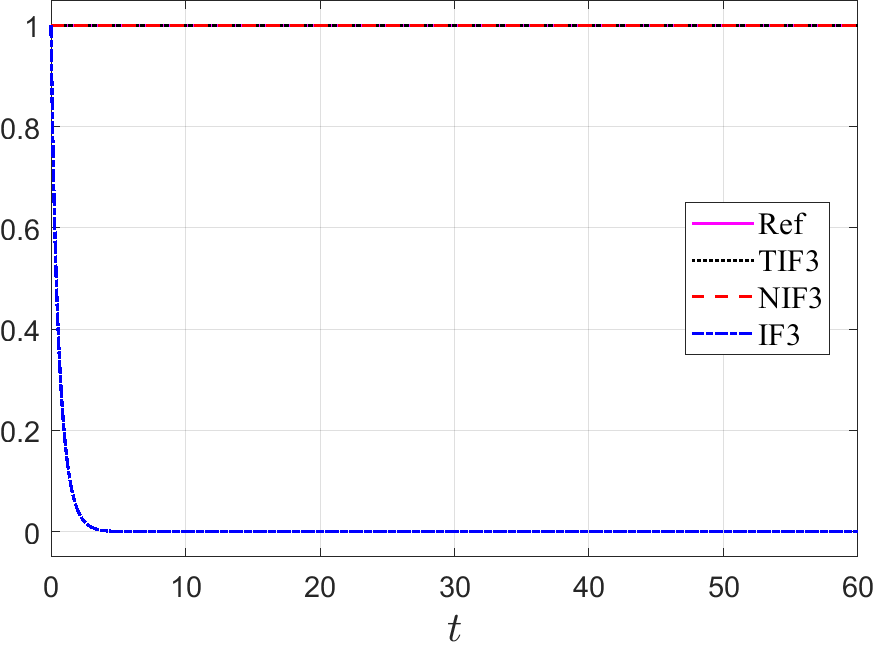}}
   \caption{Discrete energy and maximum norm of NIF3-Ralston scheme for Example \ref{example 2D}}
   \label{fig: 2D energy and max}
 \end{figure}

We now apply the IF3-Ralston scheme and two corrected versions to simulate the merging of bubbles up to $T=60$ with the spatial length $h_x=h_y=1/32$, the time-step $\tau=0.1$ and the stabilized parameter $\kappa=6$. Figures \ref{fig: 2D evolution process of IF3-Ralston} and \ref{fig: 2D evolution process of NIF3-Ralston} present the phase profiles at times $t = 0, 5$ and $10$ by using the IF3-Ralston and NIF3-Ralston scheme \eqref{def: Ralston NIF3 coefficients}, respectively. The solution profiles of the TIF3-Ralston scheme are omitted here since they are difficult to distinguish from the profiles in Figure \ref{fig: 2D evolution process of NIF3-Ralston}. The discrete curves of the associated energy $E[u_h^n]$ and maximum norm $\|u_h^n\|_{\infty}$ generated by the three numerical schemes are depicted in Figure \ref{fig: 2D energy and max}. The reference energy (maximum norm) curve is computed by the NIF3-Ralston scheme with $\tau=0.01$. As seen in Figures \ref{fig: 2D evolution process of IF3-Ralston}(b)-(c), the magnitudes of solution are about $O(10^{-4})$ and $O(10^{-8})$ at the times $t=5$ and $t=10$, respectively. That is to say, the numerical solution $u_h^n$ computed by IF3-Ralston scheme rapidly collapses to the metastable state solution, cf. the energy and the maximum norm curves in Figure \ref{fig: 2D energy and max}. As expected, the initial separated four bubbles gradually merge into a single bubble and the single bubble gradually shrinks as the time escapes, due to that the Allen-Cahn model dose not conserve the volume. Note that, they are accordant with the numerical results in previous studies \lan{\cite{ChenLiuYiYin:2024,LiLiJuFeng:2021SISC,LiaoTangZhou:2020SINUM}}.


\subsection{Tests of corrected IF4 schemes}\label{subsec:test IF4}

This subsection uses Example \ref{example 3} to examine the numerical behaviors of the IF4-Kutta method and the corresponding corrections, including the TIF4-Kutta \eqref{def: Kutta TIF4 coefficients} and NIF4-Kutta \eqref{def: Kutta NIF4 coefficients} schemes.
 \begin{figure}[htb!]
   \centering
   \subfigure[Final solution $u_h^N$]{
     \includegraphics[width=2.05in,height=1.5in]{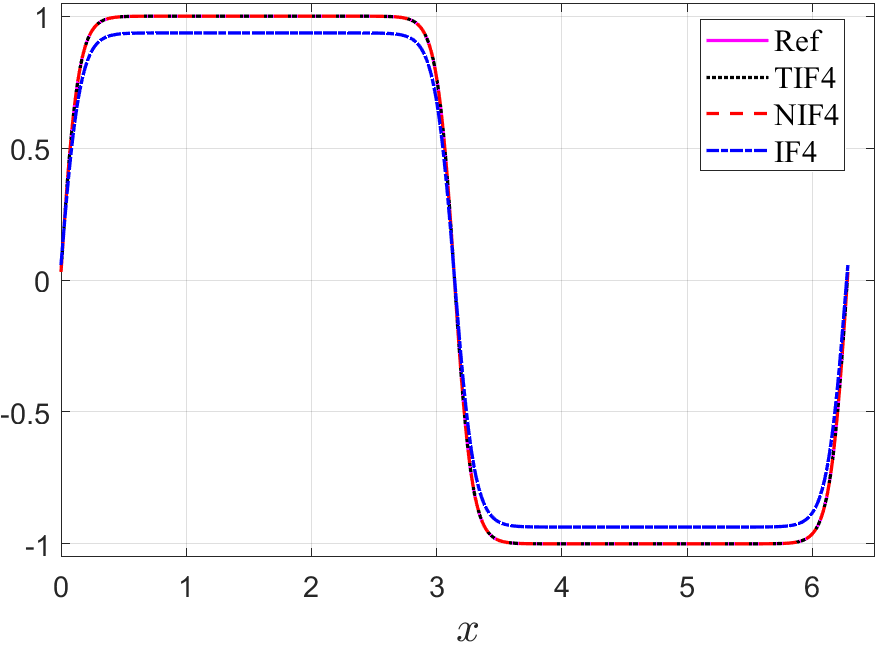}} 
   \subfigure[Discrete energy $E\kbrat{u_h^n}$]{
     \includegraphics[width=2.05in,height=1.5in]{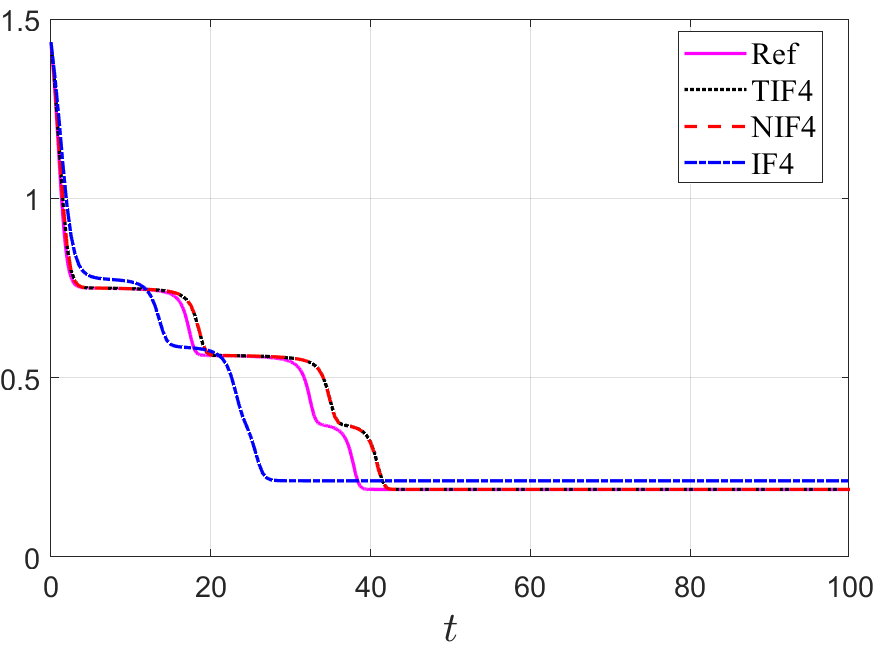}}
   \subfigure[Maximum norm $\|u_h^n\|_{\infty}$]{
     \includegraphics[width=2.05in,height=1.5in]{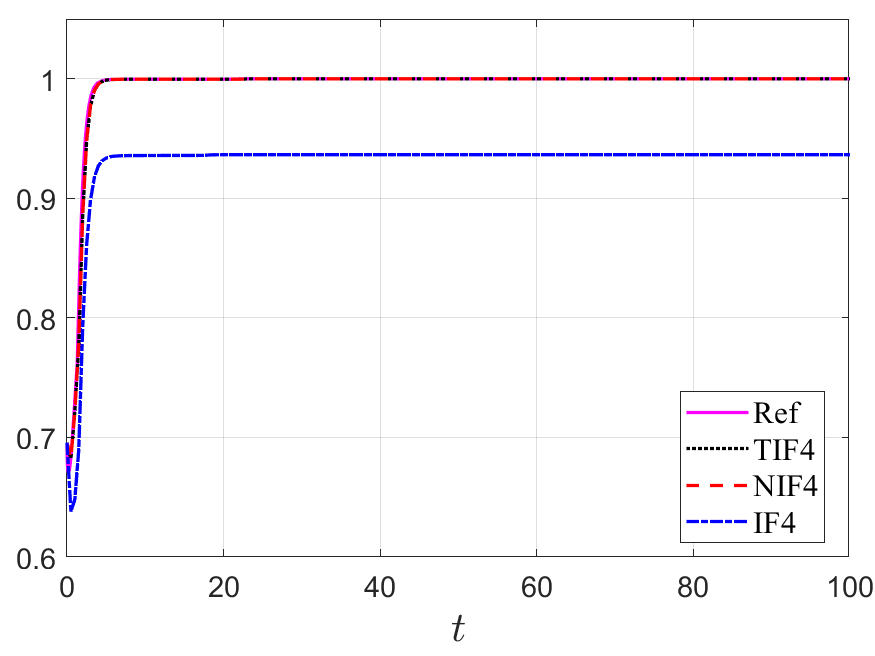}}
   \caption{Comparisons of IF4-Kutta and corrected IF4-Kutta schemes for $\tau=0.5$.}
   \label{fig: Compare IF4 and CIF4 with tau = 0.5}
 \end{figure}
 \begin{figure}[htb!]
   \centering
   \subfigure[Final solution $u_h^N$]{
     \includegraphics[width=2.05in,height=1.5in]{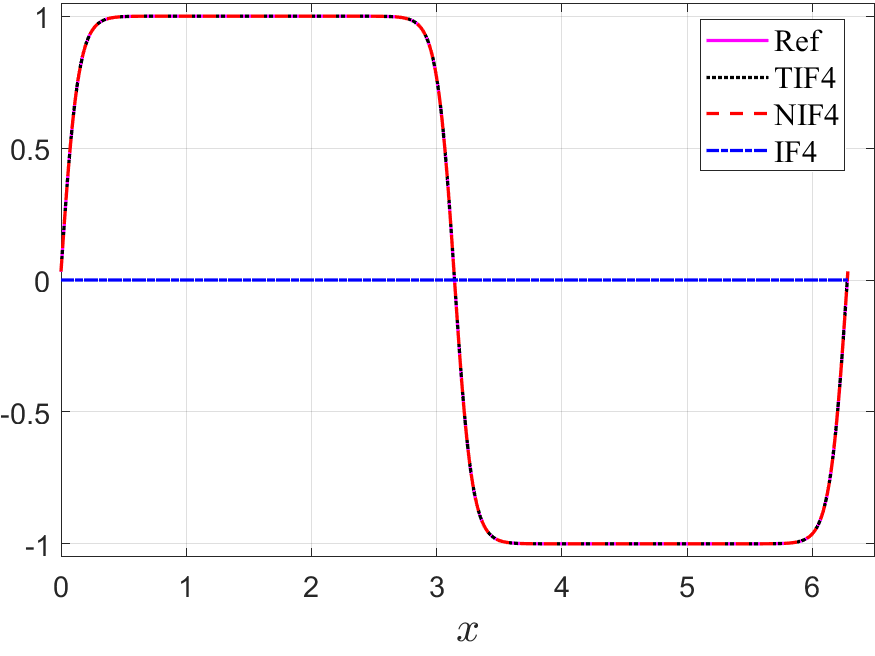}} 
   \subfigure[Discrete energy $E\kbrat{u_h^n}$]{
     \includegraphics[width=2.05in,height=1.5in]{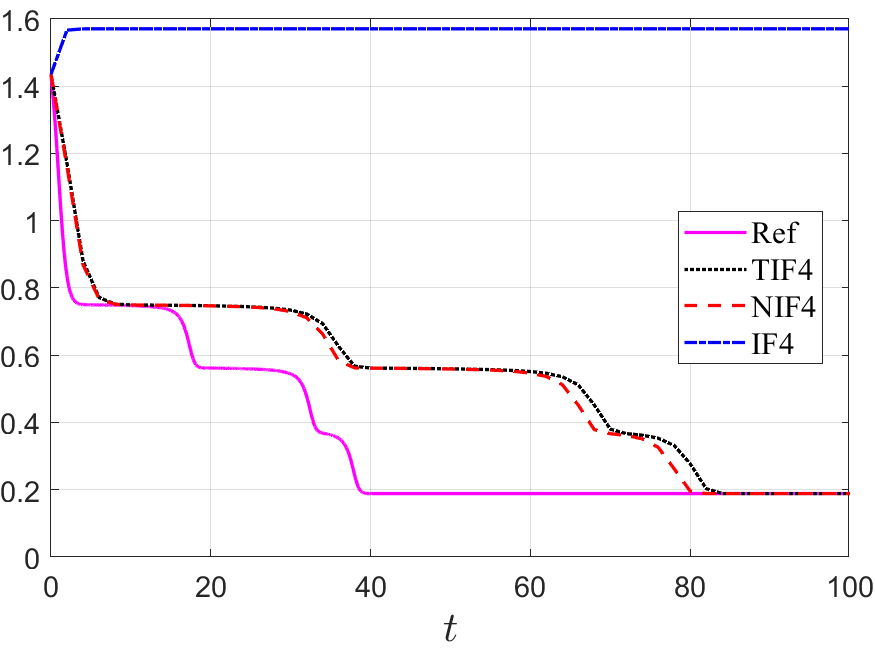}}
   \subfigure[Maximum norm $\|u_h^n\|_{\infty}$]{
     \includegraphics[width=2.05in,height=1.5in]{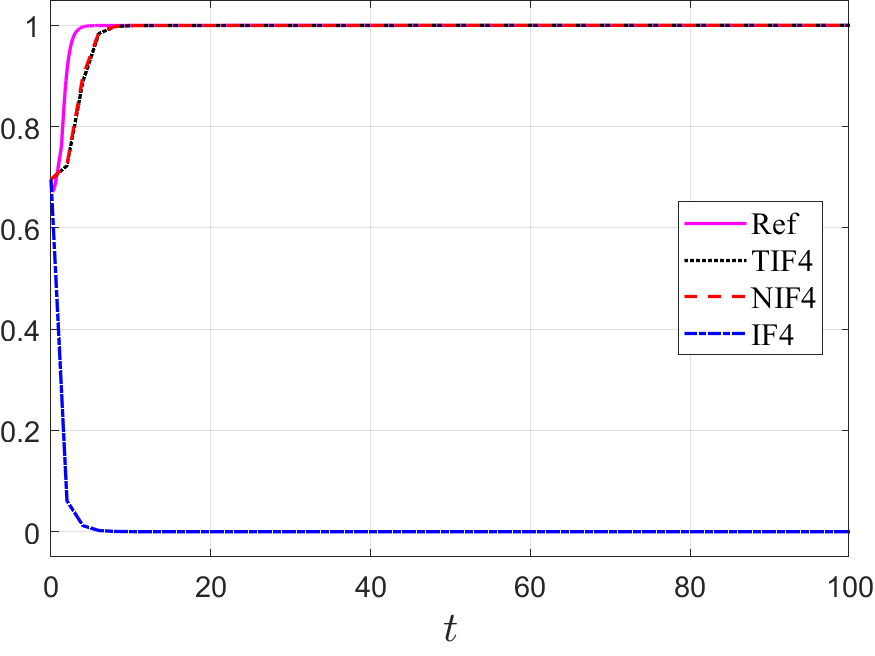}}
   \caption{Comparisons of IF4-Kutta and corrected IF4-Kutta schemes for $\tau=2$.}
   \label{fig: Compare IF4 and CIF4 with tau = 2}
 \end{figure}

We run the IF4-Kutta, TIF4-Kutta \eqref{def: Kutta TIF4 coefficients} and NIF4-Kutta \eqref{def: Kutta NIF4 coefficients} schemes up to $T=100$ with the stabilized parameter $\kappa=4$ for different time-steps $\tau=0.5$ and 2, as shown in Figures \ref{fig: Compare IF4 and CIF4 with tau = 0.5}-\ref{fig: Compare IF4 and CIF4 with tau = 2}, respectively, together with the reference solution computed by the NIF4-Kutta scheme with $\tau=0.01$. The two groups of figures list the final solution $u_h^N$, discrete energy $E[u_h^n]$, and maximum norm $\|u_h^n\|_{\infty}$ of these methods. Again, we see that the IF4-Kutta solution $u_h^N$ collapses to the trivial solution as $\tau$ increases to 2, while the solutions $u_h^N$ of the corrected IF4-Kutta schemes maintain the steady-state well. Although Section \ref{sec:CIF4} shows that the TIF4-Kutta and NIF4-Kutta schemes can not preserve the original energy dissipation law, the numerical phenomena in Figures \ref{fig: Compare IF4 and CIF4 with tau = 0.5}-\ref{fig: Compare IF4 and CIF4 with tau = 2} are not mysterious to us. Actually, the corresponding differentiation matrices $D_{\cte}^{(4)}(z)$ and $D_{\cnt}^{(4)}(z)$ seem to be positive definite for $z\le -0.5$, cf. Figure \ref{fig: corrected IF4-Kutta determinants}, although we do not verify it in mathematical manner.

\section{Concluding remarks}\label{sec:conclusion}

A simple steady-state preserving idea is proposed to overcome the main defect of IF methods in solving gradient flow problems. 
With two classes of difference correction, including the telescopic correction and nonlinear-term translation correction, 
this idea leads to many new IF methods. To distinguish the effectiveness of these corrected methods, the original energy dissipation properties of the new methods are examined 
by using the associated differentiation matrices. Some corrected IF methods up to third-order maintaining the original energy dissipation law are constructed by applying the difference correction strategies to some popular IF methods, including the widespread IF1, IF2-Heun, IF2-Ralston, IF3-Heun and IF3-Ralston schemes. We have the following results:
\begin{itemize}[itemindent=-0.5cm]
  \item[(i)] Theoretically, a total of 10 corrected IF schemes up to third-order are shown to preserve the original energy dissipation law \eqref{eq: continuous energy law} for a properly large stabilized parameter $\kappa>0$, see two first-order corrections in Theorem \ref{thm: energy stability CIF1}, four second-order corrections in Theorems \ref{thm: energy stability CIF2-Heun} and \ref{thm: energy stability CIF2-Ralston}, and four third-order corrections in Theorems \ref{thm: energy stability CIF3-Heun} and \ref{thm: energy stability CIF3-Ralston}.
  \item[(ii)] Experimentally, if the correction coefficients $\chi_{i+1}^{(s)}(z)\in(0,1)$ for $1\le i\le s$ and $z\le0$, the resulting corrected IF scheme seems to preserve the maximum bound principle in solving the Allen-Cahn model, cf. the numerical tests in subsections 2.3, 3.4 and 5.1. 
\end{itemize}
 
At the same time, our theory is far away from complete. There are many interesting issues that we have not yet addressed. Some of them are listed as follows:
\begin{itemize}[itemindent=-0.5cm]
   \item[(a)] \lan{How to prove theoretically that if the correction coefficients $\chi_{i+1}^{(s)}(z)\in(0,1)$ for $1\le i\le s$ and $z\le0$, the resulting corrected IF scheme can maintain the maximum bound principle of Allen-Cahn type model? It seems that certain constraint of the stabilized parameter $\kappa$ (or time-step size) would be required for the preservation of maximum bound principle.}
   \item[(b)] \lan{How can we choose the best scheme among numerous corrected IF schemes preserving the original energy dissipation law? In other words, some new practical or theoretical criteria would be required to distinguish various corrected IF methods for the long-time simulation of gradient flow problems. Maybe, the concept of average dissipation rate in our recent work \cite{Liao:2024arxiv} would be somewhat helpful and we will investigate it in a separate report.}
  \item[(c)] \lan{Although the energy dissipation properties of the presented corrected IF scheme are derived for Allen-Cahn type model, they are applicable to other gradient flow models, cf. \cite{FuShenYang:2024arxiv-ETDRK,Liao:2024arxiv}. In general, the best scheme for Allen-Cahn type model would be not necessarily the best scheme for other gradient flows. At least, we do not know which of the two types of algorithms, the corrected IF schemes and the EERK methods in \cite{HochbruckOstermann:2005SINUM,Liao:2024arxiv}, is better for the long-time simulations under the same accuracy.}
  \item[(d)] We are not able to find (or prove the non-existence of) a fourth-order correction preserving the energy dissipation law \eqref{eq: continuous energy law} for the IF4-Kutta method. 
\end{itemize}



\section*{Acknowledgments}
H.-L. Liao's work is supported by the National Natural Science Foundation of China (No. 12071216).

\appendix
\section{Differentiation matrices of corrected IF3-Heun schemes}\label{app:IF3-Heun}

(\textbf{$\mathrm{T}$-type correction})
Using the formula \eqref{eq: general T-type differentiation matrix}, the coefficient matrix \eqref{def: Heun TIF3 coefficients} yields the associated differentiation matrix
\begin{align*}    
  D_{\cte}^{(3,H)}(z)=&\;\begin{pmatrix}
 3-\frac{z}{2}  \\[2pt]
 \frac{3}{2} e^{\frac{z}{3}} & \frac{3}{2} e^{\frac{z}{3}}-\frac{z}{2} \\[2pt]
 \frac{1}{3} e^{\frac{2 z}{3}} & \frac{1}{3} e^{\frac{2 z}{3}} (4-z) & \frac{1}{3} e^{\frac{2 z}{3}} (4-z)-\frac{z}{2}
  \end{pmatrix}\,.
\end{align*}
The first leading principal minor of $\mathcal{S}(D_{\cte}^{(3,H)};z)$ is 
$$\mathrm{Det}[\mathcal{S}(D_{\cte,1}^{(3,H)};z)] =3-\tfrac{z}{2} \ge 3\quad\text{ for $z\le0$},$$
and the second one 
$$\mathrm{Det}[\mathcal{S}(D_{\cte,2}^{(3,H)};z)] =\tfrac{z^2}{4}-\tfrac{3z}{4} (e^{\frac{z}{3}}+2)+\tfrac{9}{16} e^{\frac{z}{3}} (8-e^{\frac{z}{3}})\ge\tfrac{63}{16}\quad\text{ for $z\le0$.}$$
Also, it is straightforward that the third leading principal minor
\begin{align*}
  \mathrm{Det}[\mathcal{S}(D_{\cte}^{(3,H)};z)] =&\; -\tfrac{z^3}{72}  (6 e^{\frac{2 z}{3}}-e^{\frac{4 z}{3}}+9)+\tfrac{z^2}{72}  (27 e^{\frac{z}{3}}+60 e^{\frac{2 z}{3}}+18 e^z-14 e^{\frac{4 z}{3}}+54) \\
  &\; -\tfrac{ z}{288} e^{\frac{z}{3}} (495 e^{\frac{z}{3}}+720 e^{\frac{2 z}{3}}-314 e^z+12 e^{\frac{4 z}{3}}+648)\\
  &\;+\tfrac{1}{24} e^z (-50 e^{\frac{z}{3}}+3 e^{\frac{2 z}{3}}+144) \ge \tfrac{47}{12} \quad\text{for $z\le0$},
\end{align*}
such that the differentiation matrix $D_{\cte}^{(3,H)}(z)$ is positive definite for $z\le0$.

(\textbf{$\mathrm{N}$-type correction})
By using \eqref{eq: general N-type differentiation matrix}, one has the associated differentiation matrix of the NIF3-Heun scheme \eqref{def: Heun NIF3 coefficients} 
\begin{align*}    
  D_{\cnt}^{(3,H)}(z)=\begin{pmatrix}
 \frac{z}{e^{\frac{z}{3}}-1}+\frac{z}{2} \\[3pt]
 \frac{z}{e^{\frac{2 z}{3}}-1}+ z & \frac{z}{e^{\frac{2 z}{3}}-1}+\frac{z}{2} \\[5pt]
 \frac{e^z [e^{\frac{z}{3}} (z-4)+4] z}{(e^{\frac{z}{3}}-1) [e^z (z-4)+4]} & \frac{e^z (z-4) z}{e^z (z-4)+4} & \frac{[e^z (z-4)-4] z}{2 e^z (z-4)+8}
  \end{pmatrix}\,.
\end{align*}
The first leading principal minor 
$$\mathrm{Det}[\mathcal{S}(D_{\cnt,1}^{(3,H)};z)] =\tfrac{z}{e^{\frac{z}{3}}-1}+\tfrac{z}{2} \ge 3\quad\text{ for $z\le0$},$$ 
and the second leading principal minor 
$$\mathrm{Det}[\mathcal{S}(D_{\cnt,2}^{(3,H)};z)] =\tfrac{z^2}{4(e^{\frac{2 z}{3}}-1)^2} (2 e^{\frac{z}{3}}+2 e^{\frac{2 z}{3}}+2 e^z+1) \ge \tfrac{9}{16}\quad\text{ for $z\le0$.} $$
The third leading principal minor of $\mathcal{S}(D_{\cnt}^{(3,H)};z)$ reads
\begin{align*}
  \mathrm{Det}[\mathcal{S}(D_{\cnt}^{(3,H)};z)] =&\; \frac{z^3 [e^z (z-4)+4]^{-2}}{8 (e^{\frac{z}{3}}-1)^3 (e^{\frac{z}{3}}+1)^2}  g_{\cnt}^{(3,H)}(z)\ge 0 \quad\text{for }z\le0,
\end{align*}
where the auxiliary function $g_{\cnt}^{(3,H)}(z)$ is defined by
\begin{align*}
  g_{\cnt}^{(3,H)}(z): =&\; e^{2 z} z^2 (-2 e^{\frac{z}{3}}+e^{\frac{2 z}{3}}-e^z-2) - 8 e^{2 z} z (-e^{\frac{z}{3}}+3 e^{\frac{2 z}{3}}-2) \\
  &\; + 16 (e^{\frac{z}{3}}-2 e^{\frac{4 z}{3}}-3 e^{2 z}-e^{\frac{7 z}{3}}+4 e^{\frac{8 z}{3}}+1) \ge 0 \quad\text{for }z\le0\,.
\end{align*}
Thus the differentiation matrix $D_{\cnt}^{(3,H)}(z)$ is positive semi-definite for $z\le0$.

\section{Differentiation matrices of corrected IF3-Ralston schemes}\label{app:IF3-Ralston}

(\textbf{$\mathrm{T}$-type correction}) By using \eqref{eq: general T-type differentiation matrix}, we have the differentiation matrix of TIF3-Ralston scheme \eqref{def: Ralston TIF3 coefficients}
\begin{align*}    
  D_{\cte}^{(3,R)}(z)=\begin{pmatrix}
 2-\frac{z}{2}  \\
 \frac{4 e^{\frac{z}{2}}}{3} & \frac{4}{3} e^{\frac{z}{2}}-\frac{z}{2}  \\[1ex]
  \frac{1}{4} e^{\frac{3 z}{4}} & \frac{1}{4} e^{\frac{3 z}{4}} (5-2 z) & -\frac{3}{4} e^{\frac{z}{4}} z-\frac{z}{2} +\frac{1}{4}e^{\frac{3 z}{4}} (9-2 z)
  \end{pmatrix}\,.
\end{align*} 
The first leading principal minor 
$$\mathrm{Det}[\mathcal{S}(D_{\cte,1}^{(3,R)};z)] =2-\tfrac{z}{2} \ge 2\quad\text{ for $z\le0$},$$ 
 the second leading principal minor 
$$\mathrm{Det}[\mathcal{S}(D_{\cte,2}^{(3,R)};z)] =-\tfrac{2z}{3} e^{\frac{z}{2}}+\tfrac{z}{4} (z-4) + \tfrac{8}{3} e^{\frac{z}{2}}-\tfrac{4 }{9}e^z \ge \tfrac{20}{9}\quad\text{ for $z\le0$},$$ 
and the third leading principal minor
\begin{align*}
  \mathrm{Det}[\mathcal{S}(D_{\cte}^{(3,R)};z)] =&\, -\tfrac{z^3}{32} \brab{6 e^{\frac{z}{4}}+4 e^{\frac{3 z}{4}}-e^{\frac{3 z}{2}}+4}  \\
  &\,+\tfrac{z^2}{96} \brab{72 e^{\frac{z}{4}}+32 e^{\frac{z}{2}}+150 e^{\frac{3 z}{4}}+32 e^{\frac{5 z}{4}}-27 e^{\frac{3 z}{2}}+48} \\
  &\, -\tfrac{z}{576} e^{\frac{z}{2}} \brab{768+2448 e^{\frac{z}{4}}-128 e^{\frac{z}{2}}+1440 e^{\frac{3 z}{4}}-477 e^{z}-128 e^{\frac{5 z}{4}}+24 e^{\frac{3 z}{2}}} \\
  &\,+\tfrac{1}{96}e^{\frac{5 z}{4}} \brab{-75 e^{\frac{z}{4}}-96 e^{\frac{z}{2}}+8 e^{\frac{3 z}{4}}+576}\ge \tfrac{413}{96}\qquad\text{for $z\le0$.}
\end{align*}
Then the differentiation matrix $D_{\cte}^{(3,R)}(z)$ is positive definite for $z\le0$.

(\textbf{$\mathrm{N}$-type correction}) By using \eqref{eq: general N-type differentiation matrix}, one has the following differentiation matrix of the NIF3-Ralston scheme \eqref{def: Ralston NIF3 coefficients}
\begin{align*}    
  \kbrab{D_{\cnt}^{(3,R)}(z)}^T=\begin{pmatrix}
 \frac{z}{e^{\frac{z}{2}}-1}+\frac{z}{2} & \frac{z}{e^{\frac{3 z}{4}}-1}+z & \frac{e^z z [5 e^{\frac{z}{2}} z+3 e^{\frac{z}{4}} (z+3)+e^{\frac{3 z}{4}} (2 z-9)+e^z (2 z-9)+9]}{(e^{\frac{z}{4}}-1) (e^{\frac{z}{4}}+1) (e^{\frac{z}{4}}+e^{\frac{z}{2}}+1) [3 e^{\frac{z}{2}} z+e^z (2 z-9)+9]} \\[5pt]
  & \frac{z}{e^{\frac{3 z}{4}}-1}+\frac{z}{2} & \frac{e^z z [3 e^{\frac{z}{4}} z-2 z+e^{\frac{3 z}{4}} (2 z-9)+9]}{(e^{\frac{z}{4}}-1) (e^{\frac{z}{4}}+e^{\frac{z}{2}}+1) [3 e^{\frac{z}{2}} z+e^z (2 z-9)+9]} \\[5pt]
  & & \frac{z}{2}-\frac{9z}{3 e^{\frac{z}{2}} z+e^z (2 z-9)+9}
  \end{pmatrix}\,.
\end{align*}
The first leading principal minor 
$$\mathrm{Det}[\mathcal{S}(D_{\cnt,1}^{(3,R)};z)] =\tfrac{z}{e^{\frac{z}{2}}-1}+\tfrac{z}{2}\ge2\quad\text{ for $z\le0$},$$
and the second leading principal minor 
$$\mathrm{Det}[\mathcal{S}(D_{\cnt,2}^{(3,R)};z)] =\tfrac{(2 e^{\frac{z}{2}}+2 e^z+1) z^2}{4 (e^{\frac{3 z}{4}}-1)^2}\ge\tfrac{4}{9}\quad\text{ for $z\le0$.} $$
The third leading principal minor
\begin{align*}
   \mathrm{Det}[\mathcal{S}(D_{\cnt}^{(3,R)};z)] =&\; \frac{z^3\kbrab{3 e^{\frac{z}{2}} z+e^z (2 z-9)+9}^{-2}}{8 (e^{\frac{z}{4}}-1)^3 (e^{\frac{z}{4}}+1) (e^{\frac{z}{4}}+e^{\frac{z}{2}}+1)^2 }g_{\cnt}^{(3,R)}(z)\ge0\quad\text{for $z\le0$,}
\end{align*}
where the auxiliary function $g_{\cnt}^{(3,R)}(z)$ is nonnegative, that is,
\begin{align*}
  g_{\cnt}^{(3,R)}(z):=&\; e^z z^2 (-21 e^{\frac{z}{2}}-20 e^z+12 e^{\frac{5 z}{4}}-8 e^{\frac{3 z}{2}}+8 e^{\frac{7 z}{4}}-4 e^{2 z}-9) \\
  &\;-18z e^{\frac{3 z}{2}} \brat{-7 e^{\frac{z}{2}}+6 e^{\frac{3 z}{4}}-2 e^z+6 e^{\frac{5 z}{4}}-3}\\
  &\;+81 (e^{\frac{z}{2}}-2 e^{\frac{3 z}{2}}-3 e^{2 z}-e^{\frac{5 z}{2}}+4 e^{\frac{11 z}{4}}+1) \ge0\quad\text{for $z\le0$.}
\end{align*}
Then the matrix $D_{\cnt}^{(3,R)}(z)$ is positive semi-definite for $z\le0$.

\end{document}